\documentclass[11pt]{amsart}
\usepackage{amssymb,amsmath,enumerate,graphics}  %maple2e,}
\date{\today}
\numberwithin{equation}{section}

\newcommand{\be}{\begin{equation}}
\newcommand{\ee}{\end{equation}}
\newcommand{\bl}{\begin{lemma}}
\newcommand{\el}{\end{lemma}}
\newtheorem{theorem}{Theorem}[section]
\newtheorem{lemma}{Lemma}[section]
\newtheorem{corollary}{Corollary}[section]

\newtheorem{conjecture}{Conjecture}[section]
\newtheorem{lemmaA}{Lemma}[section]
\newcommand{\blemmaA}{\begin{lemmaA}}
\newcommand{\elemmaA}{\end{lemmaA}}
\newtheorem{proposition}{Proposition}[section]
\newcommand{\bprop}{\begin{proposition}}
\newcommand{\eprop}{\end{proposition}}

\newcommand{\bt}{\begin{theorem}}
\newcommand{\et}{\end{theorem}}
\newcommand{\bc}{\begin{corollary}}
\newcommand{\ec}{\end{corollary}}
\newcommand{\bcon}{\begin{conjecture}}
\newcommand{\econ}{\end{conjecture}}
\newtheorem{definition}{Definition}[section]
\newcommand{\bd}{\begin{definition}}
\newcommand{\ed}{\end{definition}}
\newtheorem{remark}{Remark}[section]
\newcommand{\brem}{\begin{remark}}
\newcommand{\erem}{\end{remark}}
\newcommand{\ep}{\varepsilon}

\newcommand{\al}{\alpha}

\newcommand{\e}{\varepsilon}
\newcommand{\Om}{\Omega}

\newcommand{\area}{{\rm{area}\,}}
\newcommand{\length}{{\rm{length}\,}}
\newcommand{\closure}{{\rm{closure}\,}}
\newcommand{\CAP}{{\rm{cap}\,}}

\newcommand{\esssup}{{\rm{ess\,sup}\,}}
\newcommand{\supp}{{\rm{supp}\,}}

\begin{document}

\title[Continuous symmetrization via polarization]{ Continuous symmetrization via polarization}  % regions}

%%%%%    Information for first author
%\author[F. Brock]{Friedemann Brock}
%\address{Departmento de  Mathematicas, Facultad de Ciencias,  Universidad de Chile, Casilla 653,
%Santiago, Chile}  %
%\email{fbrock@abello.dic.uchile.cl}

\author[A. Yu. Solynin]{Alexander Yu. Solynin}
\address{Department of Mathematics and Statistics, Texas
Tech University, Box 41042, Lubbock, TX 79409, USA}  %
\email{alex.solynin@ttu.edu}  %

\begin{abstract}
We discuss  a one-parameter family of transformations which
chan\-ges sets and functions  continuously into their
$(k,n)$-Steiner symmetrizations. Our  construction consists of two
stages. First, we employ a continuous symmetrization introduced by
the author in 1990 to transform sets and functions into their
one-dimensional Steiner symmetrization. Some of our proofs in this
stage rely on a simple rearrangement called {\it polarization.}

In the second stage, we use an approximation theorem due to
Blaschke and Sarvas to give an inductive definition of the
continuous $(k,n)$-Steiner symmetrization for any $2\le k \le n$.
This transformation provides us with the desired continuous path,
along which all basic characteristics of sets and functions vary
monotonically. The latter leads to continuous versions of several
convolution type inequalities and Dirichlet's type inequalities as
well as to continuous versions of comparison theorems for
solutions of some elliptic and parabolic partial differential
equations.
 \end{abstract}

\maketitle

%shbottom
%%%%ewtheorem{corollary}{Corollary}
%\newtheorem{fig}{Box for Figure}
%\newtheorem{definition}{Definition}
%\begin{document}

\section{Introduction}\label{Introduction}
\setcounter{equation}{0}

The first geometric transformation bearing the name {\it
symmetrization} was introduced by Jacob  Steiner in 1836 \cite{St}
in one of his attempts to find a rigorous proof for the classical
isoperimetric problem. Let $C$ be a closed contour on
$\mathbb{R}^2$ enclosing a domain $D$ and let $m_D(x)$ denote the
Lebesgue measure of the intersection of $D$ with the vertical line
$v_x=\{(x,y)\in \mathbb{R}^2:\,-\infty<y<\infty\}$. Then Steiner's
symmetrization of $D$ with respect to the $x$-axis is
defined by %
\be  \label{1.1-equation} %
D^*=\{(x,y)\in \mathbb{R}^2:\,|y|<(1/2)m_D(x)\}.  %
\ee  %
This implies, in particular, that $D^*$ is symmetric with respect
to the $x$-axis and convex in the $y$-direction. Let $C^*=\partial
D^*$ be the boundary of
$D^*$. Steiner used his symmetrization to show that %
%\begin{enumerate}  %
$$ %    %
{\mbox{(a)}} \ \  \area D =\area D^*, %
\quad \quad \quad \ \quad \quad \quad \quad {\mbox{(b)}} \ \
\length C^* \le \length C, \ \ \ \ \ \ \ \ \ \ \ \ \  %
$$  %
which implies the classical isoperimetric inequality %
$$ %\be \label{1.2-equation}%
\frac{\area D}{(\length C)^2}\le \frac{1}{4\pi}, %
$$  %\ee  %
assuming the existence of a minimizer.

This ingenious idea of Steiner has appeared to be extremely
fruitful and was exploited over the years by many authors, who
proved numerous, so-called, {\it isoperimetric inequalities} for
several important geometrical and physical quantities
characterizing the shape of planar and solid regions. We want to
mention the following four such inequalities, for the transfinite
diameter $d(\bar D)$ that is  equal to the logarithmic capacity
$\CAP (\bar D)$, for the inner radius $r(D,a)$ of $D$ at its point
$a\in D$, for the torsional rigidity $P(D)$, and for the principal
frequency
$\lambda(D)$:  %
$$  %
{\mbox{(c)}} \ \ d({\overline{D^*}})\le d(\overline{D}), \qquad
\qquad {\mbox{(d)}}\ \ r(D^*,a^*)\ge r(D,a)\quad
 {\mbox{for every $a\in D$,}} %
$$  %
$$  %
{\mbox{(e)}} \ \ P(D^*)\ge P(D), \quad \quad \quad \ \
{\mbox{(f)}}\ \ \lambda(D^*)\le \lambda(D). \ \ \ \ \ \ \ \ \
\ \ \ \ \ \ \ \ \ \ \ \ \ \ \ \ \ \ \ \ \  %
$$  %

The first period of history of symmetrization was summarized in
the classical monograph of George P\'{o}lya and Gabor Szeg\"{o}
``Isoperimetric inequalities in mathematical physics'' \cite{PS}.
This fundamental study of isoperimetric inequalities was filled
with new ideas, results, and problems, some of them still remain
open. The story of continuous symmetrization also has its source
in this book. The following question was raised in Note~B of
\cite{PS}:

{\it Is it possible to define a transformation $T_\lambda$ on $D$
depending in a continuous way on the parameter $\lambda$, $0\le
\lambda \le
1$, such that the following conditions are satisfied: %
\begin{enumerate}  %
\item[{\rm{I.}}] %
$T_0$ is the identity. %
\item[{\rm{II.}}] %
$T_1$ is the transformation replacing $C$ by the symmetrized curve
$C^*$, that is, $T_1$ is Steiner's
symmetrization with respect to the line $l$. %
\item[{\rm{III.}}] %
For every $\lambda$, $0\le \lambda \le 1$, $T_\lambda$ has the same effect as described under
(a)--(d).  %
\end{enumerate}  }  %
Although  P\'{o}lya and Szeg\"{o} mentioned only relations
(a)--(d), a similar question about inequalities (e) and (f) falls
into the same context.

The authors of \cite{PS} did not explain explicitly what
motivated them to study this problem. Among obvious reasons we want to mention the following three: %
\begin{enumerate} %
\item[$\bullet$] %
First, from the point of view  of classical mechanics it is
interesting to embed $D$ and $D^*$ into a continuous path, along
which all basic geometrical and physical
characteristics of  shape vary continuously and monotonically.  %
\item[$\bullet$] %
Often better estimates, than those provided by inequalities
(b)--(f), are needed. Indeed, even for simple shapes, for example,
for the rhombus having  angle
$\alpha=\pi/100$ centered at $a=0$,  the gap in each of the inequalities (b)--(f) exceeds  40 \%.  %
\item[$\bullet$]  %
Steiner's symmetrization changes a given shape globally into a
symmetric one. %, see its effect on a domain $D$ in Figure~\ref{1.1-figure}.
So, this transformation will not work in
problems, where the minimizer does not possesses a global symmetry
and in problems concerned with  local minimality.
\end{enumerate}

\medskip

P\'{o}lya and Szeg\"{o} themselves studied this problem. In
collaboration with M.~Schiffman \cite[Note B]{PS}, they presented
one such continuous transformation and proved the relations
(a)--(d) for the case of convex domains. Suppose that $D$ is a
convex domain on $\mathbb{R}^2$  bounded from below and above by
the graphs of functions $y=y_1(x)$ and $y=y_2(x)$, respectively,
such that $y_1(x)<y_2(x)$ for all $a<x<b$. For a given continuous
function $\phi:[\alpha,\beta]\to [0,1]$,
let  %
\be \label{1.3-equation}  %
y_1^t=y_1-\phi(t)\frac{y_1+y_2}{2}, \quad \quad \quad
 y_2^t=y_2-\phi(t)\frac{y_1+y_2}{2}  %
\ee  %
and   %
\be  \label{1.4-equation}  %
D^{t,\phi}=\{(x,y)\in \mathbb{R}^2:\,a<x<b,\ y_1^t<y<y_2^t\}\,. %
\ee  %

If $\phi$ is an increasing homeomorphism from $[\alpha,\beta]$
onto $[0,1]$, then it is clear that formulas (\ref{1.3-equation}),
(\ref{1.4-equation}) define a continuous path from $D$ into $D^*$.
Choosing $\phi(t)=t$, P\'{o}lya and Szeg\"{o} \cite{PS} proved
that the relations (a)--(d) hold true for all convex domains.
S.~Abramovich \cite{Ab} used a variant of P\'{o}lya-Szeg\"{o}'s
continuous symmetrization to prove monotonicity of eigenvalues of
certain second order differential equations in one variable.

Another continuous transformation, again for smooth convex
domains, was introduced by A.~McNabb  in 1967 \cite{McN}. His
transformation, called {\it the partial
Steiner symmetrization}, can be defined as follows. We quote from \cite{McN}: %

`` A partial Steiner symmetrization of $D$ may be performed in the
following way. If the constant $\alpha$ lies between certain
limits $(\alpha_L<\alpha<\alpha_R)$, the line $x=\alpha$ will
intersect the curve defined by midpoints of the line segments
composing $D$. If just those line segments which have their
midpoints to the left of $x=\alpha$ are translated parallel to
themselves until these central points lie on $x=\alpha$, the ends
of the segments now define a curve $C_\alpha$ bounding a partially
symmetrized region $D_\alpha$. It is as though the line $x=t$
swept across the $x$-$y$-plane from $t=-\infty$ to $t=\alpha$ and
the midpoints of the line segments became attached to the line as
it passed over them. As $t$ increases from $\alpha=\alpha_L$ to
$\alpha=\alpha_R$, $D$ continuously evolves through a sequence
$D_t$ of partially symmetrized regions to its Steiner
symmetrization $D^*$. ''

As the author noted in \cite{McN}, his goal was to demonstrate on
simple examples how his transformation works. So, the treatment in
\cite{McN} was heuristic and technical ``difficulties were glossed
over'' there.

We also want to mention two interesting continuous transformations
discovered in \cite{Sz} and \cite{M}, but those are not related,
at least not directly,  to the problem raised in \cite{PS}.

\medskip

The first continuous transformation into Steiner symmetrization,
which works for non-convex domains and satisfies all the
requirements of the P\'{o}lya-Szeg\"{o} problem, was introduced by
this author  \cite{S1}. Our continuous symmetrization, which we
will abbreviate as {\it SC symmetrization}, can be considered as
an extension of  McNabb's partial Steiner symmetrization for the
case of non-convex domains.\footnote{The paper \cite{S1} does not
refer to McNabb's work \cite{McN} since at that time the author
was not aware of McNabb's publication.} We want to emphasize here
that the approaches used in \cite{McN} and \cite{S1} are
different.

It is interesting to mention that, eventually, the original idea
of P\'{o}lya and Szeg\"{o} was developed by F.~Brock \cite{B1},
\cite{B2}, who defined a continuous symmetrization, called {\it BC
symmetrization} in this paper, which works for non-convex domains.
This was achieved by choosing a parametrization
$\phi(t)=1-e^{-t}$, $-\infty<t<\infty$, in (\ref{1.3-equation}),
(\ref{1.4-equation}), combined with some other innovations.
Instead of abbreviations SC and BC, we sometimes write ``Solynin's
continuous symmetrization'' and ``Brock's continuous
symmetrization'', respectively.
%Figure~\ref{1.1-figure} shows how the SC symmetrization and BC symmetrization affect the same shape $D$.
One particular difference between SC symmetrization and BC
symmetrization is that under the first transformation the change
of the shape is localized near some boundary arcs while the second
transformation changes the boundary globally.  %

%\begin{figure}  \label{1.1-figure}  %
%{} \caption{SC and BC symmetrizations.}
%\end{figure}  %
\smallskip

Although the exposition in \cite{S1} was given for planar domains,
in the final Remark~5  \cite{S1}, the author emphasized that all
definitions and proofs can be extended without substantial changes
to $n$-dimensional spaces and that all major results of the paper
have $n$-dimensional counterparts.

\smallskip

The primary  goal of the present paper is to give a full scale
account of Solynin's continuous symmetrization in the
$n$-dimensional setting. Since the paper \cite{S1} is practically
unknown to the experts (its English translation is often
inadequate, actually it looks like a computer translation), we
want to mention here the major innovations introduced in
\cite{S1}. First of all, the {\it polarization} was used for the
first time in \cite{S1} in the context of continuous
symmetrizations. Then, an analog of the semigroup property was
applied to prove some results about the continuous
SC~symmetrization. Later on, F.~Brock \cite{B2} used a similar
property as a part of the definition of his continuous
symmetrization. Uniqueness results were treated in \cite{S1} in
all their generality. The latter leads, under certain conditions,
to strict monotonicity of the domain characteristics under
consideration as the functions of the parameter of symmetrization.
Finally, the  SC~symmetrization was applied in \cite{S1} to prove
local symmetry in some problems on Green's functions and harmonic
measures. A similar approach to local symmetry in a more general
context was also used in the papers \cite{B1} and \cite{B2}.

\smallskip

This paper is organized as follows. Section~\ref{Preliminaries}
contains our basic notations. In particular, we introduce there
necessary spaces of functions and classes of domains. In
Section~\ref{Steiner symmetrization and polarization},  we remind
the reader of  basic properties of the Steiner
$(k,n)$-symmetrization and polarization. The exposition in
Sections~\ref{Preliminaries} and \ref{Steiner symmetrization and
polarization} follows the lines of our paper \cite{BS} joint with
F.~Brock. Sections~\ref{ Continuous $(1,n)$-Steiner
symmetrization}--\ref{Rescaling and limit cases} are devoted to
geometric aspects of  SC $1$-symmetrization.

In Section~\ref{ Continuous $(k,n)$-Steiner symmetrizations}, we
give an inductive definition of the continuous $(k,n)$-Steiner
symmetrization for any $2\le k\le n$.

Sections~\ref{Integral inequalities} and \ref{Comparison
theorems}, where we again follow the lines of the paper \cite{BS},
contain our main applications. In Section~\ref{Integral
inequalities}, we show that many integral inequalities known for
the Steiner symmetrization have their continuous counterparts for
the continuous $(k,n)$-symmetrization as well. In
Section~\ref{Comparison theorems}, we give a similar treatment of
the comparison theorems for solutions of some elliptic and
parabolic PDE's. Many proofs in Sections~10 and 11 related to the
$L^p$-classes and Sobolev classes are based on ideas suggested by
F.~Brock, when we worked on Sections~9 and 10 of our joint paper
\cite{BS}, and which he developed further in \cite{B2}.

\smallskip

In the present paper we combined and extended the ideas and
methods developed in \cite{S1}, \cite{S2}, and \cite{BS}.
Preparing this article for publication, the author  used his notes
written in the Fall semester, 1995 during his stay at the
Mathematisches Forschungsinstitut Oberwolfach under the financial
support of Volkswagen-Stiftung, RiP-program for Friedemann Brock
and Alexander Solynin. Our intention at that time was to present
in a joint paper our results for both SC $k$-dimensional
continuous symmetrization and BC $k$-dimensional continuous
symmetrization. Since 1995, Brock's continuous symmetrization and
its applications were already discussed in several publications.
So, in this paper,
 we are concentrating on the Solynin's continuous
symmetrization only. Although the original plan for this paper was
changed, this work remains closely related to the paper \cite{BS}
joint with F. Brock, where such a possible continuation was
referenced as ``An approach to continuous symmetrization via
polarization''.

\section{Preliminaries }  \label{Preliminaries} %
\setcounter{equation}{0}  %

The following notations will be used  throughout the paper. Let
$\mathbb{R}^n$ be the Euclidean space,
$\mathbb{R}_+^n=\{(x_1,\ldots,x_n)\in \mathbb{R}^n:\, x_i>0,\ 1\le
i\le n\}$. %[0,\infty)$ and $\mathbb{R}^+=(0,\infty)$.

For $A\subset \mathbb{R}^n$, let $\overline{A}$ and $\partial A$
denote the closure and the boundary of $A$, respectively. If
$A,B\subset \mathbb{R}^n$  then $A+B:=\{z:\,z=x+y,\ x\in A,\ y\in
B\}$ denotes the Minkowski sum of $A$ and $B$.  For $x,y\in
\mathbb{R}^n$ , by $|x|$ and  $\langle x,y\rangle $ we denote the
norm of $x$ and the scalar product of $x$ and $y$, respectively.
Then $H(a,n)$ and $\Sigma(a,n)$ will denote the half-space $\{x\in
\mathbb{R}^n:\,\langle (x-a),n\rangle>0\}$ and the hyperplane
$\{x\in \mathbb{R}^n:\,\langle (x-a),n\rangle=0\}$ defined by the
point $a\in \mathbb{R}^n$ and the unit vector $n\in \mathbb{R}^n$.
For $M\subset \mathbb{R}^n$, by ${\mathcal{L}}^n(M)$ we denote the
$n$-dimensional Lebesgue measure of $M$. By ${\mathcal{M}}_n$,
${\mathcal{F}}_n$, and ${\mathcal{G}}_n$ we denote  the sets of
all measurable, compact, and open subsets of $\mathbb{R}^n$,
respectively. Then ${\mathcal{M}}_{n,b}$ and ${\mathcal{G}}_{n,b}$
will denote collections of all bounded subsets of
${\mathcal{M}}_n$ and
${\mathcal{G}}_n$.  %

Generally, we treat measurable sets only in an a.e. sense, i.e. we
write %
\begin{eqnarray*}  %
M=N &{\mbox{\ if and only if\ }}& {\mathcal{L}}^n(M\bigtriangleup
N)=0, \\
M\subset N &{\mbox{\ if and only if\ }}&
{\mathcal{L}}^n(M\setminus
N)=0.  %
\end{eqnarray*}  %

By $B_r^{(n)}(x_0)$ we denote the open ball in $\mathbb{R}^n$ with
radius $r>0$ centered at $x_0$ and we write
$B^{(n)}_r=B_r^{(n)}(0)$, $B^{(n)}=B_r1^{(n)}$. If $A\subset
\mathbb{R}^n$ and $\e>0$, then we denote by $A_\e:=A+\e
{\overline{B^{(n)}}}$ the exterior parallel set of $A$. The
Hausdorff
distance between compact sets $A$ and $B$ is defined by  %
$$  %
d(A,B):=\inf \{\e>0:\, A\subset B_\e,\ B\subset A_\e\}. %
$$  %
It is well known that $d$ is a metric on ${\mathcal{F}}_n$. We
define, by the metric $d$, the convergence of a sequence of sets
$F_i\in {\mathcal{F}}_n$, $i=1,2,\ldots$, to a set $F\in
{\mathcal{F}}_n$ by %
$$     %
\lim_{i\to \infty} F_i=F \quad \quad {\mbox{if and only if
\ $d(F_i,F)\to 0$ \ as $i\to \infty$.}} %
$$  %

If $\Om$ is an open set in $\mathbb{R}^n$ and $p\in[1,\infty]$
then $\|\cdot\|_p$ denotes the usual norm in the space $L^p(\Om)$.
For functions $u\in C(\mathbb{R}^n)$ we define the modulus of
continuity by %
$$  %
\omega_u(\delta):=\sup \{|u(x)-u(y)|:\,|x-y|<\delta\}, \quad
\delta>0. %
$$  %
 By $W^{1,p}(\Om)$ we denote the Sobolev space of functions $u\in
L^p(\Om)$ having generalized partial derivatives $u_{x_i}\in
L^p(\Om)$, $i=1,\ldots,n$, and we write %
\be  \label{2.2-equation}  %
\|u\|_{W^{1,p}(\Om)}:=\|u\|_p+\sum_{i=1}^n \|u_{x_i}\|_p  %
\ee  %
for the norm in this space. By $W_0^{1,p}(\Om)$ we denote the
completion of the set of infinitely differentiable functions with
compact support in $\Omega$, denoted by $C_0^\infty(\Om)$, under
the norm (\ref{2.2-equation}). Usually we extend measurable
functions $u:\,\Om\to \mathbb{R}_0^+$ by zero outside $\Om$ so
that
$W_0^{1,p}(\Om)\subset W^{1,p}(\mathbb{R}^n)$ in that sense. %
By $C_0^{0,1}(\Om)$ we denote the space of Lipschitz functions
with compact support in $\Om$. For any function space the lower
subscript ``$+$'' denotes the corresponding subspace of
nonnegative functions, e.g. $L^p_+(\Om)$, $W^{1,p}_{0+}(\Om)$,
etc.

Let ${\mathcal{S}}_n$ denote the class of real measurable
functions $u$ satisfying %
$$  %
{\mathcal{L}}^n(\{u >c\})<\infty \quad
{\mbox{for all $c>\inf u$.}}  %
$$  %
Here and in the following we use the following abbreviation:
$\{u>c\}=\{x\in \mathbb{R}^n:\,u(x)>c\}$.  Note that the spaces
$L_+^p(\mathbb{R}^n)$, $C_{0+}^{0,1}(\mathbb{R}^n)$, and the space
$W_+^{1,p}(\mathbb{R}^n)$ with $1\le p<\infty$ are subspaces of
${\mathcal{S}}_{n+}$.
 The space of measurable functions with
bounded variation is
denoted by $BV(\mathbb{R}^n)$ and we write  %
$$  %
\|Du\|_{BV}:=\sup \left\{\int_{\mathbb{R}^n} u\sum_{i=1}^n
\frac{\partial \psi_i}{\partial x_i}\,dx\,:\, \sum_{i=1}^n
\psi^2_i \le 1,\ \psi_i \in C_0^\infty(\mathbb{R}^n),\
i=1,\ldots,n\right\}\,.  %
$$  %
Recall also that if $u\in W^{1,1}(\mathbb{R}^n)$, then
$\|Du\|_{BV}=\|\nabla u\|_1$. Furthermore, if $M$ is a Caccioppoli
set in $\mathbb{R}^n$, then $\|D\chi(M)\|_{BV}$ is the perimeter
of $M$ in the sense of De Giorgi, see \cite{T}.

Finally, a function $j:\,\mathbb{R}_0^+\to \mathbb{R}_0^+$ is
called a Young function if  $j$ is continuous and convex with
$j(0)=0$.

\section{ Steiner symmetrization and polarization}  \label{Steiner symmetrization and
polarization} %
\setcounter{equation}{0}  %

First we discuss some general properties of rearrangements. We
remind the reader that a set transformation $T$ (defined on
${\mathcal{M}}_n$) is called {\it a rearrangement} if it is
monotone and measure preserving, i.e. if $T(A)\subset T(B)$ for
all $A$ and $B$ such that $A\subset B$ and
${\mathcal{L}}^n(T(A))={\mathcal{L}}^n(A)$ for every measurable
set $A$.

The class ${\mathcal{S}}_n$ introduced in the previous section is
the natural class of functions for which a rearrangement can be
defined. If $T$ is a rearrangement and $u$ is %a continuous function
in ${\mathcal{S}}_n$, then the relations %
\be  \label{3.1-equation}  %
Tu(x):=\esssup\left\{c>\inf u:\, x\in T(\{u>c\})\right\} \quad
\inf
Tu:=\inf u, %
\ee  %
define a function $Tu$ on $\mathbb{R}^n$. If $u\in{\mathcal{S}}_n$
 is   continuous, then ``$\esssup$''   in (\ref{3.1-equation}) can be replaced by ``$\sup$''.
 Clearly the function $Tu$ is uniquely
determined almost
everywhere. Since $T$ is measure preserving, %
$$  %
{\mathcal{L}}^n(T(\{u>c\}))={\mathcal{L}}^n(\{Tu>c\}) \quad \quad
{\mbox{for all $c>\inf u$. }}  %
$$  %
Thus $Tu\in {\mathcal{S}}_n$ if $u\in {\mathcal{S}}_n$. The
mapping $T:\,{\mathcal{S}}_n\to {\mathcal{S}}_n$ constructed in
this way is again called {\it a rearrangement}. The following
non-expansivity lemma will be very useful in
Sections~\ref{Integral inequalities}
and \ref{Comparison theorems}, see \cite[Theorem~3.1]{BS}. %
\bl \label{ Non-expansivity of rearrangements}  %
Let $T$ be a rearrangement. Then for every Young function $j$, we
have %
\be \label{3.2-equation}  %
\int_{\mathbb{R}^n} j\left(|Tu-Tv|\right)\,dx \le
 \int_{\mathbb{R}^n} j\left(|u-v|\right)\,dx \quad \quad {\mbox{for
all $u,v\in {\mathcal{S}}_n$,}} %
\ee  %
whenever either one of the integrals in (\ref{3.2-equation})
converges.  %
\el  %

Sometimes we will say that two functions $u,v\in{\mathcal{S}}_n$
are {\it rearrangements of each other} if $\inf u=\inf v$ and
${\mathcal{L}}^n(\{u>c\})={\mathcal{L}}^n(\{v>c\})$ for all
$c>\inf u$.

We will also use  some additional properties of rearrangements. A
set transformation $T$ is called {\it open or compact} if $T(A)$
is open or compact whenever $A$ is of the same kind, respectively.
We say that $T$ is {\it continuous from the inside} if $\cup_i
T(G_i)=T(\cup_i G_i)$ for every increasing sequence
$\{G_i\}\subset {\mathcal{G}}$. Similarly we say that $T$ is {\it
continuous from the outside} if $\cap_i T(F_i)=T(\cap_i F_i)$ for
every decreasing sequence $\{F_i\}\subset {\mathcal{F}}$.

Finally, a rearrangement $T$ is called {\it smoothing} if
$T(F_r)\supset (T(F))_r$ for every $F\in {\mathcal{F}}$ and $r>0$.
Smoothing rearrangements were introduced by Sarvas \cite{Sa}.

\bigskip

Let us now recall the definitions of the $(k,n)$-Steiner
symmetrizations (for further information see \cite{St}, \cite{L},
and \cite{Sa}).

\bd \label{Definition of $(k,n)$-Steiner symmetrization}  %
Every
$(n-k)$-dimensional plane $\Sigma\subset \mathbb{R}^n$ with $1\le
k\le n$ defines a {\it $(k,n)$-Steiner symmetrization $S$}
as follows: \\
For every $x\in \Sigma$ let $\Lambda(x)$ denote the
$k$-dimensional plane through $x$ and orthogonal to $\Sigma$.

{\bf {1)}} Let $M\in ({\mathcal{F}}_n\cup {\mathcal{G}}_n)\cap
{\mathcal{M}}_n$. If ${\mathcal{L}}^k(M\cap \Lambda(x))=0$, then
$S(M)\cap \Lambda(x)$ is empty or the point $\{x\}$ according to
whether $M\cap \Lambda(x)$ is empty or nonempty. If
${\mathcal{L}}^k(M\cap \Lambda(x))>0$, then  %
\be \label{3.3-equation} %
S(M)\cap \Lambda(x)=\left\{ %
\begin{array}{ll}  %
 B_r(x)\cap \Lambda(x)
& {\mbox{if $M$ is open,}} \\ {\overline{B_r(x)}}\cap \Lambda(x) &
{\mbox{if $M$ is compact,}}  %
\end{array}  %
\right.  %
\ee  %
where $r>0$ is defined by the condition
${\mathcal{L}}^k(B_r(x)\cap \Lambda(x))={\mathcal{L}}^k(M\cap
\Lambda(x))$.

{\bf{2)}} Let $M\in {\mathcal{M}}_n$ where $M$ is {\it neither
open nor compact.} Then the sets $S(M)\cap \Lambda(x)$ are defined
in an a.e.
sense by either one of the equations in (\ref{3.3-equation}).  %
\ed  %

\medskip

From Definition~\ref{Definition of $(k,n)$-Steiner symmetrization}
one deduces immediately that the $(k,n)$-Steiner symmetrization is
a rearrangement which is continuous from the inside and from the
outside. Note also that in case {\bf 2)} Fubini's Theorem implies
that the sets $M\cap \Lambda(x)$ are measurable with finite
${\mathcal{L}}^k$-measure for a.e. $x\in \Sigma$.

The $(n,n)$-Steiner symmetrization is often called the {\it
Schwarz symmetrization} or the {\it symmetric decreasing
rearrangement}, and we will denote it by $S^\star$.

For our purposes it will often be helpful to use a special
coordinate system in $\mathbb{R}^n=\mathbb{R}^m\times \mathbb{R}^k$, where $1\le k\le n$, $m=n-k$ and   %
$$  %
x=(x_1,\ldots,x_n)=(x',y), \quad x'=(x_1,\ldots,x_{n-k}), \quad
y=(x_{n-k+1},\ldots,x_n),  %
$$  %
in which the plane $\Sigma$ of symmetry becomes simply $\{y=0\}$.
If $M\in {\mathcal{M}}_n$, we introduce the ``$k$-slices'' of $M$
at $x'$ by %
$$  %
M(x')=\{y\in \mathbb{R}^k:\,(x',y)\in M\}, \quad \quad x'\in
\mathbb{R}^{n-k}. %
$$ %
For instance, if $x'\in\mathbb{R}^{n-k}$ with $1\le k<n$, then
$B_r^{(n)}(x')$ will denote the $k$-slice of the ball $B_r^{(n)}$
at $x'$ and not the ball in $\mathbb{R}^n$ centered at $x'$. Let
$S^\star (M(x'))$ denote the Schwarz symmetrization of $M(x')$,
{\it taken in $\mathbb{R}^k$.} Then (\ref{3.3-equation})
reads %
\be  \label{3.4-equation} %
S(M):=\{x=(x',y):\,y\in S^\star (M(x')),\ x'\in
\mathbb{R}^{n-k}\}.  %
\ee  %
If $u\in {\mathcal{S}}_n$, then we obtain from
(\ref{3.3-equation}) that the $(k,n)$-Steiner symmetrization
${\mathcal{S}}u$ of $u$ is
given by the relations  %
\be  \label{3.5-equation}  %
{\mathcal{S}}u(x',y)=\sup \{c>\inf u:\,x\in
S(\{u(x',\cdot)>c\})\}. %
\ee  %
(Here and in the following for simplicity $\{u(x',\cdot)>c\}$
denotes $\{y\in \mathbb{R}^k:\,u(x',y)>c\}$.)

Let us mention again that the equations (\ref{3.4-equation}) and
(\ref{3.5-equation}) have to be understood in the pointwise sense
if and only if $u$ is continuous. Note also that ${\mathcal{S}}u$
is ``radially symmetric and decreasing in $|y|$'', i.e.  %
\be \label{3.6-equation} %
{\mathcal{S}}u(x',y)={\mathcal{S}}u(x',z_1)\ge
{\mathcal{S}}u(x',z_2) \quad {\mbox{if $|y|=|z_1|\le |z_2|$,}}  %
\ee  %
where $x'\in \mathbb{R}^{n-k}$ and  $y$,$z_1$,$z_2\in
\mathbb{R}^k$.

Sometimes we will write $S(M)=M^*$ and ${\mathcal{S}}u=u^*$ for
the symmetrized objects.

\smallskip

There is an approach due to Schwartz and Blaschke (see, for
instance,  \cite{Bl}) reducing a $k$-dimensional symmetrization to
$(k-1)$-dimensional symmetrizations, see Theorem~4.32 in
\cite{Sa}. We will use a slightly refined version of this theorem.

Let $S$ be a $(k,n)$-Steiner symmetrization in
$\mathbb{R}^n=\mathbb{R}^m\times \mathbb{R}^k$ with the symmetry
plane $\Sigma=\{y=0\}$. Let ${\overrightarrow{v}}_1$ and
${\overrightarrow{v}}_2$ be unit vectors orthogonal to $\Sigma$
which form an angle $\gamma \pi$, where $\gamma\in(0,1)$ is
irrational. Let $S_i$ be the $(k-1,n)$-Steiner symmetrization with
the symmetry plane $\Sigma_i$ defined by the plane $\Sigma$ and
the unit vector ${\overrightarrow{v}}_i$, $i=1,2$.

For positive integer $j$ and $\Om\in {\mathcal{F}}_n\cup
{\mathcal{G}}_{n,b}$, let %
\be  \label{3.13-equation} %
\Om_j=(S_2\circ S_1)^m(\Om) \quad {\mbox{if $j=2m$ is even,}} %
\ee  %
\be \label{3.14-equation} %
\Om_j=S_1\circ (S_2\circ S_1)^m(\Om) \quad {\mbox{if $j=2m+1$ is
odd.}}  %
\ee  %
Here $(S_2\circ S_1)^0$ is the identity transformation.

\bt \label{Sarvas approximation lemma} %
Let $S$, $S_1$, and $S_2$ be the symmetrizations defined above and
let
$\Om^*=S(\Om)$. Then %
\be  \label{3.15-equation}  %
\lim\limits_{j\to \infty} d(\Om_j,\Om^*)=0 \quad \quad {\mbox{for
every compact set $\Om\in {\mathcal{F}}_n$,}}%
\ee  %

\be  \label{3.16-equation}  %
\lim\limits_{j\to \infty} d(\partial\Om_j,\partial \Om^*)=0 \quad
\quad {\mbox{for
every bounded open set $\Om\in {\mathcal{G}}_{n,b}$,}}%
\ee  %
and %
\be  \label{3.17-equation}  %
\lim\limits_{j\to \infty} {\mathcal{L}}^k((\Om_j(x')\bigtriangleup
\Om^*(x'))=0 %
\ee  %
for every  $\Om\in {\mathcal{F}}_n \cup {\mathcal{G}}_{n,b}$ and every $x'\in \mathbb{R}^{n-k}$.%
\et %

For compact sets this theorem is a part of Theorem~4.32 in
\cite{Sa}. For bounded open sets the proof will be given in the
Appendix. In Section~\ref{ Continuous $(k,n)$-Steiner
symmetrizations}, we will use the approximation scheme of
Theorem~\ref{Sarvas approximation
lemma} to give an inductive definition of our continuous $(k,n)$-Steiner symmetrization. %
%Figure~\ref{2-figure} shows how this approximation scheme effects some simple shapes.

%\begin{figure}
%$$\includegraphics[scale=.3,angle=0]{iceberg} $$
%\caption{$(2,2)$-Steiner symmetrization via $(1,2)$-Steiner
%symmetrizations.} \label{Blaschke Scheme}
%\end{figure}

\medskip

During the last decade we have seen increase of activity in the
theory of symmetrization, partly triggered by the paper \cite{BS},
that is related to the polarization. This simplest rearrangement
was introduced for sets by V.~Wolontis \cite{W} in 1952 who
attributed some ideas of his paper to L.~Ahlfors. Ahlfors himself
used polarization in \cite{A}, where he introduced this
transformation for functions. The term {\it polarization} was
suggested by V.~N.~Dubinin \cite{D1}.

Let $\Sigma$ be a hyperplane in $\mathbb{R}^n$ and let $H$ be one
of the open halfspaces into which $\mathbb{R}^n$ is divided by
$\Sigma$. Let $\sigma_H$ denote the reflection in $\Sigma$. We
write $\overline{x}=\sigma_H(x)$ for points $x\in \mathbb{R}^n$
and $\sigma_H(u)=u(\overline{x})$ for all $x\in \mathbb{R}^n$ for
functions $u\in{\mathcal{S}}$. %

\bd \label{Defenition of polarization of function} %
If $u\in {\mathcal{S}}_n$, then its polarization $Pu$ with the
{\it
polarizer} $H$ is given by %
\be \label{3.7-equation} %
Pu(x):=\left\{  %
\begin{array}{ll} %
\max\{u(x),u(\overline{x})\} & \ \ {\mbox{if $x\in H$,}}  \\ %
\min\{u(x),u(\overline{x})\} & \ \ {\mbox{if $x\in
\mathbb{R}^n\setminus H$}}. %
\end{array}  %
\right. %
\ee  %
\ed %

If $M\in {\mathcal{M}}_n$, then the polarization $P(M)$ is given
by its characteristic function via (\ref{3.7-equation}), i.e.  %
\be \label{3.8-equation} %
\chi(P(M)):=P(\chi(M)). %
\ee  %
In the case that $u$ is continuous and $M$ is open or closed,
equations (\ref{3.7-equation}) and (\ref{3.8-equation}) have to be
understood in the pointwise sense.

Equations (\ref{3.7-equation}) and (\ref{3.8-equation}) can also
be written in the following more precise form %
\be  \label{3.9-equation} %
P(M)=((M\cup \sigma_H(M))\cap H)\cup(M\cap \sigma_H(M)), \quad
\quad M\in {\mathcal{M}}_n, %
\ee  %
and %
\be \label{3.10-equation} %
Pu(x)=\esssup \left\{c>\inf u:\, x\in P(\{u>c\})\right\}, \quad
\quad
x\in \mathbb{R}^n, \quad u\in {\mathcal{S}}_n. %
\ee  %
Of course, if $u\in {\mathcal{S}}_n$ is continuous, then
``$\esssup$'' in (\ref{3.10-equation}) can be replaced by
``$\sup$''. From the representations
(\ref{3.7-equation})--(\ref{3.10-equation}) we see that the
polarization $P$ is an open and compact rearrangement which is
continuous from the inside and from the outside.

For the sake of simplicity, we will often use the subscript
``${}_H$'' to denote any one of the polarized objects, i.e. we
write $u_H$ and $M_H$ for $Pu$ and $P(M)$, respectively.

It is worth mentioning, that the polarization of a connected set
is not necessarily connected and may contain one multiply connected component. %In Figure~\label{3-figure}, we show
%an example of a simply connected domain, the polarization of which
%is disconnected and has a multiply connected component. %

%\begin{figure}
%$$\includegraphics[scale=.3,angle=0]{iceberg} $$
%\caption{Polarization.} \label{Polarization}
%\end{figure}

\smallskip

There are three major approaches to polarization. The first one
initiated in \cite{W}, \cite{A}, and \cite{BaT} uses convolution
type inequalities. In the most powerful and general form this
approach culminated in Baernstein's fundamental work \cite{Ba}.

The second approach was introduced by V.~N.~Dubinin \cite{D1} who
used the following representation of polarized functions. Let
$v(x)=u(\overline{x})$, $w(x)=u_H(\overline{x})$, $x\in H$ with
$x\in H$. Then  %
$$  %
u_H(x)=(u(x)-v(x))_+ \quad \quad {\mbox{and}} \quad \quad
w(x)=u(x)-(u(x)-v(x))_+. %
$$  %
The latter under certain conditions leads to %
\begin{eqnarray*}%
\nabla u_H(x)&=&\left\{ %
\begin{array}{ll} %
\nabla u(x) & \ \ {\mbox{a.e. on $\{u>v\}\cap H$,}} \\
\nabla v(x) & \ \ {\mbox{a.e. on $\{u\le v\}\cap H$,}} %
\end{array}  %
\right. %
\\  %
\nabla w(x)&=&\left\{ %
\begin{array}{ll} %
\nabla v(x) & \ \ {\mbox{a.e. on $\{u>v\}\cap H$,}} \\
\nabla u(x) & \ \ {\mbox{a.e. on $\{u\le v\}\cap H$,}} %
\end{array}  %
\right. %
\end{eqnarray*}  %
which easily leads to several Dirichlet type inequalities, see
\cite[Lemma~5.3]{BS}.

\medskip

The third approach suggested by A.~Solynin \cite{S2} and developed
further in \cite{BS} rests on the direct application of the
maximum principle to prove comparison theorems for solutions of
two related boundary value problems for certain partial
differential equations defined in a given domain $\Om$ and in the
 polarized domain $\Om_H$. This approach will be used in
Section~\ref{Comparison theorems}.

\section{ Continuous $(1,n)$-Steiner symmetrization}   \label{ Continuous $(1,n)$-Steiner symmetrization}%
\setcounter{equation}{0}

First we define a continuous transformation on $\mathbb{R}$. For
$M\in {\mathcal{M}}_1$ and $-\infty <t\le \infty$, the {\it
measuring  function of $M$} is defined by  %
\be  \label{m.f.} %
m_M(t)=\mathcal{L}^1((-\infty,t)\cap M).  %
\ee  %
This definition shows that $m_M(t)$ is nondecreasing and Lipschitz
continuous with constant $1$:  %
\be  \label{Lip for m.f.}  %
0\le m_M(t_2)-m_M(t_1) \le t_2-t_1 \quad {\mbox{for all $t_1\le
t_2$.}}  %
\ee  %
This implies that, for every $t\in \mathbb{R}$, the equation %
\be \label{eq. for l.f.} %
y-(1/2)m_M(y)=t  %
\ee  %
has a unique solution $y=y_M(t)\in [t,\infty)$. The function
$y=y_M(t)$, called {\it the separating function of $M$,} will play
an important role in this study. Two basic properties of $y_M(t)$
given by the following lemma are immediate consequences of the
above definitions and inequalities (\ref{Lip for m.f.}).  %
\bl  \label{properties of l.f.}  %
{\bf (a)}  If $M\subset N$, then $y_M(t)\le y_N(t)$ for all $t\in
\mathbb{R}$.  \\%
\ \ \ \ {\bf (b)} If $t_1<t_2$, then %
\be  \label{bounds for l.f.}  %
t_2-t_1\le y_M(t_2)-y_M(t_1)\le 2(t_2-t_1).  %
\ee  %
\el  %

\smallskip

For $t\in \mathbb{R}$, the $T^t$-transformation %symmetrization
of $M$ is defined as %
\begin{alignat}{10}     %
M^t &= &(y_M(t)-m_M(t),y_M(t))\cup (M\cap [y_M(t),\infty))
{\mbox{\ if $M$ is open, }}\label{M^t-open} \ \ \ \  \\ %
%\ee  %
%\be   %
M^t&=&\ [y_M(t)-m_M(t),y_M(t)]\cup (M\cap [y_M(t),\infty))
{\mbox{\  if $M$ is compact. }}\label{M^t-compact}%
\end{alignat}  %
If $M\in {\mathcal{M}}_n$ but is neither open or compact then
$M^t$ is defined in the a.e. sense by either one of the equations
(\ref{M^t-open})
or (\ref{M^t-compact}).   %Otherwise, $M^t$ is defined by (\ref{M^t-compact}) in a vague sense.

\bd  \label{SCS in R}  %
The family of mappings  $T^t:{\mathcal{M}}_{1}\to
{\mathcal{M}}_1$, $t\in \mathbb{R}$, defined by $T^t(M)=M^t$ is
called
the continuous symmetrization on $\mathbb{R}$.  %
\ed  %
Now we turn to $\mathbb{R}^n=\mathbb{R}^{n-1}\times \mathbb{R}$.
 For $\Om\in {\mathcal{M}_n}$ and $x'\in \mathbb{R}^{n-1}$,
let $\Om(x')=\{y\in \mathbb{R}:\,(x',y)\in \Om\}$  be  the
$1$-slice of $\Om$ at $x'$. If $\Om\in {\mathcal{M}}_n$ is open or
compact then $\Om'$ will denote the orthogonal projection of $\Om$
onto $\mathbb{R}^{n-1}$. Otherwise we put $\Om'=\{x'\in
\mathbb{R}^{n-1}:\,{\mathcal{L}}^1(\Om(x'))\not= 0\}$. The
measuring function
of $\Om$ is defined by %
\be  \label{m(x',t)} %
m_\Om(x',t)=\left\{ \begin{array}{cl} %
m_{\Om(x')}(t) & {\mbox{if $x'\in \Om'$}} \\ 0 &
{\mbox{otherwise.}}
\end{array} \right. %
\ee  %
One can easily show that $m_\Om(x',t)$ is lower semicontinuous in
$x'$ if $\Om$ is open and it is upper semicontinuous in $x'$ if
$\Om$ is compact.

\bd \label{def of l.f.}  %
The function $y_\Om:\mathbb{R}^{n-1}\times \mathbb{R}\to \mathbb{R}$ defined by %
\be  \label{l.f. in R^n}  %
y_\Om(x',t)=\left\{ \begin{array}{cl} y_{\Om(x')}(t)& \quad {\mbox{if $x'\in \Om'$}}\\ t& \quad {\mbox{otherwise}}  %
\end{array} \right.  %
\ee  %
is called the separating function and the graph
$F_\Om(t)=\{(x',y_\Om(x',t)):\,x'\in \mathbb{R}^{n-1}\}$ is called
the frontier of symmetrization
 of $\Om$ at $t$.   %
\ed  %

To simplify notation we will skip  the symbol of the set if its
meaning is clear from the context. Thus we often write $m(x',t)$,
$y(x',t)$, etc. instead of $m_\Om(x',t)$, $y_\Om(x',t)$, etc.
%An example of the frontier of symmetrization is shown in Figure~\ref{4-figure}, which also illustrates some other notation of this section.

%\begin{figure}
%$$\includegraphics[scale=.3,angle=0]{iceberg} $$
%\caption{One dimensional continuous symmetrization.}
%\label{4-figure}
%\end{figure}

\bl \label{semicont. of l.f.}  %
For a fixed $t\in \mathbb{R}$ and $\Om\in{\mathcal{M}}_n$, the
separating function $y(x',t)$ is lower semicontinuous  on
$\mathbb{R}^{n-1}$ if $\Om$ is open and
it is upper semicontinuous on $\mathbb{R}^{n-1}$ if $\Om$ is  compact.  %
\el %

{\it Proof.} Let $\Omega $ be open and $x'_0\in \mathbb{R}^{n-1}$.
If $m(x_0',t)=0$, then $y(x',t)\ge t=y(x_0',t)$ for all $x'$ and
the lower semi-continuity follows.

Assume that $y(x',t)$ is not lower semicontinuous at $x_0'$ such
that $m(x_0',t)>0$. Then for some $\delta>0$ and some sequence
$x'_k\to x_0'$,
 \be \label{y(x'_k)<y(x_0')}  %
 y(x'_k,t)\le y(x'_0,t)-\delta,\qquad k=1,2,\ldots .
\ee %
From (\ref{eq. for l.f.}), (\ref{Lip for m.f.}), and
(\ref{y(x'_k)<y(x_0')}) we obtain %
\begin{alignat}{10} %  %
m(x'_k,y(x_0',t)-\delta)- 2(y(x'_k,t)-t)   &=&\
m(x'_k,y(x_0',t)-\delta) - m(x'_k,y(x'_k,t))\label{lemma 2, formula 2} \\
{} &\le& y(x'_0,t)-y(x'_k,t)-\delta\,. \nonumber   \ \ \ \ \ \ \ \
\ \ \ \ \ \ \ \ \ \ %
\end{alignat}  %
This combined with ({\ref{y(x'_k)<y(x_0')}) gives %
\begin{equation} \label{limsup< in Lemma 2} %
\limsup_{k\to \infty} m (x'_k,y(x'_0,t)-\delta)\le
2(y(x'_0,t)-t-\delta)\,.
\end{equation}
Since $\Omega $ is open  the function $m(x',t)$ is lower
semicontinuous in $x'$. Hence  %
 \be \label{liminf > in Lemma 2}  %
 \liminf_{k\to\infty}  m(x'_k,y(x'_0,t)-\delta)\ge
m(x'_0,y(x'_0,t)-\delta).  %
\ee  %
Using (\ref{Lip for m.f.}) and (\ref{eq. for l.f.}) we obtain  %
\be \label{formula 5 Lemma 2}  %
 m(x'_0,y(x'_0,t)-\delta)\ge m(x'_0,y(x'_0,t))-\delta =
2(y(x'_0,t)-t)-\delta  %
\ee  %
Now (\ref{liminf > in Lemma 2}) and %
(\ref{formula 5 Lemma 2})  %
 yield the inequality
$$ %
\liminf_{k\to \infty} m(x'_k,y(x'_0,t)-\delta)\ge 2(y(x'_0,t)-t)
-\delta,
$$  %
which contradicts  (\ref{limsup< in Lemma 2}). This proves
Lemma~\ref{semicont. of l.f.} for the case of open sets. If $\Om$
is compact the proof is similar and is left to the reader. \hfill
$\Box$

Now we are ready to define what the  $1$-dimensional continuous symmetrization is. %
\bd  \label{SCSym}  %
A family of set transformations  %
$T^t:{\mathcal{M}}_{n}\to {\mathcal{M}}_n$,  $t\in \mathbb{R}$,  %
\   defined by
\be \label{T^t}  %
T^t(\Om)=\Om^t:=\{(x',y):\,x'\in \Om',\  y\in \Om^t(x')\} %
\ee  %
will be called a continuous $1$-symmetrization or  SC
$1$-symmetrization. Any single transformation of this family will
be called a $T^t$-transformation. For $t=-\infty$, we define
$T^{-\infty}$ to be the identity transformation.
\ed %

The set $\Om^t$ itself will be called the {\it
$T^t$-transformation of $\Om$} or {\it the partial symmetrization
of $\Om$ }with respect to the plane $\{y=t\}$. Sometimes we will
refer to this plane $\{y=t\}$ as {\it the moving plane of
symmetrization}. Accordingly, the parameter $t$ will be called
{\it the height} of the moving plane.

\medskip

The frontier of symmetrization $F_\Om(t)$ divides $\mathbb{R}^n$
into two parts %
\be \label{subspases of Sym}  %
H_+(t)=\{(x',y):\,y>y(x',t)\} \quad {\mbox{and}} \quad
H_-=\{(x',y):\,y<y(x',t)\}  %
\ee  %
called the upper and lower subspaces of symmetrization,
respectively.  Let %
\be \label{Om_(t) for open} %
\Om_-(t)=\Om\cap H_-(t), \quad \Om_+(t)=\Om\setminus H_-(t) \quad
{\mbox{if $\Om$ is open }} %
\ee  %
and %
\be \label{Om_(t) for compact} %
\Om_+(t)=\Om\cap H_+(t), \quad \Om_-(t)=\Om\setminus H_+(t) \quad
{\mbox{if $\Om$ is compact. }} %
\ee  %
If $\Om$ is neither open or compact then $\Om_-(t)$ and $\Om_+(t)$
are defined a.e. by (\ref{Om_(t) for compact}).

In Lemma~\ref{Basic properties of SCSym} below we list some useful
properties of SC symmetrization, which follow directly from the
above definitions. By $\Om^*_-(t)$ we denote the $(1,n)$-Steiner
symmetrization of $\Om_-(t)$ with respect to $\{y=t\}$.  %
\bl \label{Basic properties of SCSym} %
\begin{enumerate} %
\item[(a)]  %
%\begin{equation}
$\Omega^t=\Omega^*_- (t)\cup \Omega_+(t)$ for all $t\in
{\mathbb{R}}$ and every $\Omega \in {\mathcal{M}}_n$.
%\end{equation}  %
\item[(b)]  %
 ${\mathcal{L}}^1(\Omega ^t(x'))={\mathcal{L}}^1(\Omega(x'))\ $
for all $x'\in \mathbb{R}^{n-1} $ and
 ${\mathcal{L}}^n(\Omega^t)={\mathcal{L}}^n(\Omega)$. %
 \item[(c)]  %
If $\Omega $ is open, then $\Om_-(t)$, $\Om_-^*(t)$, and
$\Omega^t$ are open sets and \\
$\Omega^*_-(t)\cap\{y\ge t\}=H_-(t)\cap\{y\ge t\}$. %
\item[(d)]  %
If $\Omega $ is compact, then $\Om_-(t)$, $\Om_-^*(t)$, and
$\Omega^t $ are compact sets and \\
$\Omega^*_-(t)\cap\{y\ge t\}=\{y\ge t\}\setminus H_+(t)$.  %
\end{enumerate}  %
\el  %
 We note that Lemma~\ref{Basic properties of SCSym}(a) gives an equivalent definition
of $\Om^t$ and Lemma~\ref{Basic properties of SCSym}(b) shows that
the $T^t$-transformation  preserves measures.

%%          The following lemma gives a stability property of the $T^t$-transformation. %%

\medskip

Now we will  show that SC symmetrization  possesses all
basic geometric properties of classical symmetrizations. %
\bl \label{Basic geometric properties of SCSym} %
 For a fixed $t\in \mathbb{R}$, the $T^t$-transformation is an open, compact,
 and smoothing rearrangement,  which is continuous from the inside and from the outside.
\el %

{\it Proof.} Lemma~\ref{properties of l.f.}(a) and
Definition~\ref{SCSym} imply that $T^t(M)\subset T^t(N)$ if
$M\subset N$ and therefore the $T^t$-transformation is monotone.
Since, by Lemma~\ref{Basic properties of SCSym}(b), the
$T^t$-transformation preserves the measure of sets, it follows
that the $T^t$-transformation is a rearrangement.

\smallskip

 From Definition~\ref{SCSym} and  Lemma~\ref{Basic properties of SCSym}(c) and (d),
   one can easily  verify that the $T^t $-transformation is continuous from the inside and from the
outside.  %

\smallskip

To prove that the $T^t$-transformation is smoothing, we fix a
compact set $K$ and $\ep>0$. Since, by Lemma~\ref{Basic properties
of
SCSym}(a), $K^t=K^*_-(t)\cup\overline{K_+(t)}$, we have  %
\be \label{formula 1 of Lemma 4}  %
 K^t+\e \overline{B^{(n)}}=(K^*_-(t)+\e \overline{B^{(n)}}) \cup
(\overline{K_+(t)}+\e \overline{B^{(n)}}). %
\ee  %
 Using  Lemma~\ref{Basic properties of SCSym}(d) and the smoothing property of
the Steiner symmetrization we obtain %
\be \label{formula 2 of Lemma 4}  %
T^t(K+\e \overline{B^{(n)}})\supset T^t(K_-(t)+\e
\overline{B^{(n)}})=(K_-(t)+\e \overline{B^{(n)}})^*\supset
K^*_-(t)+\e \overline{B^{(n)}}, %
\ee  %
where $(\cdot)^*$ denotes the Steiner symmetrization with respect
to $\{y=t\}$. Since $K_+(t)$ remains unchanged under
the $T^t$-transformation of $K$ one can easily see that %
\be \label{formula 3 of Lemma 4}  %
T^t(K+\e \overline{B^{(n)}})\supset T^t(K_+(t)+\e \overline{B^{(n)}})
=\overline{K_+(t)}+\e \overline{B^{(n)}}.  %
\ee  %
Now  (\ref{formula 1 of Lemma 4}) -- (\ref{formula 3 of Lemma 4})
yield:  %
$$  %
T^t(K+\e \overline{B^{(n)}})\supset T^t(K)+\e \overline{B^{(n)}},  %
$$ %
which shows that the $T^t $-transformation is smoothing.  %

\smallskip

It was shown by Sarvas, see Lemmas~3.2 and 3.3 in \cite{BS}, that
the properties established above imply that the
$T^t$-transformation is open and compact. The proof is complete.
$\hfill \Box$

\medskip

 Since for any fixed $t\in \mathbb{R}$,  the $T^t $-transformation is a  rearrangement of sets
 there is a standard way  to define a corresponding transformation of
 functions:  %

\bd \label{SCSym of functions}  %
SC $1$-symmetrization of functions: Let $u\in {\mathcal{S}}_n$.
Then the  family of functions $u^t$,
$t\in \mathbb{R}$, defined for $x\in \mathbb{R}^n$ by %
$$ % \label{SCSym of functions formula} %
 u^t(x)%\equiv T^tu(x)
 :=  \left\{ \begin{array}{ll} %
 \esssup \left\{ c>\inf u:\ x\in T^t( \{ u>c\})\right\} &
{\mbox{if \ }} x\in \bigcup_{c>\inf u} T^t(\{ u>c\})\\
\inf u & {\mbox{if \ }} x\not\in \bigcup_{c>\inf u} T^t(\{ u>c\}) %
\end{array} \right.%
$$  %\ee  %
is called the SC $1$-symmetrization of $u$. %
\ed  %

\medskip

\section{Continuity properties of SC symmetrization}   \label{Continuity properties of SC Symmetrization}%
\setcounter{equation}{0}  %

 Now we  prove that the  $T^t$-transformation
depends continuously on  the parameter $t$, which justifies the
use of
the word "continuous" in the name of this transformation. %
\bl \label{Continuity of T^t} %
 Let $t_k$, $k=1,2,\ldots$ be an increasing sequence such that
 $t_k\to  t'<\infty$ and let $\Omega'=\Omega^{t'}$ and $\Omega_k=\Omega^{t_k}$.
 Then  %
\be \label{formula 1 of Lemma 6} %
 d(\partial \Omega_k,\partial \Omega')\to 0 \quad {\rm {as}}
\quad k\to\infty \ee %
 if $\Omega$ is an open bounded
set and  %
\be \label{formula 2 of Lemma 6}  %
 d(\Omega_k,\Omega')\to 0 \quad {\rm {as}} \quad k\to \infty
\ee  %
if $\Omega$ is a compact set. %
\el  %

{\it Proof.} We will prove the lemma for open sets. The proof for
compact sets follows the same lines, it is easier and is left to
the reader.

If (\ref{formula 1 of Lemma 6}) does not hold, then there exist
subsequences of the sequences $t_k$ and $\Om_k$, which we still call $t_k$ and $\Omega_k$, such that  %
\be \label{formula 3 of lemma 6}  %
 d(\partial
\Omega_k,\partial\Omega')\ge \varepsilon, \quad \quad k=1,2,\ldots  %
\ee  %
for some $\varepsilon>0$. Since $\Om$ is bounded
equation~(\ref{formula 3 of lemma 6}) implies that  there are
subsequences (again denoted by $t_k$ and $\Omega_k$) and  a
sequence of
 points $x_k\in \mathbb{R}^n$ with
$ x_k\to x_0\in \mathbb{R}^n $ such that one of the following two
conditions is satisfied:

\begin{enumerate}  %
\item[(a)]%
 $x_k\in \partial\Omega'$ and ${\rm
dist}(x_k,\partial\Omega_k)\ge \varepsilon$ for all
$k=0,1,2,\ldots$  %
\vspace{0.2 cm} %
\item[(b)]  %
 $x_k\in \partial \Omega_k$ for $k=1,2,\ldots$ and
${\rm dist}(x_k,\partial \Omega')\ge \varepsilon$ for
$k=0,1,2,\ldots$  %
\end{enumerate}  %

The case (a) contradicts Lemma~\ref{Stability Lemma} below, which
characterizes  the stability property of the boundary of
$\Omega^t$.

Now we consider the case (b). Let  %
$$  %
x_k=(x'_k,y_k)\to (x'_0,y_0)=x_0.  %
$$  %
Then  %
$$  %
x_k\in \partial(\Omega^*_-(t_k))\cup \partial(\Omega_+(t_k)).  %
$$  %
Taking  a subsequence if necessary we may restrict ourselves to
 two cases:

(1) In the first case we assume that  $x_k\in \partial
\Omega_+(t_k)$ for all $k=1,2,\ldots$.
 Since $x_k\not\in \Omega_k$ we have $x_k\not\in\Omega_+(t_k)$.
 Now if  $x_k\in \Omega$ then $x_k\in \Omega_-(t_k)$. Since $\Om_-(t_k)$ is open
 and $\Om_-(t_k)\cap \Om_+(t_k)=\emptyset$ we conclude that $x_k\not\in
\partial\Omega_+(t_k)$, which  contradicts the above assumption.
 Thus in the case under consideration
  $x_k\not\in \Omega$ and therefore  $x_k\in \partial\Omega$ for all $k=1,2,\ldots$. Taking the limit we get
\be \label{formula 4 of Lemma 6}  %
x_0\in \partial \Omega.  %
\ee  %
 Since $x_k\in\partial \Omega_+(t_k)$ using (\ref{bounds for l.f.}) we obtain
$$  %
y_k\ge y(x'_k,t_k)\ge y(x'_k,t')-2(t'-t_k).  %
$$  %
Passing to the limit  we get  %
\be \label{formula 5 of Lemma 6} %
y_0\ge y(x'_0,t')  %
\ee  %
since $y(x',t)$ is lower semicontinuous in $x'$.
Equations~(\ref{formula 4 of Lemma 6}) and (\ref{formula 5 of
Lemma 6}) imply that  %
\be \label{formula 6 of Lemma 6}%
x_0\not\in \Omega'.  %
\ee  %

Now we will  investigate possible dispositions of $y_0$,
$y(x'_0,t')$, and $t'$.

(i) \, If $y_0=y(x'_0,t')>t'$, then $x_0\in
\partial \Omega'$  contradicting the assumption~(b).

\smallskip

(ii) \, The second possibility is %
$$  %
y_0>y(x'_0,t')\ge t'.  %
$$  %
It follows from the assumption~(b) and (\ref{formula 6 of Lemma
6}) that there exists a ball $B^{(n)}_{\e_1}(x_0)$ with a small
radius $\varepsilon_1>0$ such that  %
\be \label{formula 7 of Lemma 6}  %
 \overline{B^{(n)}_{\e_1}}(x_0)\subset \left( {\mathbb{R}}^n
\setminus \Omega' \right) \cap H_+(t').  %
\ee  %

By  (\ref{formula 4 of Lemma 6}), there exists a point $\hat x\in
\Omega\cap B^{(n)}_{\e_1}(x_0).$ Since $B^{(n)}_{\e_1}(x_0)\subset
H_+(t')$ we must have $\hat x \in \Omega'.$ But this contradicts
(\ref{formula 7 of Lemma 6}) and therefore contradicts~(b).

\smallskip

(iii) \, Let $y_0=t'.$  The assumption~(b) and (\ref{formula 6 of
Lemma 6}) imply that there is a ball $B^{(n)}_{\e_1}(x_0)$ such
that $B^{(n)}_{\e_1} (x_0)\subset \mathbb{R}^n\setminus \Omega'$.
Then it follows from the definition of $\Om^{t'}$ that
$B^{(n)}_\e(x_0)\subset \mathbb{R}^n\setminus \Om$, which
contradicts (\ref{formula 4 of Lemma 6}). This finishes the proof
of the lemma in the case (1) under the
assumption~(b).  %

\medskip

(2) In the second case we assume that %
\be \label{formula 9 of Lemma 6}  %
x_k\in \partial (\Omega^*_-(t_k))\setminus \partial \Omega_+(t_k)
\qquad {\mbox{for all $k=1,2,\ldots$}}  %
\ee  %
As in the case (1) we consider possible dispositions of $y_0$,
$y(x'_0,t')$, and $t'$.

\smallskip

(i) \, If $y_0<t'$, then $y_k<t_k$ for all sufficiently large $k$.
Since $x_k\in \partial (\Om_-^*(t_k))$ and since the open interval
between the points $(x'_k,2t_k-y(x'_k,t_k))$  and
$(x'_k,y(x'_k,t_k))$ is in
$\Om_-^*(t_k)$ we have  %
\be \label{formula 8 of Lemma 6} %
y_k\le 2t_k-y(x'_k,t_k)\le 2t'-y(x'_k,t')  %
\ee %
for all $k$ large enough.  The second inequality here follows from
(\ref{bounds for l.f.}). The latter inequality implies that
$x_k\not\in \Omega'$ and therefore $x_0$ is not an inner point of
$\Omega'$.

Let us show that $x_0$ cannot be an outer point of $\Omega'$, too.
 Since $x_k\in \partial \Omega^*_-(t_k)$ and $x_k\to x_0$  there is a sequence of points
 $\hat x_k=({\hat x}'_k,\hat y_k)\in
\Omega^*_-(t_k)$ such that $\hat x_k \to x_0$. Then we have
$2t_k-y(\hat{ x}'_k,t_k)<\hat{y}_k$. Therefore using (\ref{bounds
for l.f.}) we obtain  %
$$  %
 2t'-y(\hat{ x}'_k,t')\le 2t_k-y(\hat{x}'_k,t_k)+(t'-t_k)<(t'-t_k)+\hat{
y}_k.  %
$$  %
This yields $(\hat{ x}'_k,(t'-t_k)+\hat{y}_k)\in \Omega'$. Since
$x_0\not\in \Om'$ taking the limit we get  $x_0\in \partial
\Omega'$, which  contradicts the assumption~(b).

\smallskip

(ii) \, Let $y_0>t'$. Then $y(x'_k,t_k)\le y_k$ for all $k$ large
enough. Therefore by (\ref{bounds for l.f.}), %
$$  %
 y(x'_k,t')\le 2(t'-t_k)+y(x'_k,t_k)\le 2(t'-t_k)+y_k.  %
$$  %
 Passing to the limit and using the semi-continuity property of the separating function
 we obtain  %
\be \label{formula 10 of Lemma 6}  %
y(x'_0,t')\le y_0.  %
\ee  %

The case of equality $ y(x'_0,t')=y_0\ (>t')$ leads to the
inclusion $x_0\in
\partial \Omega'$, which contradicts the assumption~(b). Indeed,
if $x_0\not\in \partial \Om'$, then $x_0\in \Om'$. Therefore
$x_k\in \Om$ for all sufficiently large $k$. Since $x_k\in
\partial \Om_-^*(t_k)$ and $x_k\in \Om$, it follows that $x_k\in
\Om_+(t_k)$. Therefore $x_k\in \partial \Om_+(t_k)$ for all $k$
large enough, which contradicts (\ref{formula 9 of Lemma 6}).

Now let  $y(x'_0,t')<y_0$. If $x_0\not\in \Omega'$, then it
follows from the assumption~(b) that there is  a ball
$B^{(n)}_{\e_2}(x_0)$ with
some small radius $\varepsilon_2>0$ such that %
$$  %
 \overline{B^{(n)}_{\e_2}}(x_0)\subset \left( {\mathbb{R}}^n \setminus  \Omega'
\right) \cap H_+(t').  %
$$  %
 Then $\overline{B^{(n)}_{\e_2}}(x_0)\cap
\Omega=\emptyset$. Since $\overline{B^{(n)}_{\e_2}}(x_0)\subset
H_+(t_k)$ we must have  $\overline{B^{(n)}_{\e_2}}(x_0)\cap
\Omega_k=\emptyset$ for all $k=1,2,\ldots$ The latter is
impossible since $x_k\in
\partial \Omega_k$ and $x_k\to x_0$.

\smallskip

If $x_0 \in \Omega'$ and $y(x'_0,t')<y_0$, then $x_0\in \Omega.$
Hence, $x_k\in \Omega$ for all sufficiently large $k$. Since
$x_k\in \partial \Omega_-^*(t_k)$ this implies that $x_k\in
\Omega_k$ contradicting the assumption~(b) if $k$ is large enough.
This completes the proof of the lemma in the case $y_0>t'.$

\smallskip

(iii) \, Let $y_0=t'$. If $x_0\not\in \Omega'$ then by the
assumption~(b), $\overline{B^{(n)}_\e}
(x_0)\cap\Omega'=\emptyset$. This implies that
$\overline{B^{(n)}_\e}(x_0)\cap\Omega^*_-(t_k)=\emptyset$ for all
$k=1,2,\ldots$ The latter contradicts   the assumptions $x_k\in
\partial\Omega^*_-(t_k)$ and $x_k\to x_0$.

Let $x_0\in \Omega'$. Then by the assumption~(b), $B^{(n)}_\e(x_0
) \subset \Omega^*_-(t').$ This together with the inequality
(\ref{bounds for l.f.}) implies that $B^{(n)}_{\e/2}(x_0)\subset
\Omega^*_-(t_k)$ for all $k$ large enough, again   contradicting
the assumption~(b). The proof of the lemma for open sets is now
complete. $\hfill \Box$

\medskip

Figure~3 shows the results of continuous symmetrization of a
domain $\Omega$ for a few values of $t$. Comparing domains in
Figure~3c and Figure~3d, the reader may see that the
$T^t$-transformation is not continuous from the right in the
Hausdorff metric.

\begin{figure}  %

\caption{$T^t$-transformation of sets.}
\end{figure} %

\medskip

Now we will prove the lemma  characterizing  a stability property
of the boundary
$\partial \Omega^t$.  %
\bl \label{Stability Lemma}  %
 Let $-\infty<s<t<\infty$. Then %
\be \label{formula 1 of lemma 7} %
\partial \Omega^t\subset \partial \Omega^s +(t-s)\overline{B^{(n)}} %
\ee %
 if $\Omega$ is open and  %
\be \label{formula 2 of lemma 7} %
\Omega^t\subset \Omega^s +(t-s)\overline{B^{(n)}} %
\ee  %
if $\Omega$ is a compact.   %
\el  %

{\it Proof.}  Once more we  give a proof for open sets. The easier
case of compact sets is left to the reader.

Suppose that $\Omega$ is open and that (\ref{formula 1 of lemma
7}) is not satisfied. Then there exists a point $x_0\in
\partial \Omega^t$ such that %
\be  \label{formula 3 of lemma 7} %
{\rm dist}(x_0,\partial \Omega^s)>(t-s).  %
\ee %
 By Lemma~\ref{Basic properties of SCSym}(a), we have
$\Omega^t=\Omega^*_-(t)\cup\Omega_+(t).$ Hence, %
$$  %
\partial\Omega^t\subset \partial \Omega^*_-(t)\cup \partial
\Omega_+(t).%
$$  %

\smallskip

As in the proof of Lemma~\ref{Continuity of T^t} we consider two
cases.

(1) \, In the first case we assume that  %
\be  \label{formula 4 of lemma 7} %
x_0\in \partial\Omega_+(t).  %
\ee  %
 By Lemma~\ref{Basic properties of SCSym}(c), the set  $\Omega_-(t)$ is open.
 By the definitions of $\Om_-(t)$ and $\Om_+(t)$,
 we have
 $\Omega_-(t)\cap \Omega_+(t)=\emptyset$.  Therefore (\ref{formula 4 of lemma 7}) implies that
  $x_0\not\in\Omega_-(t)$.

Since $\Omega_+(t)\subset \Omega$ we must have $x_0\in \overline{
\Omega }$. If $x_0\in \Omega$, we conclude that $x_0\in
\Omega_+(t)$ and therefore $x_0\in \Omega^t$. Since by
Lemma~\ref{Basic properties of SCSym}(c), $\Om^t$ is open, the
latter contradicts the assumption $x_0\in \partial \Om^t$.
Therefore we must have %
\be \label{formula 5 of lemma 7}  %
x_0\in \partial\Omega.  %
\ee  %
Since  $ \Omega_+(t)\subset \Omega_+(s)\subset \Omega$, equations
 (\ref{formula 4 of lemma 7}) and (\ref{formula 5 of lemma 7})
 imply %
\be \label{formula 6 of lemma 7}  %
 x_0\in \partial \Omega_+(s).  %
\ee  %
 Hence, $x_0\in \overline{\Omega^s}$.  This implies that $x_0\in \partial \Om^s$.
 Indeed, let $x_0\in \Omega^s$.
 Since $x_0\not\in \Omega$ we conclude that $x_0\not\in \Omega_+(s).$
Therefore in the case under consideration we must have  %
\be \label{formula 7 of lemma 7}  %
x_0\in \Omega^*_-(s).  %
\ee  %
  By Lemma~\ref{Basic properties of SCSym}(c), the set
$\Omega^*_-(s)$ is open and it follows from the definitions that
$\Omega^*_-(s)\cap\Omega_+(s)=\emptyset$. Therefore (\ref{formula
7 of lemma 7}) contradicts  (\ref{formula 6 of lemma 7}) and we
must have $x_0\in \partial \Om^s$ contradicting (\ref{formula 3 of
lemma 7}).

\medskip

(2) \, In the second case we assume that $x_0=(x'_0,y_0)\in
\partial\Omega^*_-(t)\setminus\partial\Omega_+(t)$. Now we
consider two subcases.

\smallskip

(i) \, Suppose that $y_0\le t$. Since $x_0\in \Om^t$ we have  %
$$  %
y_0\le 2t-y(x'_0,t).  %
$$  %
This together with (\ref{bounds for l.f.}) implies  %
$$  %
y_0-(t-s)\le 2s-y(x'_0,s).  %
$$   %
This inequality yields  %
\be \label{formula 8 of lemma 7}  %
\tilde x:=( x'_0,y_0-(t-s))\not\in \Omega^s.  %
\ee  %
 Since $x_0\in\partial\Omega^*_-(t)$ then for each
$\varepsilon>0$ there exists a point $\hat x=(\hat x',\hat
y)\in\Omega^*_-(t)$ such that %
\be \label{formula 9 of lemma 7}  %
 |\hat x-x_0|<\varepsilon.  %
\ee  %

Let $y_0<s$. Then $m(\hat x',\hat y)>0$ and therefore  $m(\hat
x',s)>0$ if $\hat x$ is close enough to $x_0$. Therefore the
interval  %
$$  %
((\hat x',y):\ 2s-y(\hat x',s)<y<y(\hat x',s))  %
$$  %
is nonempty and belongs to $\Omega^s$.

Using (\ref{bounds for l.f.}) and the relation $\hat x \in
\Om_-^*(t)$ we conclude that %
$$  %
2s-y(\hat x',s)\le 2t-y(\hat x',t)<\hat y<s  %
$$  %
if $\ep>0$ in (\ref{formula 9 of lemma 7}) is small enough. Then
we have $\hat x\in \Omega^s.$  Therefore the closed interval
$[\hat x,\tilde x]$ contains a point $z\in
\partial \Omega^s$. Hence, %
$$  %
|z-x_0|\le \max\{|x_0-\tilde x|,\ | x_0-\hat x|\}=(t-s)  %
$$  %
if $\ep>0$ in (\ref{formula 9 of lemma 7}) is $\le t-s$, which
contradicts inequality~(\ref{formula 3 of lemma 7}).

Let $y_0=s$ and let $\hat y\le s$. Then $m(\hat x',s)>0$ and
$(\hat x',s)\in \Omega^s$. The closed interval $[x_0,(\hat x',s)]$
contains a point $z\in \partial \Omega^s.$ Hence, $|x_0-z|\le
\varepsilon$, which again  contradicts (\ref{formula 3 of lemma
7}) if $\ep <t-s$.

Let $y_0=s,\ \hat y>s$. If $m(\hat x',s)>0$, then we argue as in
the previous case.

Let $m(\hat x',s)=0$. Since $m(\hat x',\hat y)>0$ then there
exists a point $z=(\hat x',\zeta)\in \Omega$ such that $s<\zeta\le
\hat y$. But in this case $z\in \Omega^s$ and $|z-x_0|\le |\hat
x-x_0|\le \varepsilon$. Since $\ep>0$ is arbitrary small this
easily leads to  a contradiction to (\ref{formula 3 of lemma 7}).

Let $s<y_0\le t$. If $m(\hat x',s)>0$, then $(\hat x',s)\in
\Omega^s$ and  %
$$  %
|x_0-(\hat x',s)|\le |x_0-(x'_0,s)|+|(x'_0,s)-(\hat x',s)|\le
(t-s)+\varepsilon.  %
$$  %
Since $\varepsilon>0$ can be chosen  arbitrary small  the latter
inequality contradicts  (\ref{formula 3 of lemma 7}).

Let $m(\hat x',s)=0$. Since $m(\hat x',\hat y)>0$ then there
exists a point $z=(\hat x',\zeta)\in \Omega$ such that
$s<\zeta<\hat y$. Then $z\in \Omega^s$. In addition, %
$$   %
|z-x_0|\le (t-s)+\varepsilon.  %
$$  %

Since $\varepsilon>0$ can be chosen  arbitrary small the latter
inequality again contradicts (\ref{formula 3 of lemma 7}).

\smallskip

(ii) \,  Let $y_0>t$. Since  $x_0\not\in \Omega^t$ we have $
y_0\ge y(x'_0,t)\ge y(x'_0,s). $ This implies  that
 $x_0\not\in \Omega^s$.
 Since $x_0\in \partial \Om_-^*(t)$ we conclude  that for every $\ep>0$
  there exists a point $ \hat x=(\hat x',\hat y)\in \Omega^*_-(t)$ such
  that  %
\be \label{formula 10 of lemma 7} %
 |\hat x-x_0|<\varepsilon.  %
\ee  %
If $\ep>0$ is small enough we have $ t<\hat y<y(\hat x',t).$ %

If $y(\hat x',s)\ge \hat y$, then $\hat x \in
\overline{\Omega^s}$. Therefore (\ref{formula 10 of lemma 7})
contradicts  (\ref{formula 3 of lemma 7}).

Now we consider the case  $y(\hat x',s)< \hat y$. Let
$I:=((\hat{x'},y(\hat{ x'},t)),\ (\hat{ x'}, y(\hat{ x'},s))$ and
let $\Om_I=I\cap \Om$. It follows from (\ref{eq. for l.f.}) that
the linear measure of $\Om_I$ equals   $ 2(y(\hat{ x'},t)-y(\hat{
x'},s))-2(t-s). $ Since $\hat{ x}\in I$ we have %
$$  %
{\rm dist}(\hat{ x},\Omega_I)\le
{\mathcal{L}}^1(I)-{\mathcal{L}}^1(\Omega_I).  %
$$  %
This inequality together with (\ref{bounds for l.f.}) implies %
$$  %
{\rm dist}(\hat{ x},\Omega_I)\le 2(t-s)-(y(\hat{ x'},t)-y(\hat{
x'},s))\le t-s.  %
$$ %
Since $\Omega_I\subset \Omega^s$ there exists a point $z\in
\overline{\Omega^s}$ such that $|\hat{ x}-z|\le t-s$. Finally, we
have %
$$  %
|x_0-z|\le |x_0-\hat{ x}|+|\hat{ x}-z|\le (t-s)+\varepsilon.  %
$$ %
Since $\varepsilon>0$ can be chosen  arbitrary small the latter
inequality contradicts (\ref{formula 3 of lemma 7}). The proof of
the lemma for open sets is complete.   $\hfill \Box$

\section{Stability and semigroup properties of
$T^t$-transformation}  \label{Stability and semigroup properties
of $T^t$-transformation}  %
\setcounter{equation}{0} %

First we prove a lemma characterizing stability with respect to
polarization, cf.  \cite[Theorem~3]{S1}.  %

\bl \label{Stability for polarization}  %
For $s\in \mathbb{R}$, let $P_s$ denote the polarization with the polarizer $H_s=\{(x',y)\in \mathbb{R}^n:\,y>s\}$.
If  $t_1<t_2$ and $s\le (1/2)(t_1+t_2)$, then  %

\be \label{formula 1 of lemma 5.1}  %
 T^{t_2}\circ P_s\circ T^{t_1}=T^{t_2} \quad \quad {\mbox {on ${\mathcal{M}}_n$}}. %
 \ee  %
 \el  %

 {\it Proof}. Let $\Om \in {\mathcal{M}}$ be an open set and let
 $D=P_s(\Om^{t_1})$. We have to show that  %
 \be \label{formula 2 of lemma 5.1}  %
 D^{t_2}=\Om^{t_2}.  %
 \ee  %
 By Lemma~\ref{Basic properties of SCSym}(a), %
 \be  \label{formula 3 of lemma 5.1}  %
 D^{t_2}=D^*_-(t_2)\cup D_+(t_2), \quad \Om^{t_2}=\Om^*_-(t_2)\cup
 \Om_+(t_2).  %
 \ee  %
 Since the $T^t$-transformation and polarization both preserve measure
 in $1$-slices, it follows from (\ref{formula 3 of lemma 5.1}) that
 (\ref{formula 2 of lemma 5.1}) holds true if and only if %
 \be  \label{formula 4 of lemma 5.1}  %
 D_+(t_2)=\Om_+(t_2).  %
 \ee  %

Since $\Om^{t_1}\cap \{(x',y):\,
y<2t_1-y_\Om(x',t_1)\}=\emptyset$, it follows from the definition
of polarization and $T^t$-transformation that %
\begin{alignat}{10} %
P_s(\Om^{t_1})\ &\cap & \{(x',y):\, y\ge 2(s-t_1)+y_\Om(x',t_1)\} \label{formula 5 of lemma 5.1}\\
\nonumber =\Om^{t_1}&\cap & \{(x',y):\, y\ge
2(s-t_1)+y_\Om(x',t_1)\}
 \\ %
\nonumber =  \Om \  \ &\cap & \{(x',y):\, y\ge
2(s-t_1)+y_\Om(x',t_1)\}\, . %
\end{alignat}  %

Since $s\le (1/2)(t_1 +t_2)$, (\ref{bounds for l.f.}) implies %
\be  \label{formula 6 of lemma 5.1}  %
y_\Om(x',t_2)\ge y_\Om(x',t_1)+2(s-t_1).  %
\ee  %
Now from (\ref{formula 5 of lemma 5.1}) and (\ref{formula 6 of
lemma 5.1}) we conclude that %
$$  %
P_s(\Om^{t_1})\cap \{(x',y):\,y\ge y_\Om(x',t_2)\}=\Om\cap
\{(x',y):\, y\ge y_\Om(x',t_2)\},  %
$$  %
which implies (\ref{formula 4 of lemma 5.1}).

For open sets the proof is complete. If $\Om$ is not open the
proof follows the same lines.
   $\hfill \Box$

\medskip

Lemma~\ref{Stability for polarization} and Definition~\ref{SCSym
of functions} immediately imply:

\bc  \label{Stability of functions for polarization}  %
Let $u(x)\in {\mathcal{S}}_n$,  $t_1<t_2$, and let  $s\le
(1/2)(t_1+t_2)$. Then  %
$$  %
\left((u^{t_1})_{H_s}\right)^{t_2}=u^{t_2}.  %
$$  %
\ec  %

In the particular case $s=t_1$, Lemma~\ref{Stability for
polarization} and its corollary show that our $T^t$-transformation
possesses the following ``presemigroup property'', which was first
mentioned in \cite{S1}. A similar property was used in \cite{B2}
to define the BC symmetrization.

\bc  \label{semigroup property}  %
Let $T^t$ denote a $T^t$-transformation on ${\mathcal{M}}_n$ or on
${\mathcal{S}}_n$.  Then  %
\be  \label{semigroup formula} %
T^{t_2}\circ T^{t_1}=T^{t_2} \qquad {\mbox{for all  \ \
$t_1<t_2$.}}  %
\ee  %
\ec  %

To justify the name ``presemigroup property'', we consider the
transformation ${\widehat{T}}^t=St^t\circ T^t$, where $t\ge 0$ and
$St^t$ denotes a continuous shift of $t$  units in the direction
of the negative $y$-axis. Then one can easily see that equation
(\ref{semigroup formula}) is equivalent to the following familia
semigroup property of ${\widehat{T}}^t$:  %
$$ %
{\widehat{T}}^{t+s}={\widehat{T}}^t\circ {\widehat{T}}^s \quad
{\mbox{for all $t\ge 0$, $s\ge 0$.}} %
$$ %
In Section~8, a similar composite transformation will be used to
replace the parameter $t$, the range of which is $-\infty< t<
\infty$, by a parameter $\tau$ with the standard range $0\le
\tau\le 1$.

\medskip  %

Now we prove a criterion for $\Om$ and $u$ to possess a partial
symmetry .

\bl  \label{Criterion of partial symmetry for functions}  %
(a) If $\Om\in {\mathcal{M}}_n$ is open or compact, then
$\Om=\Om^t$ if
and only if $P_s(\Om)=\Om$ for all $s\le t$. %

(b) If $u\in{\mathcal{S}}_n$, then $u\equiv u^t$ if and only if $u_{H_s}\equiv u$ for all $s\le t$.  %
\el  %

{\it Proof}.  %
(a) We work with an open set $\Om$. The necessary part is obvious:
If $\Om=\Om^t$, then Lemma~\ref{Basic properties of SCSym}(a)
shows that $P_s(\Om)=\Om$ for all $s\le t$.

\smallskip

To prove the sufficient part, we assume that $P_s(\Om)=\Om$ for
all $s\le t$ but $\Om\not=\Om^t$. Since $\Om_+(t)=(\Om^t)_+(t)$ we
conclude that $\Om_-(t)\not=\Om_-^*(t)$. Then by
Lemma~6.3 in \cite{BS},  %
there is a polarizer $H_{s_0}=\{y>s_0\}$ with $s_0<t$ such that %
\be  \label{formula 1 of lemma 5.2} %
P_{s_0}(\Om_-(t))\not=\Om_-(t).  %
\ee  %
Then there is $x_0=(x'_0,y_0)$ with $y_0<s_0$ such that %
\be \label{formula 2 of lemma 5.2} %
x_0\in \Om_-(t) \quad {\mbox{and}} \quad \hat
x:=(x'_0,2s_0-y_0)\not \in \Om_-(t)\, .%
\ee  %

If $2s_0-y_0<y(x'_0,t)$, then (\ref{formula 2 of lemma 5.2})
contradicts the assumption $P_{s}(\Om)=\Om$ for $s=s_0<t$.

Let $2s_0-y_0\ge y(x'_0,t)$. Since $P_s(\Om)=\Om$ for all $s\le t$
and $x_0\in \Om$, it follows that the closed interval
$\{(x'_0,y):\,y_0\le y\le 2t-s_0\}$ is in $\Om$. Then %
$$  %
m_\Om(x'_0,y(x'_0,t))>y(x'_0,t)-y_0\ge 2(y(x'_0,t)-s_0)\ge
2(y(x'_0,t)-t)  %
$$  %
contradicting the equation~(\ref{eq. for l.f.}). The proof of
Lemma~\ref{Criterion of partial symmetry for functions}(a) for
open sets is complete. If $\Om$ is not open the proof is left to
the reader.

The proof of Lemma~\ref{Criterion of partial symmetry for
functions}(b) easily follows from Lemma~\ref{Criterion of partial
symmetry for functions}(a) and Definition~\ref{SCSym of
functions}. \hfill  $\Box$

Our next result will be used in the proof of
Theorem~\ref{11-3-Theorem} in Section~11.

\bc \label{Corollary 6.3} %
Let $\Omega$ be a bounded domain in $\mathbb{R}^n$ and let $f\in
L_+^2(\Omega)$, $f>0$ on $\Omega$. If for some $T\in \mathbb{R}$,
$f^T$ is not a translation of $f$ in the direction of $y$-axis,
then we can find a polarizer $H_{t_0}=\{y>t_0\}$ with $t_0<T$ such
that $\left(f_{H_{t_0}}\right)^T=f^T$, $f\not=f_{H_{t_0}}$, and
$\sigma_{H_{t_0}}(f)\not=f_{H_{t_0}}$.

In particular, if $\Omega^{T}$ is not a translation of $\Omega$ in
the direction of $y$-axis, then we can find a polarizer
$H_{t_0}=\{y>t_0\}$ with $t_0<T$ such that
$\left(\Omega_{H_{t_0}}\right)^T=\Omega^T$, $\Omega\not=
\Omega_{H_{t_0}}$, and
$\sigma_{H_{t_0}}(\Omega)\not=\Omega_{H_{t_0}}$.
\ec %

{\it Proof.} Let $t_1=\inf\,t$, where the infimum is taken over
all $t\in \mathbb{R}$ such that $\Omega_+(t)=\emptyset$. Then for
all $t\ge t_1$, $f^t$ is a translation of $f^{t_1}$. Therefore,
proving this corollary, we may assume that $T\le t_1$. In this
case, $f^T$ is not a translation of $f$ if and only if
$f^T\not=f$.

Let $t_\Omega=\sup t$, where the supremum is taken over all $t\in
\mathbb{R}$ such that $P_\tau(f)=f$ for all $\tau\le t$. Then,
$f^{t_\Omega}=f$ by Lemma~\ref{Criterion of partial symmetry for
functions}. Since $f^T\not=f$,
it follows from the same lemma %Lemma~\ref{Criterion of partial symmetry for functions}(a)
that $t_\Omega<T$. Combining this with Corollary~\ref{Stability of
functions for polarization}, we obtain that
$$ %
\left(f_{H_s}\right)^T=\left(\left(f^{t_\Omega}\right)_{H_s}\right)^T=f^T
\quad {\mbox{for all $s\le (1/2)(t_\Omega+T)$.}}
$$ %
 Since $t_\Omega<T$,
$f^{t_\Omega}=f$, and $\Omega_+(t_\Omega)\not= \emptyset$, we must
have $\sigma_{H_s}(\Omega)\not=\Omega_{H_s}$, and therefore
$\sigma_{H_s}(f)\not=f_{H_s}$, for all $s$ sufficiently close to
$t_\Omega$. Therefore, every polarizer $H_s$ with
$t_\Omega<s<(1/2)(t_\Omega+T)$ such that $f_{H_s}\not= f$ will satisfy the requirements of the corollary.  %

 By the definition of $t_\Omega$, there is a sequence $\{s_k\}_{k=1}^\infty$, which converges to
 $t_\Omega$ from the right, such that $f_{H_{s_k}}\not= f$ for all $k$.
 Choosing $k$ sufficiently large such that $s_k\le (1/2)(t_\Omega+T)$, we obtain the desired polarizer $H_{t_0}=H_{s_k}$.
This  proves the first assertion of the corollary. Since
$\Omega=\{f>0\}$ the second assertion is a particular case of the
first one. \hfill $\Box$

\begin{lemma}
Let $\Om\in {\mathcal{M}}_n$ and let $(\cdot)^*$ denote the
$(1,n)$-Steiner symmetrization with respect to the plane
$\{y=t\}$. If  $T^t (\Omega ) =\Omega ^*$ and $\Om\in
{\mathcal{F}}\cup{\mathcal{G}}_b$ then $T^t (\Omega _{\varepsilon}
)=(\Omega _{\varepsilon })^*$ for every $\varepsilon >0$.
 \end{lemma}

{\it Proof.}   Assume that $\Omega \in {\mathcal{G}}_b $. For a
fixed $\ep>0$, let $y_{\varepsilon}(x',t)$ and
$m_{\varepsilon}(x',y)$ denote, respectively,  the separating
function  and measuring
function of $\Omega_\ep $. For $x'\in (\Om_\ep)'$, let  %
$$  %
 y_{\varepsilon}(x')=\sup\{y:\ (x',y)\in T^t(\Omega_\varepsilon
)\}.   %
$$   %
Then $y_\ep(x')\ge y_\ep(x',t)$.  If $
y_{\varepsilon}(x')=y_{\varepsilon}(x',t)$ for all $x' \in (
\Omega_\ep)' $,  the conclusion of the lemma is immediate.

Suppose that  %
$$  %
 y_{\varepsilon}(x'_0)> y_{\varepsilon} (x'_0,t)  %
$$  %
 for some $x'_0\in (\Om_\ep)'$. Then (\ref{eq. for l.f.}) implies %
\be  \label{formula 1 of Lemma 3} %
 m_{\varepsilon}(x'_0,y_{\varepsilon}(x'_0))<
2(y_{\varepsilon}(x'_0)-t).   %
\ee  %
Now we fix $\delta>0$ sufficiently small. Then we can find two
points $\hat x=(x'_0,\hat y)$ and $\tilde x=(\tilde x',\tilde y)$
satisfying the following conditions: %
\be \label{formula 2 of Lemma 3} %
\hat x\in T^t(\Omega_\ep) \qquad {\rm and} \quad
y_{\varepsilon}(x'_0)-\delta<\hat y< y_{\varepsilon}(x'_0),  %
\ee  %
\be \label{formula 3 of Lemma 3} %
\tilde x\in \Omega \qquad {\rm and }\qquad |\hat x-\tilde
x|<\varepsilon. %
\ee  %
 Since $T^t(\Omega)=\Omega^*$ and $\tilde x\in \Om$ we have:
\be \label{formula 4 of Lemma 3}  %
m(\tilde x',\tilde y)\ge 2(\tilde y-t).  %
\ee  %
Using (\ref{formula 3 of Lemma 3}) and (\ref{formula 4 of Lemma
3}), after an elementary geometric argument we obtain the inequality  %
$$  %
 m_{\varepsilon}(x'_0,\hat y)\ge 2(\hat y-t),  %
$$  %
which together with (\ref{formula 2 of Lemma 3}) implies %
$$  %
 m_{\varepsilon}(x'_0,y_{\varepsilon}(x'_0))\ge
m_{\varepsilon}(x'_0,\hat y)\ge
2(y_{\varepsilon}(x'_0)-t)-2\delta.  %
$$  %
Since $\delta>0$ can be chosen arbitrary small, the latter
inequality contradicts (\ref{formula 1 of Lemma 3}).  This proves
the lemma for open sets. If $\Om$ is a compact set the proof
follows the same lines.  $\hfill \Box$

\section{Approximation by  polarizations}  \label{Approximation by  polarizations}%
\setcounter{equation}{0} %

A compact set $K\subset \mathbb{R}^n$ is said to be {\it simple}
if it can be decomposed into a finite number of blocks $\bar
R_i=\bar R(a^i,b^i)$ such that $R_i\cap R_j=\emptyset$ if $i\not=
j$. Here $\bar R(a,b)$ denotes the closure of an open block  %
$$  %
R(a,b):=\{x\in \mathbb{R}^n:\  b_k<x_k<b_k+a_k, 1\le k\le n\}, %
$$  %
where  %
$a=(a_1,\ldots,a_n)\in\mathbb{R}^n_+$ and $b=(b_1,\ldots,b_n)\in
\mathbb{R}^n$.

If $K$ is a simple compact set in
$\mathbb{R}^n=\mathbb{R}^{n-1}\times \mathbb{R}$, $n\ge 2$, then
its projection $K'$ is a simple compact set in $\mathbb{R}^{n-1}$.
Let $K'=\cup_{i=1}^N {\bar R}'_i$ be a block decomposition of
$K'$.

A collection of open sets $Cu(K):=\left\{C_i=R'_k\times
\mathbb{R}\right\}_1^N$ will be called \,{\it a cubicle structure
of} $K$. If $C_i$ is a cubicle in $Cu(K)$, then $C_i\cap
K=R_i'\times K^i$, where $K^i$ is a simple compact set in
$\mathbb{R}$. Therefore $C_i\cap K$ consists of a finite number of
disjoint blocks $K_{i,j}=R'_i\times[\al_{ij},\beta_{ij}]$, $1\le
j\le j_i$. We always enumerate these blocks such that %
$$  %
\al_{i1}<\beta_{i1}<\al_{i2}<\cdots <\al_{ij_i}<\beta_{ij_i},
\qquad 1\le i\le N.  %
$$  %
  Thus every simple compact set $K$ admits
{\it a block decomposition} of the form %
$$  %
K=\bigcup _{i=1}^N \overline{R'}_i \times
K^i=\bigcup_{i,j}\overline{K}_{ij}.  %
$$  %

 In this section we  show that the $T^t$-transformation of sets and functions can be approximated by a sequences of
 polarizations. This approach was first used by V.~Wolontis \cite{W}. The method was developed
 further by V.~N.~Dubinin \cite{D2}, \cite{D3}. So we will call this method {\it the Wolontis-Dubinin approach}.
 In the context of continuous symmetrization the Wolontis-Dubinin approach was first used in
 \cite{S1}. %
 %Figure~\ref{5-figure} illustrates some notation used in the proof of Lemma~\ref{approximation of compact sets}, which are related to to this approach.

\bl  \label{approximation of compact sets}  %
For every $t\in \mathbb{R}$ and every collection of simple compact
sets $K_1,\ldots,K_m$ in $\mathbb{R}^n$ there exist a finite
number of polarizations
$P_k$, $1\le k\le N$ with polarizers $\{ y>y_k \}$ such that  %
\be \label{formula 1 of lemma 6.1}  %
y_1 <y_2 <\ldots <y_N<t  %
\ee  %
and %
\be  \label{formula 2 of lemma 6.1}  %
 \bigcirc _{k=1}^N P_k (K_i
)=K_i ^t  \qquad {\mbox{for all \ $1\le i\le m$.}} %
\ee %
\el  %

{\it Proof}. %
(1) \ Let $t\in \mathbb{R}$ be fixed.  First we consider the  case
$n=1$. %To unify our further notation we will write $K_i^{(0)}$ for $K_i$.
Then the block decomposition of  $K_i$ consists of a finite number
of
disjoint segments:  %
\be \label{formula 3 of lemma 6.1}  %
K_i =\bigcup_{j=1}^{j_i} \ [a_{ij}, b_{ij}],  %
\ee   %
 where  %
\be  \label{formula 4 of Lemma 6.1}   %
a_{i1} < b_{i1} <a_{i2} < b_{i2}
<\cdots <a_{ij_i} < b_{ij_i} ,\qquad j_i\ge 1, \qquad 1\le i\le m.  %
\ee   %

\smallskip

For every index $i$ we define two numbers:  %
\be   \label{formula 5 of lemma 6.1} %
\alpha _i = (1/2) \left( t+ (1/2)\left(a_{i1}
+b_{i1}\right)\right) ,  %
\ee  %
and   %
\be  \label{formula 6 of lemma 6.1}  %
\beta _i = \left\{ \begin{array}{ll}  %
(1/2) \left( a_{i1} +a_{i2} \right)  & \qquad{\mbox{if $j_i
>1$}} \\ %
t  & \qquad {\mbox{otherwise.}} %
\end{array} \right. %
\ee  %
 Let ${\mathcal{A}}_t^{(0)}$ denote the set of all numbers
$\alpha _i$ and all numbers  $\beta _i$ that are less than $t$. It
is easy to verify that $T^t\left((K_i)\right)=K_i$ for some  $i$
if and only if $\al_i\ge t$. Therefore if
${\mathcal{A}}_t^{(0)}=\emptyset$, then $\al_i\ge t$ for all $i$
and there is nothing to prove.

\smallskip

If ${\mathcal{A}}_t^{(0)}\not=\emptyset$, then we set %
$$  %
y_1 := \min \{ \gamma :\ \gamma \in {\mathcal{A}}_t^{(0)} \}.  %
$$  %

Now we consider a polarization  $P_1$ with the polarizer $
\{y>y_1\}$. Let $ K_i^{(1)} = P_1 (K_i)$,  $1\le i\le m$. For each
$i$, the compact set $K_i^{(1)}$ satisfies the following
conditions:  %
\begin{enumerate}  %
\item[(a)] %
$K_i^{(1)}$ is a  simple compact set whose block decomposition
$K_i^{(1)}=$ \\ $\bigcup_{j=1}^{j_i^{(1)}} [a_{ij}^{(1)},
b_{ij}^{(1)}]$ consists of $j_i^{(1)}$ disjoint segments, where
$j_i^{(1)}\le
j_i$.   %
\item[(b)] %
 $\left(K_i^{(1)}\right)_+ (t)=\left(K_i\right)_+ (t)$ and therefore $T^t(K_i^{(1)})=T^t(K_i)$  for all $1\le i\le
 m$. %
\end{enumerate}  %

\smallskip

Now starting with simple compact sets  $K_i^{(1)}$, $1\le i\le m$,
we define a set ${\mathcal{A}}_t^{(1)}$ of elements
$\al_i^{(1)}<t$ and $\beta_i^{(1)}<t$ in the same way as
${\mathcal{A}}_t^{(0)}$ was defined for $K_i$, see formulas
(\ref{formula 3 of lemma 6.1})--(\ref{formula 6 of lemma 6.1}). As
we noticed above, if ${\mathcal{A}}_t^{(1)}=\emptyset$, then
$K_i^{(1)}=T^t(K_i^{(1)})=T^t(K_i)$, $1\le i\le m$, where the
second equality follows from the condition~(b). In this case our
construction is finished.

If ${\mathcal{A}}_t^{(1)}\not=\emptyset$, we continue our
construction further to get simple compact sets $K_i^{(2)}$,
$K_i^{(3)}$, $\ldots$\, , $1\le i\le m$. Let
${\mathcal{A}}_t^{(k)}$, $\al_i^{(k)}$, $\beta_i^{(k)}$, and
$y_{k+1}$ correspond $K_i^{(k)}$ in the process of this
construction. One can easily see that in each step of our
construction $K_i^{(k)}$, ${\mathcal{A}}_t^{(k)}$, and $y_{k+1}$
satisfy the following
conditions:  %
\begin{enumerate}  %
\item[(i)]  %
 If for some $k\ge 1$, ${\mathcal{A}}_t^{(k)}\not=\emptyset$, then $y_1<y_2<\cdots <y_k<t$.  %
 \item[(ii)]  %
 If for some $k\ge 1$ and some $s$, $1\le s\le m$, $y_k=\al_s^{(k-1)}$, then $\#{\mathcal{A}}_t^{(k)}
 \le \#{\mathcal{A}}_t^{(k-1)}-1$, where $\#$ denotes the cardinality of the corresponding set and
 ${\mathcal{A}}_t^{(0)}={\mathcal{A}}_t$.  %
 \item[(iii)] %
 If $y_k=\beta_s^{(k-1)}$ for some $s$, $1\le s\le m$, then $j_s^{(k)}\le j_s^{(k-1)}-1$,
 where $j_i^{(l)}$ denotes the number of disjoint segments in $K_i^{(l)}$ if $l\ge 1$ and
 $j_i^{(0)}=j_i$. %
 \end{enumerate}  %

Since the total number of segments constituting $K_1,\ldots,K_m$
is finite and $\#{\mathcal{A}}_t<\infty$ the conditions (ii) and
(iii) show that the above construction of compact sets $K_i^{(k)}$
terminates after a finite number of steps. Thus there is an
integer $N\ge 0$ such that ${\mathcal{A}}_t^{(N)}=\emptyset$.
Since in each step of our construction  the conditions (a) and (b)
are satisfied we have %
$$  %
K_i^{(N)}=T^t\left(K_i^{(N)}\right)=T^t\left(K_i^{(N-1)}\right)=\cdots=T^t\left(K_i\right),
\qquad 1\le i\le m,  %
$$  %
which proves the lemma  in the case under consideration.

\medskip

(2) \  Now we consider the case $n\geq 2$. Since $K_i$ is a simple
compact set in $\mathbb{R}^n$ it can be  decomposed as follows %
$$ %
K_i =\bigcup _{j=1}^{j_i} \overline{R'}_{ij} \times K_i^j ,  %
$$  %
where $\overline{R'}_{ij}$ are closed blocks in $\mathbb{R}^{n-1}$
such that $R'_{ij}\cap R'_{il}=\emptyset$ if $j\not=l$  and
$K_i^j$ are simple compact sets in ${\mathbb{R}}$. %

By part~(1) of this proof, for every $t\in \mathbb{R}$ there are
polarizations $P_k$, $1\le k\le N$,  in $\mathbb{R}$ with the polarizers $\{y>y_k\}$ such that %
$$  %
y_1<y_2<\cdots <y_N<t  %
$$  %
and  %
$$  %
\bigcirc _{k=1}^N P_k (K_i^j )=T^t(K_i^j ) \qquad {\mbox{for all
$1\le j\le j_i$ and $1\le i\le m$.}}  %
$$  %
Considering $P_k$, $1\le k\le N$, as polarizations in
$\mathbb{R}^n$ with the corresponding polarizers $\{y>y_k\}$ we
obtain   %
$$
\bigcirc _{k=1}^N P_k (K_i )= \bigcup_{j=1}^{j_i}
\overline{R'}_{ij} \times \left(T^t(K_i^j)\right)
=T^t\left(K_i\right) ,\qquad 1\le i\le m.  %
$$ %

The proof of the lemma is complete.  \hfill $\Box $

%\begin{figure}
%$$\includegraphics[scale=.3,angle=0]{iceberg} $$
%\caption{Approximation by polarizations.} \label{5-figure}
%\end{figure}

\bc  \label{approximation of open sets by polarization}  %
 Let
$t\in {\mathbb{R}}$ and let $D$, $\Omega$ be bounded  open  sets
such that $\bar D \subset\Omega$. Then there is a finite number of
polarizations $P_k $,  $1\le  k  \le N$, with polarizers
$\{y>y_k\}$, where $y_1<y_2<\ldots<y_N\le t$, such that %
$$  %
\bar{D}_N\subset \Omega^t,  %
\qquad {\mbox{where $\bar{D}_N=\bigcirc _{k=1} ^N P_k(\bar{D})$.}}  %
$$  %

In addition,  %
$$   %
 {\rm {dist}}(\partial D,\partial \Omega)\le
{\rm{dist}}(\partial D_N,\partial\Omega^t), \qquad {\mbox{where $D_N=\bigcirc _{k=1} ^N P_k(D)$.}}  %
$$  %
\ec   %

{\it Proof}. Let $\rho_0={\rm{dist}}(\partial D,\partial\Omega)$.
 Then $\rho_0>0$ and for every $\rho$,  $0\le\rho<\rho_0$, there exists a simple
compact set $K=K(\rho)$ such that  %
$$  %
(\bar{D})_\rho\subset K \subset \Omega .  %
$$  %
 By Lemma~\ref{approximation of compact sets},  we can find a finite number of polarizations
$P_k$, $1\le k\le N$ such that  %
$$   %
 K^t=K_N:= \bigcirc_{k=1}^N P_k(K).  %
$$  %

Since the $T^t$-transformation and polarization both are monotone
and smoothing transformations we have  %
\be \label{formula 2 of corollary 3}  %
 \Omega^t\supset K^t=K_N\supset \bigcirc_{k=1}^N P_k((\bar{D})_\rho) \supset
 (\bar{D}_N)\rho. %
\ee  %
 Since  (\ref{formula 2 of corollary 3}) holds true for every
 $\rho$, $0\le \rho < \rho_0$, the corollary follows.  %
 \hfill $\Box$

\smallskip

\bd \label{Simple function}  %
A function $f\in {\mathcal{S}}_n$ is said to be  simple if $f$ can
be represented as %
\be \label{formula of a simple function}  %
 f=\varepsilon \left( -n_0
+\sum\limits_{i=1}^m \chi (K_i )\right)    %
\ee  %
with some $\varepsilon
>0$, $n_0\in {\mathbb{N}} $ and some simple compact sets $K_i$,  $1\le i\le
m$, such that $K_1 \supset \cdots \supset K_m $. %
\ed  %

We note that simple functions are dense in  $L^p _+
({\mathbb{R}}^n )$ for every $1\le p<\infty $. The following lemma
shows that the $T^t$-transformation of a simple function can
be reduced to a finite sequence of polarizations.  %

\bl \label{approximation of a simple function by polarization}  %
Let $f:\mathbb{R}^n\to \mathbb{R}$ be a  simple function. Then for
every $t\in \mathbb{R}$  there are polarizations $P_k $ with
polarizers  $\{ y>y_k \}$,  $1\le k\le N$ such that %
\be \label{formula 1 of Lemma 6.2}  %
 y_1 <y_2<\ldots <y_N<t,  %
\ee  %
 and  %
\be  \label{formula 2 of lemma 6.2} %
 \bigcirc _{k=1}^N P_k (f)=f^t .  %
\ee  %
\el  %

{\it Proof}. Since  $f$ is simple it can be represented in the
form~(\ref{formula of a simple function}).   By
Lemma~\ref{approximation of compact sets}, we can find
polarizations $P_k$,  $1\le k\le N$, with the polarizers
$\{y>y_k\}$ satisfying (\ref{formula 1 of Lemma 6.2})  such that
$$  %
\bigcirc _{k=1}^N P_k (K_i )=K_i^t , \qquad 1\le i\le m.  %
$$  %
Now (\ref{formula 2 of lemma 6.2})  follows from the relations
$$  %
\bigcirc _{k=1}^N P_k (f) =\varepsilon \left( -n_0+
\sum\limits_{i=1}^m \chi \left( \bigcirc _{k=1}^N P_k (K_i
)\right) \right) %
$$  %
and  %
$$  %
f^t =\varepsilon \left( -n_0 +\sum\limits_{i=1}^m \chi (K_i ^t
)\right),   %
$$  %
which in their turn follow from Definition~\ref{SCSym of
functions} and Definition~\ref{Defenition of polarization of function}.  %
 \hfill   $\Box $  %

 \smallskip

 The following lemma shows that, for functions $u\in
 L^p_+(\mathbb{R}^n)$ such that $u\not=u^t$, one polarization is
 enough to find a better approximation of $u^t$. Its proof, left to
 the reader, follows the ideas of \cite[p.252]{BaT} and repeats
 the proof of Lemma~6.4 in \cite{BS}. %
 \bl  \label{Lemma - Better approximation for p-norm}  %
 Let $u\in L^p_+(\mathbb{R}^n)$, $1\le p<\infty$. If $u=u^{t_1}$ and $u\not=u^t$ for some $t_1,t\in {\mathbb{R}}$ such that $t_1<t$,
 then there is a polarizer $H_s=\{y>s\}$ with $t_1<s<t$ such that      %
$$  %
\left(u_{H_s}\right)^t=u^t  %
$$  %
and %
$$  %
\|u_{H_s}-u^t\|_p<\|u-u^t\|_p.  %
$$  %
\el  %

\bl  \label{Lemma - Approximating sequence in p-norm} %
Let $u\in L^p_+(\mathbb{R}^n)$, $1\le p<\infty$. Then for every
$t\in \mathbb{R}$ there is a sequence of polarizations $P_k$ with
polarizers $\{y>y_k\}$, where $y_k\le t$ for all $k=1,2,\ldots$,
such that the sequence of functions $u_m:=\bigcirc_{i=1}^m P_i u$
satisfies the minimality condition  %
\be  \label{equation-7.11}  %
\|u_{m+1}-u^t\|_p=\min \{\|(u_m)_H-u^t\|_p\}, \quad m=1,2,\ldots,
\
u_0=u, %
\ee  %
where the minimum is taken over all polarizations with polarizers
$H=\{y>s\}$ such that $-\infty<s\le t$. %

In addition, %
$$  %
u_m\to u^t \quad {\mbox{in $L^p(\mathbb{R}^n)$.}}  %
$$  %
\el  %

{\it Proof}. It follows from Lemma~5.2 \cite{BS} that the minimum
in (\ref{equation-7.11}) is attained for some polarizer
$\{y>y_{k+1}\}$, $-\infty<y_{k+1}\le t$. Then by Lemma~6.1
\cite{BS}, there are some function $v\in L^p_+(\mathbb{R}^n)$ and
a
subsequence $u_{m'}$ such that  %
$$  %
u_{m'}\to v \quad {\mbox{in $L^p(\mathbb{R}^n)$.}}  %
$$  %
Then using the non-expansivity Lemma~\ref{ Non-expansivity of
rearrangements} we obtain $v^t=u^t$.

Now assume that $v\not=u^t$. By Lemma~\ref{Lemma - Better
approximation for p-norm}, we can choose a polarizer $H$ such that %
$$  %
\|v_H-u^t\|_p<\|v-u^t\|_p.  %
$$  %
It follows that %
\begin{alignat}{10}  \label{equation-7.12} %
\|(u_{m'})_H-u^t\|_p-\|u_{m'}-u^t\|_p &\le &
\|(u_{m'})_H-v_H\|_p+\|v_H-u^t\|_p \ \ \ \ \ \ \ \ \ \ \
\\   %
+\|u_{m'}-v\|_p-\|v-u^t\|_p&\le&
2\|u_{m'}-v\|_p+\|v_H-u^t\|_p-\|v-u^t\|_p  \nonumber \\%
&\to& \|v_H-u^t\|_p-\|v-u^t\|_p<0 \nonumber \ \ \ \ \ \ \ \ \ \ \
\ \ \
\end{alignat}  %
as $m'\to \infty$.

On the other hand the sequence $\|u_m-u^t\|_p$ is monotonically
decreasing. Hence,  %
$$  %
\|u_{m+1}-u^t\|_p-\|u_m-u^t\|_p\to 0 \quad {\mbox{as $m\to
\infty$.}} %
$$  %
Together with (\ref{equation-7.12}) this contradicts the
minimality property~(\ref{equation-7.11}).  \hfill $\Box$

\smallskip

\bl \label{Lemma - convergence for continuous functions}  %
Let $u\in C(\mathbb{R}^n)\cap L^1(\mathbb{R}^n)$ and let $u_m$ be
defined as in Lemma~\ref{Lemma - Approximating sequence in
p-norm}, whereby the condition (\ref{equation-7.11}) is satisfied
with $p=1$. Then %
$$  %
u_m\to u^t \quad {\mbox{in $C(\mathbb{R}^n)$.}} %
$$  %
\el  %

{\it Proof}. By Lemma~\ref{Lemma - Approximating sequence in
p-norm}, $u_m\to u^t$ in $L^1(\mathbb{R}^n)$, and the functions
$u_m$ are equicontinuous in view of Lemma~5.1 \cite{BS}. Because
of Lemma~6.2 \cite{BS} we also have %
$$  %
u_{m'}\to v \quad {\mbox{in $C(\mathbb{R}^n)$}} %
$$  %
for a subsequence $u_{m'}$ and $\omega_v\le \omega_u$. Thus
$v=u^t$ and the assertion follows. \hfill $\Box$

\section{Rescaling and limit cases} \label{Rescaling and limit cases} %
\setcounter{equation}{0}  %

The $T^t$-transformation of sets and functions defined in
Section~\ref{ Continuous $(1,n)$-Steiner symmetrization} depends
on the height $t$ of the moving plane of symmetrization $\{y=t\}$,
the range of which   is $-\infty<t<\infty$. We recall  that  in
its original setting as described in the Introduction the problem
on the continuous symmetrization requires the standard
range $0\le t\le 1$ with the ``boundary conditions'':  %
\be \label{formula 1 of Section Rescaling} %
T^0(\Om)=\Om, \qquad T^1(\Om)=\Om^*,  %
\ee  %
where $\Om^*$ denotes $(1,n)$-Steiner symmetrization of $\Om$ with
respect to a fixed plane $\{y=t_0\}$. If $t_0=1$, then the
$T^t$-transformation perfectly matches these conditions when
 restricted to the family ${\mathcal{M}}_n[0,1]$, where
$$  %
{\mathcal{M}}_n[a,b]:=\{\Om\in {\mathcal{M}}_n:\, \Om\subset
\Pi(a,b)\}  %
$$  %
and  %
 $$  %
 \Pi(a,b):=\{(x,y)\in \mathbb{R}^n:\, a\le y\le b\}, \qquad
   \quad -\infty \le a<b\le \infty.  %
 $$  %

  If $\Om\in {\mathcal{M}}_n[a,b]$ with $-\infty <a<b<\infty$, the
  linear change of variables $\tau=(t-a)/(b-a)$ gives the desired
  range $0\le \tau\le 1$ with ${\tilde T}^0(\Om)=\Om$ and ${\tilde
  T}^1(\Om)=\Om^*$, where ${\tilde T}^\tau=T^t$ and the asterisk  denotes
  the $(1,n)$-Steiner symmetrization with respect to $\{y=b\}$.
   The case $\Om\in {\mathcal{M}}_n[-\infty,b]$ with $b<\infty$
   can be handled in a similar way with  a non-linear change of variables, for instance,
   $\tau=(1+b-t)^{-1}$.

   The remaining case $\Om\in {\mathcal{M}}_n[a,\infty]$ with
   $-\infty \le a<\infty$ requires some additional consideration
   since the limit set $\Om^\infty$ is
   not defined yet. One possible way to ``correct'' this defect is to
   use a post-composition ${\widehat{T}}^t=Sh^t\circ T^t$ of the $T^t$-transformation with a
   continuous shift $Sh^t$ of $(t-t_0)^+$ units in the direction of the negative
   $y$-axis. Here  $(t-t_0)^+=\max \{0,t-t_0\}$.
   In this case, the limit set
   $\Om^\infty={\widehat{T}}^\infty(\Om)$ may be identified with
    the Steiner symmetrization $\Om^*$ of $\Om$ with respect to
   the plane $\{y=t_0\}$.  Of course, ${\widehat{T}}^t$ still
possesses all basic properties of the $T^t$-transformation.
   In particular, Lemmas~\ref{Basic geometric properties of SCSym}  and \ref{Continuity of T^t} are still
   valid for ${\widehat{T}}^t$ for all $-\infty <t\le +\infty$.
   The case $t=\infty$ is included here because
   ${\widehat{T}}^\infty$ is just the Steiner symmetrization with
   respect to $\{y=t_0\}$.
   Changing the variables via $\tau=(1/\pi)\arctan t +1/2$ we will
   have a family of transformations ${\widetilde{T}}^{\tau}$
   defined on the standard interval $I=\{\tau:\,0\le \tau \le
   1\}$. As we have seen in previous sections, the $T^t$-transformations is continuous from the left for $-\infty <t<\infty$. Therefore
   the ${\widetilde{T}}^t$-transformation is continuous from the left on the open interval $0<\tau<1$.
By this reason, the ${\widetilde{T}}^t$-transformation will be
called the continuous symmetrization defined on the standard
interval.
  % {\it In what follows the symbol $\ \tilde{}$  (tilde) will always denote
%   a continuous transformation defined on the standard
%   interval.}

Being restricted to the class of bounded sets and functions with
bounded supports, the ${\widetilde{T}}^t$-transformation is
continuous at $\tau=1$ as well.
   For unbounded sets and functions with unbounded supports the
   situation is different: some functionals of domains or functions are
   continuous at $\tau=1$ while some other such functionals are not
   continuous at $\tau=1$. For instance,    Lemma~\ref{Continuity of T^t} is obviously not
   valid if $\Omega$ is an unbounded set while the relevant    Steiner symmetrization of $\Omega$ is bounded.

\medskip

To define a one-dimensional $T^t$-transformation for an arbitrary
half-space $H=H(a,n)$, we may combine the $T^t$-transformation of
Definition~\ref{SCSym} with appropriate affine transformations of
$\mathbb{R}^n$. To be more precise, let $k>0$ and let $A$ be an
orthogonal $n\times n$ matrix such that the operator $L:=k(A-Aa)$
maps $H(a,n)$ onto $\{y>0\}$. Then $T^t_H:=L^{-1}\circ T^t\circ L$
is the desired continuous symmetrization with respect to $H$. If
$M\in {\mathcal{M}}_n$ is bounded then the scaling factor $k$ can
be chosen such that $L(M)\subset B^{(n)}$ and therefore when
working with bounded sets $M$ we may assume without loss of
generality that $M\subset B^{(n)}$. Similarly if $u\in
{\mathcal{S}}_n$ has a bounded support we often assume without
loss of generality that $\supp u \subset B^{(n)}$.

\medskip

We stress once more that {\it any $T^t_H$ -transformation can be
rescaled to get a ${\widetilde{T}}^t_H$-transformation defined on
the standard interval $I$ and that all basic properties of
continuous symmetrization discussed in Sections~\ref{ Continuous
$(1,n)$-Steiner symmetrization}-\ref{Approximation by
polarizations} still remain true for any
${\widetilde{T}}^t_H$-transformation.}

\section{ Continuous $(k,n)$-Steiner symmetrizations} \label{ Continuous $(k,n)$-Steiner symmetrizations}%
\setcounter{equation}{0}  %

To define a one-parameter family of rearrangements $T^t_k$ left
continuous in $t\in [0,1]$ transforming sets and functions into
their $(k,n)$-Steiner symmetrization for any $2\le k\le n$, we
suggest the following inductive algorithm:

\medskip

Let $2\le k\le n$. Suppose that for any $(n-k+1)$-dimensional
plane $L'$ we have defined a continuous family of rearrangements
transforming sets, functions, etc. into their $(k-1,n)$-Steiner
symmetrizations with respect to $L'$. Then let $S$ be a
$(k,n)$-Steiner symmetrization with respect to an
$(n-k)$-dimensional plane $L$. By Lemma~\ref{Sarvas approximation
lemma}, we can choose two intersecting $(n-k+1)$-dimensional
planes $L_1$ and $L_2$ defining
$(k-1,n)$-Steiner symmetrizations $S_1$ and $S_2$ such that %
\be  \label{9.1-equation}  %
S=\lim_{i\to \infty} \left(S_2\circ S_1\right)^i=\lim_{i\to
\infty} S_1\circ \left(S_2\circ S_1\right)^i,  %
\ee  %
where $\left(S_2\circ S_1\right)^0$ is the identity
transformation.

By our inductive assumption there are continuous $(k-1,n)$-Steiner
symmetrizations $T^{1,t}$ and $T^{2,t}$, $t\in [0,1]$ into $S_1$
and $S_2$, respectively. Let $t(j) =(j-1)/j$ and let
$I_j=\{t:\,t(j)\le t< t(j+1)\}$, $j=1,2,\ldots$  Then
$\left(\cup_{j=1}^\infty I_j\right) \cup \{1\}$ is a disjoint
decomposition of the standard  interval $I$. % =\{t:\,0\le t\le 1\}$. into a disjoint union of sets. %
Next, for every $j=1,2,\ldots$
and every $t\in I_j$ we define the transformation $T^t_k$ by %
\be \label{9.2-equation}  %
T^t_k=\left\{ %
\begin{array}{lllll}  %
T^{1,\tau} \circ \left(S_2\circ S_1\right)^{i-1}& {\mbox{with}}&
\tau=\frac{t-t(2i-1)}{t(2i)-t(2i-1)} \ &{\mbox{if\ }} j=2i-1,
\\%
T^{2,\tau} \circ S_1\circ \left(S_2\circ S_1\right)^{i-1}&
{\mbox{with}}& \tau=\frac{t-t(2i)}{t(2i+1)-t(2i)} \
&{\mbox{if\ }} j=2i,  %
\end{array}\right.  %
\ee  %
where $i=1,2,\ldots$

\smallskip

Inductive algorithm (\ref{9.2-equation}) defines a continuous
transformation into a $(k,n)$-Steiner symmetrization $S$ as
follows. First, we choose a complete binary tree ${\mathcal{T}}$
of planes $L_{0i_1\ldots i_l}$, $1\le l\le k-1$, $i_s\in\{1,2\}$,
rooted at $L=L_0$  such that for every $1\le l\le k-2$ and every
multi-index $i_1\ldots i_l$, the planes $L_{0i_1\ldots i_l1}$ and
$L_{0i_1\ldots i_l2}$ are $(n-k+l+1)$-dimensional and the
corresponding $(k-l-1,n)$-Steiner symmetrizations  $S_{0i_1\ldots
i_l1}$ and $S_{0i_1\ldots i_l2}$ approximate $S_{0i_1\ldots i_l}$
as in (\ref{9.1-equation}).

\smallskip

\bd  \label{continuous $(k,n)$-symmetrization}  %
By SC $(k,n)$-Steiner symmetrization corresponding to a tree
${\mathcal{T}}$ we mean a transformation $T^t_k$ of sets and
functions defined inductively by (\ref{9.2-equation}), where the
initial $2^{k-1}$ continuous $1$-dimensional symmetrizations are
selected to be  SC $1$- symmetrizations with respect to the
corresponding hyperplanes $L_{0i_1\ldots
i_{k-1}}$,  where $i_1\ldots i_{k-1}\in  \{1,2\}^{k-1}$. %
\ed  %

%\smallskip

%Figure~\ref{6-figure} shows the effect of the SC $(1,2)$-Steiner symmetrization on some simple shapes.
Three remarks are in order now: %
 %
 %\brem %

\noindent %
{\bf Remark~9.1.} \label{9.1-Remark} %
The Sarvas approximation scheme used in the
algorithm~(\ref{9.2-equation}) may be replaced by any other
approximation scheme. For example, one may use a sequence $S_i$,
$i=1,2,\ldots$ of $(1,n)$- Steiner symmetrizations to approximate
a given $(k,n)$-Steiner symmetrization: $S=\lim_{n\to \infty}
\bigcirc_{i=1}^n S_i$. In this case a $T^t$-transformation can be
defined as $T^t=T_n^\tau \circ \left(\bigcirc_{i=1}^n S_i\right)$
with $\tau=(t-t(i))/(t(i+1)- t(i))$. %
%\erem  %

%\brem \label{9.2-Remark} %

\smallskip

\noindent  %
{\bf Remark~9.2.} \label{9.2-Remark} %
One can define a continuous $(k,n)$-symmetrization by replacing
Solynin's one dimensional symmetrization with Brock's one
dimensional symmetrization. Alternatively, one may want to define
a continuous $(k,n)$-symmetrization by ``mixing'' in the different
stages of the algorithm~(\ref{9.2-equation}) the SC
symmetrizations with  BC symmetrizations (or with other types of
continuous $(1,n)$-symmetrization if those will be discovered).
%\erem  %

%\brem  \label{9.3-Remark}%

\smallskip

\noindent  %
{\bf Remark~9.3.} \label{9.3-Remark}  %
 Any continuous deformation
defined by the algorithm~(\ref{9.2-equation}) is generated by two
kinds of transformations, continuous one dimensional
symmetrization and $(k,n)$-Steiner symmetrization. Such a
``hybrid'' will inherit the properties that both its ``parents''
have and often the proof of a particular property of the
transformation $T_k^t$, for any {\bf fixed} $t$, can be given via
what we call the ``standard inductive argument''. More precisely
the latter means that if a certain statement is known to be true
for any $(k,n)$-Steiner symmetrization, as well as for a
corresponding continuous one dimensional symmetrization, and if
the corresponding property is invariant under translations and
scaling, then the corresponding result will remain true for any
continuous $k$-dimensional symmetrization defined by the
algorithm~(\ref{9.2-equation}).  We will demonstrate how this
``standard inductive argument'' works in the proof of
Theorem~\ref{9-1-Theorem} below. We want to emphasize here that if
the assumptions of some theorem or lemma deal with the
\textbf{varying} parameter $t$, then the standard inductive
argument is \textbf{not valid} and the proof of such statement is,
usually, not so straightforward,
cf. the proof of Theorem~\ref{9-2-Theorem}. %
%\erem  %

%\begin{figure}
%$$\includegraphics[scale=.3,angle=0]{iceberg} $$%
%\caption{SC $(1,2)$-Steiner symmetrization.} \label{6-figure}
%\end{figure}

\bt \label{9-1-Theorem}  %
 Let $T^t_k$ be an SC $(k,n)$-Steiner symmetrization
into a $(k,n)$- Steiner symmetrization corresponding to the
decomposition ${\mathbb{R}}^n={\mathbb{R}}^{n-k}\times{\mathbb{
R}}^k$ and let $T^t_{k,x'_0}$ be a restriction of $T^t_k$ onto the
slice $\mathbb{R}(x'_0):=\{x=(x_0',y):\,  y\in{\mathbb{R}}^k\}$.
Then $T^t_k$, respectively $T^t_{k,x_0'}$, acts on $\mathbb{R}^n$,
respectively on $\mathbb{R}(x'_0)$, as an open, compact, and
smoothing rearrangement, which is continuous from the
inside and from the outside. %
\et  %

{\it Proof}\,. The result is already known for the one dimensional
symmetrization, see Lemma~\ref{Basic geometric properties of
SCSym}.

Assume now that the theorem is true for any integer $n\ge 3$ and
some $k-1$ such that $2\le k<n$. Now for a fixed $t$, $0\le t\le
1$, let $T^t_k$ be a continuous symmetrization as described in the
formulation of the theorem.  If $t=1$ or $t=t(j)$ for some $j\ge
2$, then the $T^t_k$-transformation coincides with the
$(k,n)$-Steiner symmetrization or it coincides with  a finite
composition of $(k-1,n)$-Steiner symmetrizations, respectively. In
both cases the theorem is known to be true.

Let $t(j)<t<t(j+1)$ for some integer $j\ge 2$. Then according to
the algorithm~\ref{9.2-equation}, the $T^t_k$-transformation can
be represented as a finite composition of $(k-1,n)$-Steiner
symmetrizations postcomposed by a continuous $(k-1,n)$-Steiner
symmetrization. Then using the known result for $(k-1,n)$-Steiner
symmetrizations and the inductive assumption on the continuous
$(k-1,n)$-Steiner symmetrizations we conclude that the theorem is
true for any $k$, $1\le k\le n$. \hfill $\Box$
\smallskip

The following lemma  shows that continuous symmetrizations are
monotone in measure on slices. Its proof  easily follows  from
Theorem~\ref{9-1-Theorem} and  the non-expansivity property of Lemma~\ref{ Non-expansivity of rearrangements}. %
\bl \label{9-1-Lemma} %
 Let $T^t_k$ be a continuous symmetrization as in
 Theorem~\ref{9-1-Theorem} and  let $\Omega\in
{\mathcal{G}}_b\cup{\mathcal{F}}$. If $0\le t_1<t_2<t_3\le 1$,
then  %
\be \label{9.3-equation}  %
{\mathcal{L}}^k(\Omega^{t_2}(x')\triangle\Omega^{t_3}(x'))\le
{\mathcal{L}}^k(\Omega^{t_1}(x')\triangle\Omega^{t_3}(x'))  %
\ee  %
for every  $x'\in {\mathbb{R}}^{n-k}$ and therefore  %
\be  \label{9.4-equation}  %
{\mathcal{L}}^n(\Omega^{t_2}\triangle\Omega^{t_3})\le
{\mathcal{L}}^n(\Omega^{t_1}\triangle\Omega^{t_3}). %
\ee  %
\el  %

\smallskip

\bt \label{9-2-Theorem}  %
 Let $T^t_k$ be a continuous
symmetrization as in Theorem~\ref{9-1-Theorem} and let $t_s$,
$s=1,2,\ldots$ be an increasing sequence such that $t_s\in [0,1]$
and $t_s\to t_0$. Then  %
\be \label{9.5-equation} %
d(\partial\Omega^{t_s},\partial\Omega^{t_0})\to 0 \quad {\mbox{as\
$s\to \infty$}} %
\ee  %
if $\Omega$ is a bounded open set and %
\be \label{9.6-equation} %
d(\Omega^{t_s},\Omega^{t_0})\to 0\quad {\mbox{as \ $s\to \infty$}} %
\ee %
if $\Omega$ is a compact set.  %
\et  %

{\it Proof}. To be specific, we will assume that  $\Omega$ is a
bounded open set. In the case of compact sets the proof is easier
and is left to the reader.

\smallskip

We proceed by induction. For $k=1$, the statement is true by
Lemma~\ref{Continuity of T^t} modulo some remarks in
Section~\ref{Rescaling and limit cases}. Suppose  that for every
SC $(j,n)$-Steiner symmetrization with $1\le j\le k-1$ the
conclusion of the theorem is true. Our goal now is to prove that
the theorem is true for $j=k$.

\medskip

If $t_0<1$, then $t(m)<t_0\le t(m+1)$ for some positive integer
$m$. By  Definition~\ref{continuous $(k,n)$-symmetrization} and
Algorithm~(\ref{9.2-equation}), the set $\Omega^{t(m+1)}$ is the
image of $\Omega^{t(m)}$ under some   $(k-1,n)$-Steiner
symmetrization $S_1$ and moreover there exists a continuous
$(k-1,n)$-Steiner symmetrization $T^{1,\tau}$ into $S_1$ such that %
\be  \label{9.7-equation}  %
 \Omega^{t}=T^{1,\tau}(\Omega^{t(m)}) \quad {\mbox{for all \ \ $t(m)<t\le
 t(m+1)$,}}%
\ee   %
 where $\tau=\tau(t)=(t-t(m))/(t(m+1)-t(m))$.  Since $\tau(t_s) \to  \tau(t_0)$ as
$s\to \infty$, the theorem follows from (\ref{9.7-equation}) and
the inductive assumption.

\medskip

In the case  $t_0=1$ the proof is by contradiction. If
(\ref{9.5-equation}) does not hold, then as in the proof of
Lemma~\ref{Continuity of T^t} we can find $\ep>0$ and, if
necessary, subsequences (for which we keep previous notation),
$t_s$ and $\Om_s=\Om^{t_s}$, and a sequence of points $x_s\in
\mathbb{R}^n$ with $x_s\to x_0\in \mathbb{R}^n$ such that one of
the following
two conditions is satisfied: %
\begin{enumerate} %
\item[\bf{(a)}] $x_s\in\partial\Omega^*$ and  ${\rm
dist}(x_s,\partial
\Omega_s)\ge \varepsilon$ for all $s=0,1,2,\ldots $  %
\vspace{4pt} %
\item[\bf{(b)}]  $x_s\in \partial \Omega_s$ for $s=1,2,\ldots$ and
${\rm dist}(x_s,\partial \Omega^*)\ge \varepsilon$  for all $s=0,1,2,\ldots $  %
\end{enumerate}

\medskip

In the case {\bf{(a)}} we consider two  subcases:

(i) $B^{(n)}_{\varepsilon}(x_0)\cap\Omega^{t_s}=\emptyset $ \ \
for all $s=1,2,\ldots$

(ii) $B^{(n)}_{\varepsilon}(x_0)\subset \Omega^{t_s}$ \ \  for all
$s=1,2,\ldots$

\smallskip

Case (i).  Let $t(m_s)<t_s\le t(m_s+1)$. In the case under consideration there is $\delta_0>0$ such that %
\be  \label{9.8-equation}%
{\mathcal{L}}^n\left(\Omega^{t_s}\triangle\Omega^*\right)\ge \delta_0  %
\ee  %
for  all $s=1,2,\ldots$ Since  %
$$  %
\Om^{t(m_s)}\triangle \Om^{t(m_s+1)}\subset
\left(\Om^{t(m_s)}\triangle
\Om^*\right)\cup\left(\Om^{t(m_s+1)}\triangle \Om^*\right),  %
$$  %
 by equation~\ref{3.17-equation} of Theorem~\ref{Sarvas approximation lemma} we
have  %
\be  \label{9.9-equation}%
{\mathcal{L}}^n\left(\Omega^{t(m_s)}\triangle
\Omega^{t(m_s+1)}\right)\le
{\mathcal{L}}^n\left(\Omega^{t(m_s)}\triangle
\Omega^*\right)+{\mathcal{L}}^n\left(\Omega^{t(m_s+1)}\triangle
\Omega^*\right)\to 0   %
\ee  %
as  $m\to \infty$. Next by Lemma~\ref{9-1-Lemma},  %
$$  %
{\mathcal{L}}^n\left(\Omega^{t_s}\triangle\Omega^{t(m_s+1)}\right)\le
{\mathcal{L}}^n\left(\Omega^{t(m_s)}\triangle
\Omega^{t(m_s+1)}\right),  %
$$  %
which together with (\ref{9.9-equation})  contradicts
(\ref{9.8-equation}).

\medskip

Case (ii).  Let $x_0= (x'_0,y_0)$, where $\ x'_0\in
{\mathbb{R}}^{n-k}$. Let  $\Omega^*(x'_0)$ and $\Omega^{t}(x'_0)$
denote the restrictions of $\Om^*$ and $\Om^t$ respectively in the
slice $\mathbb{R}_{x'_0}$, each of which is not empty  by the
condition (ii).  Since $\Omega^*(x'_0)$ is
 a  $(k,k)$-Steiner symmetrization of $\Omega(x'_0)$,
 it is  a nonempty $k$-dimensional ball. By assumption {\bf{(a)}},  $x_0\in
\partial \Omega^*$ and therefore %
$$  %
{\mathcal{L}}^k(B^{(k)}_{\varepsilon}(y_0)\setminus
\Omega^*(x'_0))\ge \delta_1  %
$$  %
for some $\delta_1>0$. Hence,  %
\be    \label{9.10-equation}  %
{\mathcal{L}}^k(\Omega^*(x'_0)\triangle \Omega^{t_s}(x'_0))\ge
\delta_1,  \quad s=1,2,\ldots  %
\ee  %

By inequality (\ref{9.3-equation}) of Lemma~\ref{9-1-Lemma} we have %
$$  %
{\mathcal{L}}^k(\Omega^{t_s}(x'_0)\triangle\Omega^{t(m_s+1)}(x'_0))\le
{\mathcal{L}}^k(\Omega^{t(m_s)}(x'_0)\triangle
\Omega^{t(m_s+1)}(x'_0))\,.  %
$$  %
By equation~\ref{3.17-equation} of Theorem~\ref{Sarvas approximation lemma}, %
$$  %
 {\mathcal{L}}^k(\Omega^{t(m)}(x'_0)\triangle\Omega^*(x'_0)) \to
0 \quad \mbox{as} \quad m\to \infty.   %
$$  %
As in the case (i), the latter two  relations contradict
(\ref{9.10-equation}). This completes the proof in the case
{\bf{(a)}}.

\medskip

Now we turn to the condition  {\bf{(b)}}.  In this case
equation~\ref{3.16-equation} of Theorem~\ref{Sarvas approximation
lemma} implies that there is a neighborhood ${\mathcal{U}}$ of
$x_0$
such that one of the following two conditions is satisfied: %

\begin{enumerate} %
\item[(i)]  %
\ ${\mathcal{U}}\cap \Om^*=\emptyset$ \ and ${\mathcal{U}}\cap
\Om^{t(m)}=\emptyset$ \ for all $m$ large enough, %
\item[(ii)]  %
\ ${\mathcal{U}}\subset \Om^*$ \ and ${\mathcal{U}}\subset
\Om^{t(m)}$ \ for all $m$ large enough. %
\end{enumerate}

\medskip

 For $\delta>0$, $\tau>0$, and $x_0=(x'_0,y_0)$ with $y_0\in \mathbb{R}^k$, let  %
$$   %
Q_\delta(x_0)= B^{(n-k)}_\delta(x'_0)\times{\mathbb{R}}^k, \quad
Q(x_0,\delta,\tau)=B_\delta^{(n-k)}(x'_0)\times B_\tau^{(k)}. %
$$ %
Considering case (i), we first assume that $|y_0|>0$. Since for
every $x'\in \mathbb{R}^{n-k}$ the section $\Om^*(x')$ is a ball
in $\mathbb{R}^k$ centered at the origin, condition (i) implies
that
we can find $\delta_0>0$ and $\tau_0>0$ sufficiently small such that %
\be  \label{9.11-equation}  %
 \Omega^*\cap Q_{\delta_0}(x_0)\subset
 Q(x_0,\delta_0,|y_0|-\tau_0) \quad {\mbox{and}} \quad \Omega_m\cap Q_{\delta_0}(x_0)\subset
 Q(x_0,\delta_0,|y_0|-\tau_0)%
\ee   %
for all $m$  large enough. According to (\ref{9.2-equation}) the
set $\Om^{t_s}$ is obtained as a result of the action of a certain
continuous $(k-1,n)$-Steiner symmetrization of the set
$\Om^{t(m_k)}$ into the set $\Om^{t(m_s+1)}$. Since the
intersections $\Om^{t(m_k)}\cap Q_{\delta_0}(x_0)$ and
$\Om^{t(m_k+1)}\cap Q_{\delta_0}(x_0)$ both belong to
$Q(x_0,\delta_0,|y_0|-\tau_0)$ and since the set
$Q(x_0,\delta_0,|y_0|-\tau_0)$ is invariant under the continuous
$(k-1,n)$-Steiner symmetrization mentioned above, it follows from
the property of monotonicity in slices of the continuous
symmetrization that $\Om^{t_s}\cap Q_{\delta_0}(x_0)\subset
Q(x_0,\delta_0,|y_0|-\tau_0)$ for all $s$ sufficiently large. The
latter contradicts the assumptions $x_s\in \partial \Om^{t_s}$ and
$x_s\to x_0$.  %

\smallskip

If $y_0=0$ then, in the case (i), $\Om^*\cap
Q_{\delta_0}(x_0)=\emptyset$ for some $\delta_0>0$. Hence the
intersection $\Om\cap Q_{\delta_0}(x_0)$ is also empty and
therefore $\Om^{t_s}\cap Q_{\delta_0}(x_0)=\emptyset$ for all $s$,
which again contradicts the assumptions $x_s\in \partial
\Om^{t_s}$ and $x_s\to x_0$.  This finishes the proof in the case
(i).  %

In the case (ii), we argue as above to conclude that there are
$\delta_0>0$ and $\tau_0>0$ such that %
\be  \label{9.12-equation}  %
 Q(x_0,\delta_0,|y_0|+\tau_0)\subset \Omega^*\cap Q_{\delta_0}(x_0) \quad {\mbox{and}} \quad
 Q(x_0,\delta_0,|y_0|+\tau_0)\subset \Omega_m\cap Q_{\delta_0}(x_0)%
\ee   %
for all $m$ sufficiently large. As above, the monotonicity
property of the continuous $(k-1,n)$-Steiner symmetrization from
$\Om^{t(m_s)}$ into $\Om^{t(m_s+1)}$ shows that
$Q(x_0,\delta_0,\tau_0)\subset \Om^{t_s}$ for all $s$ large
enough. The latter contradicts the assumptions $x_s\in \partial
\Om^{t_s}$ and $x_s\to x_0$. The proof of the theorem is complete.  \hfill $\Box$  %

\medskip

\bl  \label{9-2-Lemma}   %
Let $\Omega,\Omega'\in {\mathcal{G}}_{n,b}$ be such that
$\overline{\Omega'}\subset\Omega$ and let $(\cdot)^t$ denote an SC
symmetrization into the $(k,n)$-Steiner symmetrization {\it Sym}
corresponding to the decomposition
${\mathbb{R}}^n={\mathbb{R}}^{n-k}\times {\mathbb{R}}^k$. Then for
every $t$,  there exists a finite number of transformations
$P_1,\ldots,P_N$
such that %
$$  %
 \Omega'_N:= P_1\circ\cdots\circ
P_N\,(\Omega')\subset \Omega^t %
$$  %
and  %
$$  %
{\rm dist}(\partial \Omega',\partial\Omega)\le {\rm
dist}(\partial\Omega'_N,\partial\Omega^t), %
$$  %
where  each $P_j$ is either a polarization with the polarizer
$H_j$ such that $\Sigma_j:=\partial H_j\supset {\mathbb{R}}^{n-k}$
or a shift in a
direction $\vec{l}$ such that $\vec{l}\perp {\mathbb{R}}^{n-k}$. %
\el  %

 {\it Proof}\,. If $k=1$, then the considered symmetrization is
 one dimensional and the lemma follows from Corollary~\ref{approximation of open sets by
 polarization}.

 Let $k\ge 2$. If $t=1$ or $t=t(j)$ for some $j\ge 2$, then
 the $T_k^t$-transformation coincides with the $(k,n)$-Steiner
 symmetrization or with a finite composition of $(k-1,n)$-Steiner
 symmetrizations. In either case the result is known to be true,
 cf. Lemmas~7.1 and 7.2 \cite{BS}. Therefore for any fixed $t$,
 the desired conclusion follows via the standard inductive
 argument as we explained in Remark~9.3 and demonstrated in the
 proof of Theorem~\ref{9-1-Theorem}.   \hfill $\Box$

\section{Integral inequalities} \label{Integral inequalities}  %
\setcounter{equation}{0}  %

Many   integral inequalities well-known in the theory of
symmetrization remain valid
 for the continuous symmetrization as well.
 We start with the following $L^p$-continuity lemma.

\bl \label{10-1-Lemma}  %
 Let $u\in L^p _{0+} (\mathbb{R}^n )$ , $1\leq p<\infty$ and let $u^t=T^t(u)$,
 $0\le t\le 1$ denote a continuous $(k,n)$-Steiner symmetrization. Then the
mapping $t\mapsto u^t $ is continuous in $L^p$-norm; i.e.,
 \be  \label{10.1-equation}  %
\lim\limits_{t\to 0} \| u^{s+t} - u^s \|_p  =  0 \qquad {\mbox{for
all \ \ $0\le s\le 1$.}} %
\ee  %
\el  %

\noindent  %
{\it Proof}. First we assume that  $u$ is a step function. Then
$$
u=\delta \sum\limits_{i=1}^{m} \chi (E_i )
$$
with some  $\delta >0 $ and some measurable sets $E_1\supseteq
\ldots
 \supseteq E_m $. Let $E_i^\tau=T^\tau(E_i)$, $0\le \tau \le 1$. By Lemma~\ref{9-1-Lemma}, we
 conclude that  %
$$  %
\| u^{s+t} -u^s \| _p  \leq  \delta \sum\limits_{i=1}^m
{\mathcal{L}}^n (E_i ^{t+s} \Delta E_i^s )\to 0 \quad
{\mbox{ as  $t\to 0$.}}  %
$$  %

 In general, $u$ can be approximated in the $L^p$-norm by
 step functions $u_m$. Then using the inequality  %
 $$  %
 \|u^{s+t}-u^s\|\le
 \|u^{s+t}-u_m^{s+t}\|+\|u_m^{s+t}-u_m^s\|+\|u_m^s-u^s\| %
 $$  %
 and the non-expansivity Lemma~\ref{ Non-expansivity of rearrangements} we obtain (\ref{10.1-equation}).  $\hfill \Box$

\smallskip

Our next lemma  shows that a similar continuity property remains
valid for the space of continuous functions.

\bl  \label{10-2-Lemma}   %
 Let $u\in C_0({\mathbb{R}}^n)\cap L^1_+({\mathbb{R}}^n)$. Then for
 every $0\le s\le 1$, $\|u^{s+t}-u^s\|_\infty\to 0$ as $t\to 0$.  %
 \el  %

 {\it Proof}\,.  %
 Fix $0\le s\le 1$ and $\e >0$. Since $u\in C_0(\mathbb{R}^n)\cap L^1_+(\mathbb{R}^n)$,
 we can find $v\in W^{1,\infty}_0({\mathbb{R}}^n)
 \cap L^1_+({\mathbb{R}}^n)$
 such that $\|v-u\|_{\infty}<\e/3$.  Then using the non-expansivity property of Lemma~\ref{ Non-expansivity of rearrangements}
  we obtain for all $t$ small enough %
\begin{eqnarray*}  %
\|u^s-u^{s+t}\|_{\infty} &\leq &
\|u^s-v^s\|_{\infty}+\|v^s-v^{s+t}\|_{\infty}+\|v^{s+t}-u^{s+t}\|_{\infty}\\
 &\leq& 2\|u-v\|_{\infty}+\|v^s-v^{s+t}\|_{\infty}<\e\, , %
\end{eqnarray*}  %
which proves the lemma. $\hfill \Box $

\smallskip

The rest of this section contains seven theorems. The first six of
them are  monotonicity statements that correspond to the
inequalities, which are well-known in the theory of Steiner
symmetrization. Their proofs in all cases follow the same scheme.
First we use an approximation by
 polarizations to prove the required monotonicity for the continuous
SC  $1$-symmetrization. Then we apply the standard inductive
argument to
 show that the same kind of monotonicity  holds for the continuous
 $(k,n)$-symmetrization for any $k$.

 Theorem~\ref{10-1-Theorem} below is related to the well-known convolution type inequalities, cf. \cite[Lemmas~8.1, 8.2]{BS},
 \cite{Beckner}, \cite[Corollary~2]{Ba}. It shows that
 convolutions are, in fact, monotone functions of the parameter of
 the corresponding continuous symmetrization.

 \bt \label{10-1-Theorem}  %
  Let $u,v,w\in {\mathcal{S}}_+$ with $w=w^* $, where $(\ )^*$ denotes a $(k,n)$-Steiner symmetrization
   and let $j$ be a Young-function. Then the integral %
$$  %
\int \!\!\! \int\limits_{{\mathbb{R}}^{2n} } j( |u^t (x)-v^t (y)|)
w(x-y)\,dxdy  %
$$  %
decreases  in $0\le t\le 1$ provided that the integral converges for $t=0$..  %
\et  %

{\it Proof}. First we prove the theorem for the SC
$1$-symmetrization. By the semi-group property (\ref{semigroup
formula}), we have to prove the inequality %
\be \label{10.7-equation} %
\int\!\!\int_{\mathbb{R}^{2n}} j(|u^{t}(x)-v^{t}(x)|)w(x-y)\,dxdy
\le \int\!\!\int_{\mathbb{R}^{2n}}
j(|u(x)-v(x)|)w(x-y)\,dxdy %
\ee  %
for $-\infty <t<\infty$.

First observe that in view of non-expansivity Lemma~\ref{
Non-expansivity of rearrangements}, we can restrict ourselves to
the case that $u$ and $v$ are continuous functions with bounded
support. Then we define inductively two sequences $u_m$ and $v_m$
of polarizations of $u$ and $v$ as in Lemma~\ref{Lemma -
Approximating sequence in p-norm}, where the corresponding
half-spaces $H_m$ are chosen in such a way that the minimality
property (\ref{equation-7.11}) is satisfied. By Lemma~\ref{Lemma -
convergence for continuous functions}, the sequences $u_m$ and
$v_m$ converge in $C(\mathbb{R}^n)$ to $u^t$ and $v^t$,
respectively. Then (\ref{10.7-equation}) follows by applying
Lemma~8.1 in \cite{BS} inductively.

To get values of $t$ varying in the standard range $0\le t\le 1$,
we may use scaling and translation as it is explained in
Section~\ref{Rescaling and limit cases}. Of course, the latter two
operations do not change the integrals in (\ref{10.7-equation}).

Finally, it is well known that the desired monotonicity result
holds for the $(k,n)$-Steiner symmetrization for any $k$, $1\le
k\le n$, see \cite[Lemma~8.2]{BS}. Therefore for $k\ge 2$, the
proof of Theorem~\ref{10-1-Theorem} follows via our standard
inductive argument. \hfill $\Box$

\smallskip

The Dirichlet-type inequalities for functions and their
$(k,n)$-Steiner symmetrizations also admit  continuous
counterparts.

\bt  \label{10-2-Theorem}  %
Let $u\in W^{1,p}_+ ({\mathbb{R}}^n )\cap {\mathcal{S}}_+$, $1\leq
p\leq +\infty$. Then  $u^t \in W^{1,p}_+({\mathbb{R}}^n )\cap
{\mathcal{S}}_+$ for all $0\le t\le 1$ and
 $\| \nabla u^t \| _p $ decreases in $0\le t\le 1$.  %

 Furthermore, if $V$ is some linear subspace which either contains
 all ``$y$-directions" $x_{n-k+1},\ldots,x_n$, or is orthogonal to
 each of these directions, then $\|\nabla_V u^t\|_p$ decreases in
 $0\le t\le 1$.
 \et  %

 {\it Proof}. For the $(k,n)$-Steiner symmetrization this result
 is well known, see \cite[Theorem~8.2]{BS}. So, we give the proof
 for the case of the SC $1$-symmetrization. Then for $k\ge 2$,
 Theorem~\ref{10-2-Theorem} will follow via the standard inductive
 argument.

 Let $u^t$ denote the SC $1$-symmetrization of $u$. The semigroup
 property (\ref{semigroup formula}) shows that it is enough to
 prove the inequalities  %
 \be   \label{10.8-equation}  %
\|\nabla u^t\|_p\le \|\nabla\|_p \quad {\mbox{and}}\quad
\|\nabla_Vu^t\|_p\le \|\nabla_Vu\|_p  %
\ee  %
for $-\infty< t<\infty$.

For a fixed $t$, let $u_m$ be the sequence of polarizations of $u$
defined by Lemma~\ref{Lemma - Approximating sequence in p-norm},
which converges to $u^t$ in $L^p(\mathbb{R}^n)$. We consider two
cases.

{\bf (i)}  Let $1<p<\infty$. Since $\|\nabla u_m\|_p=\|\nabla
u\|_p$, by Lemma~5.3 \cite{BS} we can find a function $v\in
W^{1,p}(\mathbb{R}^n)$ and a subsequence $u_{m'}$ such that %
$$
u_{m'}\rightharpoonup v \quad {\mbox{weakly in
$W^{1,p}(\mathbb{R}^n$).}} %
$$  %
This means that for every $\varphi\in C_0^\infty(\mathbb{R}^n)$
and $i=1,\ldots,n$, %
$$  %
\int_{\mathbb{R}^n} \varphi v_{x_i}\,dx\leftarrow
\int_{\mathbb{R}^n}\varphi \frac{\partial (u_{m'})}{\partial
x_i}\,dx=-\int_{\mathbb{R}^n}\varphi_{x_i}u_{m'}\,dx \rightarrow
-\int_{\mathbb{R}^n}\varphi_{x_i}u^t\,dx, %
$$  %
that is, $v=u^t$. In view of the lower semi-continuity of the norm
it follows that %
$$  %
\|\nabla u^t\|_p\le \liminf \|\nabla (u_m)\|_p=\|\nabla u\|_p  %
$$  %
that is the first inequality in (\ref{10.8-equation}). Using
equation~(5.10) of Lemma 5.3 \cite{BS} one can prove the second
inequality in (\ref{10.8-equation}) analogously.

{\bf (ii)}  Let $p=1$. By Lemma~5.3 \cite{BS}, the functions
$|\nabla u_m|$ and $|\nabla u|$ are rearrangements of each other.
This means that for every $\delta>0$,  %
\begin{alignat}{10} %
\sup\left\{\int_E |(u_m)_{x_i}|\,dx: {\ } L^n(E)\le \delta
\right\} &\le& {\ }\sup \left\{ \int_E|\nabla u_m|\,dx: \
L^n(E)\le \delta\right\}
\nonumber \\
&=&  \sup \left\{ \int_E|\nabla v|\,dx: \ L^n(E)\le
\delta\right\}. \nonumber  \   %
\end{alignat}  %
Hence, if $E_k$ is any  sequence of measurable sets with $\lim
(L^n(E_k))=0$, we infer that %
$$  %
\sup \left\{\int_{E_k}|(u_m)_{x_i}|\,dx:\ m\in
\mathbb{N}\right\}\to 0 \quad {\mbox{as $k\to \infty$.}}  %
$$  %

Applying a well-known weak compactness criterion in
$L^1(\mathbb{R}^n)$, see \cite[p. 199]{Alt}, we again can extract
a subsequence $u_{m'}$ converging weakly in
$W^{1,1}(\mathbb{R}^n)$. Finally, proceeding as in case {\bf (i)}
the assertion follows in the case $p=1$ too.  \hfill $\Box$

\smallskip

 The following theorem gives an analog of the previous theorem for
the spaces of continuous functions.  %

 \bt \label{10-3-Theorem}  %
 Let $u\in C({\mathbb{R}}^n ) \cap {\mathcal{S}}_+$. Then
$\omega _{u^t}$ decreases in $0\le t\le 1$.  %
\et  %

{\it Proof}. Let $u\in L^1_+(\mathbb{R}^n)$. For the SC
$1$-symmetrization the result follows from the proof of
Lemma~\ref{Lemma - convergence for continuous functions}. In the
general case we choose a sequence $u_m$ of functions in
$C(\mathbb{R}^n)\cap L^1_+(\mathbb{R}^n)$ converging to $u$ in
$C(\mathbb{R}^n)$. Then applying the non-expansivity Lemma~\ref{
Non-expansivity of rearrangements} we obtain  %
$$  %
\|(u_m)^t-u^t\|_\infty \le \|u_m-u\|_\infty, \quad m=1,2,\ldots  %
$$  %
and the assertion for the SC $1$-symmetrization follows. For $k\ge
2$, we apply the standard inductive argument. \hfill $\Box$

It is also easy to prove the monotonicity property of convex
functionals.  %

\bt  \label{10-4-Theorem}  %
Let $u\in W^{1,1} _+ ({\mathbb{R}}^n )$ and let $j$ be a
Young-function such that  \\ %
$ \int_{{\mathbb{R}}^n } j(|\nabla
u|)\,dx <\infty$. Then the integral mean $\int_{{\mathbb{R}}^n}j(|
\nabla u^t |)\,dx $ decreases in $0\le t\le 1$.

 Furthermore, if $V$ is a linear subspace which either contains all ``$y$-
directions'' $x_{n-k+1},\ldots,x_n$,  or is orthogonal to each of
these directions, then $\int_{{\mathbb{R}}^n } j(| \nabla _V u^t
|)\,dx$ decreases in $0\le t\le 1$.  %
\et  %

{\it Proof}. The standard inductive argument still works.
Therefore, it is remains to prove the theorem for the SC
$1$-symmetrization.

Fix $t$, $-\infty<t<\infty$. By the semigroup property
(\ref{semigroup formula}), we have to prove the inequalities  %
\be  \label{10.9-equation}  %
\int_{\mathbb{R}^n} j(|\nabla u^t|)\,dx \le \int_{\mathbb{R}^n}
j(|\nabla u|)\,dx, \quad \quad    \int_{\mathbb{R}^n} j(|\nabla_V
u^t|)\,dx\le \int_{\mathbb{R}^n} j(|\nabla_V u|)\,dx.  %
\ee  %

Assume first that $u\in C_{0+}^{0,1}(\mathbb{R}^n)$. If we choose
a sequence of polarizations of $u$ converging to $u^t$ in
$L^1(\mathbb{R}^n)$, we conclude from Lemma~5.3 \cite{BS} that  %
$$  %
\int_{\mathbb{R}^n} j(|\nabla u_m|)\,dx =\int_{\mathbb{R}^n}
j(|\nabla u|)\,dx.  %
$$  %
Because of the weak lower semi-continuity of the integral
functional this leads to  %
$$  %
\int_{\mathbb{R}^n} j(|\nabla u|)\,dx =\liminf \int_{\mathbb{R}^n}
j(|\nabla u_m|)\,dx \ge \int_{\mathbb{R}^n} j(|\nabla u^t|)\,dx. %
$$  %
If $u\in W_+^{1,1}(\mathbb{R}^n)$, we choose a sequence $v_m\in
C_{0+}^{0,1}(\mathbb{R}^n)$ such that %
$$  %
v_m\to u \quad {\mbox{in $W^{1,1}(\mathbb{R}^n)$}}  %
$$  %
and  %
$$  %
\int_{\mathbb{R}^n} j(|\nabla v_m|)\,dx \to \int_{\mathbb{R}^n}
j(|\nabla u|)\,dx.  %
$$  %
This means that we have for a subsequence $v_{m'}$  %
$$  %
(v_{m'})^t \rightharpoonup u^t \quad {\mbox{weakly in
$W^{1,1}(\mathbb{R}^n)$,}} %
$$and we conclude by the weak lower semi-continuity of the functionals that the first inequality in
(\ref{10.9-equation}) holds true.

One can prove the second inequality in (\ref{10.9-equation})
analogously.  \hfill  $\Box$

\smallskip

An analog of Theorem~\ref{10-2-Theorem} also holds for
$BV$-functions.  %
\bt  \label{10-5-Theorem}  %
 If $u\in BV({\mathbb{R}}^n )\cap L^1_+ ({\mathbb{R}}^n)$, then $u^t \in
BV({\mathbb{R}}^n )\cap L^1_+({\mathbb{R}}^n)$ and
 $\| Du^t \| _{BV}$ decreases in $0\le t\le 1$.  %
 \et  %

 {\it Proof}. As in the previous theorems, we have to prove the
 result for the SC $1$-symmetrization. By the semigroup property
 (\ref{semigroup formula}) we have to prove that %
 \be  \label{10.10-equation} %
 \|Du^t\|_{BV}\le \|Du\|_{BV}  %
 \ee  %
 for any fixed $t$, $-\infty <t<\infty$.

 We choose a sequence of functions $u_m\in W_+^{1,1}
 (\mathbb{R}^n)$ which converges to $u$ in $BV(\mathbb{R}^n)$. By
 Lemma~5.3 \cite{BS}, the functions $(u_m)^t$ are equibounded in
 $W^{1,1}(\mathbb{R}^n)$. Therefore there are some function $v\in
 BV(\mathbb{R})^n$ and a sequence $u_{m'}$ such that %
 $$  %
 (u_{m'})^t\rightharpoonup v \quad {\mbox{weakly in
 $BV(\mathbb{R}^n)$.}}  %
 $$  %
 On the other hand from the inequalities  %
 $$  %
 \|(u_m)^t-u^t\|_1\le \|u_m-u\|_1  %
 $$  %
 we conclude that %
 $$  %
 (u_m)^t\to u^t \quad {\mbox{in $L^1(\mathbb{R}^n)$.}}  %
 $$  %

 Now let $\mu_i$ denote the Radon-measure which is associated with
 the weak partial derivative $v_{x_i}$, $i=1,\ldots,n$. Then we
 have for every $\varphi \in C_0^\infty(\mathbb{R}^n)$,  %
 $$  %
 \int_{\mathbb{R}^n} \varphi\,d\mu_i\leftarrow \int_{\mathbb{R}^n}
 \varphi \frac{\partial ((u_{m'})^t)}{\partial x_i}\,dx
 =-\int_{\mathbb{R}^n} \varphi_{x_i}(u_{m'})^t\,dx \rightarrow
 -\int_{\mathbb{R}^n}\varphi_{x_i}u^t\,dx,  %
 $$  %
 which means that $v=u^t$.

 Finally, the weak lower semi-continuity of the norm gives  %
 $$  %
 \|D u^t\|_{BV}\le \liminf \|\nabla((U_m)^t)\|_1=\lim \|\nabla
 u_m\|_1=\|Du\|_{BV}.   %
 $$  %
\hfill $\Box$  %

\smallskip

 Choosing for $u$ in Theorem~\ref{10-5-Theorem} the characteristic
 function of a set of finite perimeter we derive the following
 ``monotonic
 isoperimetric inequality'' in $\mathbb{R}^n$.  %
 \bt  \label{10-6-Theorem} %
 Let $E$ is a Caccioppoli-set in $\mathbb{R}^n$. Then the perimeter $\| D\chi_{E^t}\|_{BV}$ of $E^t $
is a decreasing function in $0\le t\le 1$.  %
\et  %

Now we prove that the mapping $t\mapsto u^t$ is continuous from
the left in  Sobolev-spaces
 $W^{1,p}_+({\mathbb{R}}^n)$.  %
 It is worth mentioning
here that an analog of this result does not hold in the space
$BV({\mathbb{R}}^n)\cap L^1_+({\mathbb{R}}^n )$. The
characteristic function $u=\chi_I$ of a single interval $I\subset
\mathbb{R}$ provides a simple counterexample in dimension $n=1$.  %

\bt \label{10-7-Theorem} %
Let $u\in W^{1,p}_+({\mathbb{R}}^n)$,  $1\leq p<\infty$ and let
 $t_m$, $m=1,2,\ldots$ be an increasing sequence in $(0,1)$ such that $t_m\to t$ as $m\to \infty$.
 Then
\be \label{10.2-equation} %
 u^{t_m } \longrightarrow u^t \qquad
\mbox{ in } \ W^{1,p} ({\mathbb{R}}^n).   %
\ee  %
\et  %

{\it Proof}. First we consider the case $p>1$. By
Lemma~\ref{10-1-Lemma}, we have $u^{t_m } \rightarrow u^t $ in
$L^p({\mathbb{R}}^n)$. From this  we conclude  that  %
 \be   \label{10.3-equation}%
 u^{t_m }\rightharpoonup u^t \qquad \mbox{ weakly in } \
 W^{1,p}({\mathbb{R}}^n). %
\ee   %
In addition,  Theorem~\ref{10-2-Theorem} combined with  the weak
lower semi-continuity of the $L^p$-norm implies  %
 \be   \label{10.4-equation}  %
 \lim_{m\to \infty} \|\nabla u^{t_m }\|_p=\|\nabla u^t\|_p .  %
\ee %
 Since for $p>1$ the spaces $W^{1,p} ({\mathbb{R}}^n)$ are uniformly
 convex, (\ref{10.2-equation}) follows from (\ref{10.3-equation})
 and (\ref{10.4-equation}).

\smallskip

Now let $p=1$. We note that (\ref{10.3-equation}) and
(\ref{10.4-equation}) remain true for $p=1$. Then we fix an index
$i\in \{ 1,\ldots ,n\} $ and set $v_m := (u^{t_m})_{x_i}$, $v:=
(u^t )_{x_i } $. Let $j(z):= \sqrt{1+z^2}-1$. Then $j(z)$ is a
Young-function which satisfies the inequality $j(z)\le z$.
Therefore  Theorem~\ref{10-4-Theorem} and the weak lower
semi-continuity property of the integral imply that  %
\be \label{10.5-equation}   %
\lim_{m\to \infty } \int_{{\mathbb{R}}^n } ( \sqrt{1+v_m^2}
-1)\,dx = \int_{{\mathbb{R}}^n}(\sqrt{1+v^2}-1)\,dx\,.  %
\ee  %
Applying Taylor's formula with a reminder we obtain  %
\begin{alignat}{10}  %
\int_{{\mathbb{R}}^n}(\sqrt{1+v_m^2}-1)\,dx  &\ge &
\int_{{\mathbb{R}}^n}(\sqrt{1+v^2}-1)\,dx+\int_{{\mathbb{R}}^n}\frac{v}{\sqrt{1+v^2}}(v_m -v)\,dx    \nonumber \\
 &+&\  \frac{1}{2} \int_{{\mathbb{R}}^n}\frac{(v_m -v)^2}{(
 1+c_m^2)^{3/2}}\,dx\,, \nonumber \ \ \ \ \ \ \ \  \ \ \ \ \ \ \ \ \ \ \ \ \ \ \ \ \ \ \ \ \ \ \ \ \ \  %
\end{alignat}  %
where $c_m^2=\max\{v^2,v_m^2\}$. Passing to the limit as $m\to
\infty $ we obtain  %
$$  %
\lim_{m\to \infty} \int_{{\mathbb{R}}^n}\frac{(v_m
-v)^2}{(1+c_m^2)^{3/2}}\,dx=0\,.  %
$$  %
The latter equation implies that  for every positive integer $k$,%
\be  \label{10.5.1-equation}  %
\lim_{m\to \infty}\int_{F_{k,m}} |v_m -v| \ dx = 0,  %
\ee   %
where $F_{k,m}=\{x\in \mathbb{R}^n:\,|v_m|\le k,\,|v|\le k\}$.
Since $v_m\rightharpoonup v$ weakly in $L^p ({\mathbb{R}}^n)$ we
find
that  %
\be  \label{10.6-equation}   %
\lim_{k\to \infty}\int_{G_{k,m}} |v_m |\,dx =0   %
\ee   %
uniformly over all $m$. Here $G_{k,m}=\{x\in
\mathbb{R}^n:\,|v_m|>k\}$.

Now (\ref{10.2-equation})  follows  from (\ref{10.5-equation}) an
(\ref{10.6-equation}) via the inequality  %
$$  %
\|v_m-v\|_1 \le \int_{F_{k,m}} |v_m -v|\,dx+2\int_{G_{k,m}}|v_m
|\,dx +2\int_{\{|v|>k\}}|v|\,dx\, ,  %
$$  %
which holds for all positive integers $m$ and $k$. $\hfill \Box $

\section{Comparison theorems}  \label{Comparison theorems}  %
\setcounter{equation}{0}  %

First we introduce a partial order $\,\prec^t\,$ related to the
continuous $(k,n)$-Steiner symmetrization $T^t$. %
\bd \label{Partial order} %
 For functions $f,g\in L^1_+({\mathbb{R}}^n)$ and $0\le t\le 1$, we
 write %
\be  \label{11.1-equation}  %
 f\prec ^t g \ \quad {\mbox{if and only if}} \quad \int_{{\mathbb{R}}^n} fh\,dx \le \int_{{\mathbb{R}}^n}
 g^th^t\,dx \quad {\mbox{for all $h\in L^{\infty}_+({\mathbb{R}}^n)$.}} %
\ee%
\ed  %

 The following theorem generalizes the well known equivalences in the theory of
 symmetrization (see \cite[Remark~10.1]{BS}, \cite{ALT1}) to the case of SC $1$-symmetrization.  %

\bt \label{11-1-Theorem}  %
Let  $f,g\in L^1_+({\mathbb{R}}^n)$ and let $(\cdot)^t$ denote the
SC $1$-symmetrization, $-\infty<t<\infty$.
Then the following relations are equivalent:  %
\begin{eqnarray} %
 & & f\prec^t g  \label{11.2-equation}\\  %
 &  &
\int_{{\mathbb{R}}^n} f^t h^t\,dx  \le  \int_{{\mathbb{R}}^n} g^t
h^t\,dx \quad {\mbox{for all $h\in L^{\infty}_+({\mathbb{R}}^n)$,}}\label{11.3-equation}\\  %
 & & j(f) \prec^t j(g) \quad {\mbox{for all Lipschitz-continuous
 Young-functions $j$,}}\label{11.4-equation}\\  %
 & & \int_{2t-y}^y f^t(x',s)\,ds \le \int_{2t-y}^y g^t(x',s)\, ds \quad {\mbox{for all $x'\in {\mathbb{R}}^{n-1}$ and every $y\ge
 t$,}}\label{11.5-equation}\\  %
 & & f^t (x',y)\leq g^t (x',y) \quad {\mbox{for all $x'\in {\mathbb{R}}^{n-1}$ and every $y\ge
 t$.}}\label{11.6-equation} %
\end{eqnarray}  %
\et  %

{\it Proof}. {\bf{(a)}} To prove that (\ref{11.2-equation})
implies (\ref{11.5-equation}), we fix $c>0$ and set $M_c:=\{ f>c
\} $. Then for fixed $x'_0\in \mathbb{R}^{n-1}$ and positive $\e$,
we define
a function %
$$  %
h_\e(x',y):=\kappa(x',y)\varphi_\e(x',y) %
$$
with $x=(x',y)\in \mathbb{R}^n$ and  %
$$  %
\varphi_\e=\varphi_\e(x',y):= \chi(\{(x',y)\,:|x'-x_0'| <\e\}), %
$$  %
$$  %
\kappa(x',y):=\left\{\begin{array}{ll} 1 & {\mbox{\  if  $f(x',y)>c$  and $y<y(x',t)$}}\\
0 &  {\mbox{\  otherwise,}} %
\end{array} \right. %
$$
where the separating function $y(x_0',t)=y_{M_c}(x'_0,t)$ is
defined in Section~\ref{ Continuous $(1,n)$-Steiner
symmetrization}.

 Now (\ref{11.2-equation}) implies
\begin{eqnarray*}
& &\int_{\{|x'-x_0'|<\e\}}\int_{2t-y(x',t)}^{y(x',t)}f^t(x',s)\,
ds
\,dx'=\int_{{\mathbb{R}}^n} f^t(h_{\e})^t\,dx = \int_{{\mathbb{R}}^n}f h_\e \,dx \le\\  %
& & \int_{{\mathbb{R}}^n} g^t(h_{\e})^t
\,dx=\int_{\{|x'-x_0'|<\e\}} \int_{2t-y(x',t)}^{y(x',t)}
g^t(x',s)\, ds\,dx'.  %
\end{eqnarray*}  %

Taking the limit here as $\e\to 0^+$, we obtain  %
$$  %
\int_{2t-y(x_0 ',t)} ^{y(x_0 ',t)} f^t (x_0 ',s)\ ds \le  %
\int_{2t-y(x_0 ',t)} ^{y(x_0 ',t)} g^t (x_0 ',s)\ ds,  %
$$  %
where $y(x'_0,t)=y_{M_c}(x'_0,t)$ depends on $c$. Since $c>0 $ can be chosen arbitrary small the latter inequality
implies (\ref{11.5-equation}).  %

\smallskip

{\bf (b)} To show that (\ref{11.2-equation}) implies
(\ref{11.6-equation}) we fix $(x_0',y_0)\in{\mathbb{R}}^n$ with
$y_0\ge t$. Then using notation of part {\bf (a)} of this proof we
note that the separating function $y(x_0',t)=y_{M_c}(x'_0,t)$
depends monotonously on the height of polarization $c$. Therefore
we can choose $c>0 $ small enough such that $y_0 \ge y(x_0',t)$.
Then we have $f^t(x_0',y)=f(x_0',y)$ for all $y\ge y_0$. Now we
choose the function $h$ in (\ref{11.1-equation}) to be the Dirac
$\delta$-function at $x_0=(x'_0,y_0)$, i.e. $h:=\delta(x_0)$. Then
we have $h=h^t$.
Finally applying (\ref{11.1-equation}) we derive  %
$$  %
f^t(x_0',y_0)=f(x_0',y_0)=\int_{{\mathbb{R}}^n}f h\,dx\le
\int_{{\mathbb{R}}^n} g^t h^t\, dx = g^t(x_0',y_0),  %
$$  %
which is (\ref{11.6-equation}).  %

\smallskip

{\bf (c)} Now we show that (\ref{11.5-equation}) and
(\ref{11.6-equation}) together imply (\ref{11.3-equation}). First
we assume that $h=\chi(M)$ is a characteristic function of $M\in
{\mathcal{M}}(\mathbb{R}^n)$. Let $y(x',t)=y_M(x',t)$ and let
$M(x',t)=M(x')\cap \{s:\,s>y(x',t)\}$. Then
(\ref{11.5-equation}) and (\ref{11.6-equation}) imply  %
\begin{eqnarray}
& &\int_{M^t }f^t\,dx
=\int_{{\mathbb{R}}^{n-1}}\left\{\int_{2t-y(x',t)}^{y(x',t)}f^t(x',s)\,ds+\int_{M(x',t)}
f^t(x',s)\,ds \right\} \, dx'\le  \label{11.7-equation}\\  %
& & \int_{{\mathbb{R}}^{n-1}}\left\{\int_{2t-y(x',t) }^{y(x',t)}
g^t (x',s)\,ds + \int_{M(x',t)} g^t(x',s)\,ds \right\} \,dx' =
\int_{M^t } g^t\,dx\,. \nonumber  %
\end{eqnarray}   %

 Next let $h$ be a step function, i.e.
$$
h:= \varepsilon \sum_{i=1}^m \chi (M_i ) %
$$  %
with some $\e>0$ and some sets $M_i \in
{\mathcal{M}}({\mathbb{R}}^n)$ such that $M_1 \supset \cdots
\supset M_m$. Applying inequality (\ref{11.7-equation}) to the
functions $\chi(M_i)$ we obtain the
desired inequality:  %
$$ %
\int_{{\mathbb{R}}^n}f^t g^t\,dx =\e \sum_{i=1}^m
\int_{M_i^t}f^t\,dx \le \e
 \sum_{i=1}^m \int_{M_i^t}g^t\,dx =\int_{{\mathbb{R}}^n} g^t
 h^t\,dx\,.%
$$  %

Finally, every  $h\in L^{\infty}_+({\mathbb{R}}^n)$ can be
approximated by step functions.   Therefore, in the general case (\ref{11.3-equation}) follows from the previous inequality. %

\smallskip

{\bf (d)}  To prove that (\ref{11.3-equation}) implies
(\ref{11.4-equation}), we may assume without loss of generality
that $j\in C^1 $. For $p\in L^{\infty}_+({\mathbb{R}}^n)$, we set
$h:=j'(f)p$. Since $j'$ is an increasing function
 we have $h^t=j'(f^t)p^t$. This equality combined with (\ref{11.3-equation}) implies  %
\be \label{11.8-equation}  %
\int_{{\mathbb{R}}^n}f^tj'(f^t )p^t\,dx \le \int_{{\mathbb{R}}^n}
g^t j'(f^t)p^t\,dx\, .  %
\ee   %
Since $j$ is a convex function we have  %
$$ %
j(g^t)-j(f^t) \ge j'(f^t)(g^t-f^t), %
$$   %
which together with (\ref{11.8-equation}) leads to (\ref{11.4-equation}). %

{\bf (e)}  Finally, (\ref{11.2-equation}) is a special case of
(\ref{11.4-equation}). \hfill  $\Box $  %

\smallskip

The following corollary shows that equivalencies
(\ref{11.2-equation}) -- (\ref{11.4-equation}) remain valid for
continuous $(k,n)$-Steiner symmetrization for all $k$, $1\le k\le
n$. Its proof follows from Theorem~\ref{11-1-Theorem} via the
standard inductive argument.  %

\bc  \label{11.01-Corollary}  %
Let $f,g\in L^1_+(\mathbb{R}^n)$ and let $(\cdot)^t$, where $0\le
t\le 1$, denote the continuous $(k,n)$-Steiner symmetrization.
Then equivalencies (\ref{11.2-equation}) -- (\ref{11.4-equation})
remain valid for all $0\le t\le 1$.  %
\ec  %

\bc  \label{11.1-Corollary}  %
 Let $f,g\in L^1_+({\mathbb{R}}^n)$ and let $(\cdot)^t$ denote
 the continuous $(k,n)$-Steiner symmetrization into a given
 $(k,n)$-Steiner symmetrization $(\cdot)^*$. Then for all $s$ and $t$ such that $0\le s<t\le 1$,
 the following implication holds: %

If $f\prec^s g$, then $f\prec^t g$.

In particular,  if $f\prec^t g$, then $f\prec^* g$. %
\ec  %

For the SC $1$-symmetrization Corollary~\ref{11.1-Corollary}
 follows immediately from the semigroup property (\ref{semigroup
formula}). Then for any $k\ge 2$, its proof follows  via the
standard inductive argument.

\smallskip

Now we prove two comparison lemmas concerning partial symmetry of
solutions of certain elliptic and parabolic PDE's. These results
and their proofs show, in fact, that the approach to comparison
theorems in partially symmetric domains, based on the continuous
symmetrization, is a closed relative of the Alexandrov's moving
plane method, see \cite{S} and \cite{GNN}. In the context of
comparison theorems the approach based on a continuous
symmetrization for the first time was
used by Solynin  \cite{S1}. %

\bl  \label{11-1-Lemma}   %
 Let $\Omega \subset {\mathbb{R}}^n$ be
a bounded domain, $c\ge 0$, $f,g\in L^2_+(\Omega)$, and let
 $u,v$ be solutions to the following boundary value problems:
\be  \label{11.11-equation}   %
u,v\in W_0^{1,2}(\Omega),\qquad -\Delta u+cu=f,\ \ -\Delta v+cv=g
\quad {\mbox{in \ $\Omega$.}}  %
\ee   %
For $0\le t\le 1$, let $(\cdot)^t$ denote a continuous
$(k,n)$-Steiner symmetrization. If for some $0\le t\le 1$, $\Omega
=\Omega^t$, $f=f^t$,
$g=g^t$, and $f\prec^t g$, then  %
\be  \label{11.12-equation}  %
 u=u^t,\qquad v=v^t  %
\ee  %
 and  %
\be  \label{11.13-equation}  %
 u\prec^t v. %
\ee  %
\el  %

{\it Proof}. First, we prove the lemma for the SC
$1$-symmetrization. Then, of course, $-\infty<t<\infty$.  Equality
(\ref{11.12-equation}) follows easily from the maximum principle.

To prove (\ref{11.13-equation}), we may assume without loss of
generality that $f,g\in C^{0,1}_{0+}(\Omega)$. Then $u$ and $v$
are classical solutions which are smooth in $\Omega $.

Let $G=\{x=(x',y):\,(x',2t-y)\in \Omega\}$ and let $w=u-v$. For
$(x',y)\in G$ and $\e\ge 0$, we define a function %
$$ %
W_{\e}(x',y):=\int_t^y \left(w(x',s)+w(x',2t-s)-\sup_{\partial
G\cap
H_t} w-\e \right)\,ds . %
$$   %
Since $f\prec^t g$ one can easily see
 that  %
 \be \label{11.14-equation}  %
 -\Delta W_{\e}+c W_{\e} \le 0 \quad {\mbox{in $G$.}} %
\ee  %
In addition, we have $W_{\e}=0$ on $\partial H_t$. Since solutions
of (\ref{11.14-equation}) satisfy the maximum principle we have:
$\sup_G W_{\e}=\sup_{\partial G\cap H_t} W_{\e}$. Since
$\Om=\Om^t$ one can easily see that the domain $G$ is convex in
the direction of $y$-axis. Therefore the positive direction of
$y$-axis points outward the domain $G$ on $\partial G\cap H_t$.
Since $\frac{\partial W_{\e}}{\partial y}\le-\e$ on $\partial
G\cap H_t$ we conclude that $W_{\e} \le 0$ in $G$ for all $\e>0$.
Taking the limit as $\e\to 0^+$, we obtain  %
\be \label{11.15-equation}   %
W_0 \le 0 \quad \quad {\mbox{in $G$}}.  %
\ee  %

%\smallskip

Let us show that $w\le 0$ in $\Om_-(t)$. If not, then
$\sup_{\Omega_-(t)} w
>0$. Since $-\Delta w +c w\le 0$ in $\Om_-(t)$
and $w=0$ on $\partial \Om \cap H_t$ the maximum principle implies
that $\sup_{\Om_-(t)}w =\sup_{\partial H_t}w$. The second supremum
here is attained at some point $x_0 \in \partial H_t$ with
$w(x_0)>0$. Since $w(x_0)\ge \sup_{\partial G \cap H_t}w$, the
inequality $w(x_0)>0$ implies that $W_0(x)>0$ for  $x\in G$ in a
small neighborhood of $x_0$, which contradicts
(\ref{11.15-equation}).

 Therefore we have $w\le 0$ in $\Om_-(t)$. Now (\ref{11.13-equation}) follows from
 (\ref{11.5-equation}), (\ref{11.6-equation}), and (\ref{11.15-equation}).

 Now for $k\ge 2$, the lemma follows via the standard inductive argument. \hfill  $\Box$

\medskip

Similar lemma holds also  for parabolic problems. Its proof
follows along the  lines of the proof of Lemma~\ref{11-1-Lemma}.
Therefore the details are left to the
reader. %

\bl \label{11-2-Lemma}  %
Let $(\cdot)^s$ denote a continuous $(k,n)$-Steiner
symmetrization, $0\le s\le 1$. Let $c\ge 0$, $T>0$ and let $\Om$
be a bounded domain in $\mathbb{R}^n$ such that $\Om=\Om^s$. Let
functions $u_0,v_0\in L^2_+(\Om)$ and $f,g\in L^2_+(\Om\times
(0,T))$ satisfy, respectively, the following conditions: %
$$  %
u_0 =u_0^s,\quad  v_0=v_0^s,\quad  u_0\prec^s v_0 %
$$  %
 and %
 $$   %
 f(\cdot,t)=f^s(\cdot,t),\quad g=g(\cdot,t)=g^s(\cdot,t),\quad
 f(\cdot,t)\prec^s g(\cdot,t) \quad {\mbox{for all $t\in (0,T)$.}}  %
 $$  %
  Further let $u,v\in L^2(0,T;W_0^{1,2}(\Om))\cap C([0,T];L^2(\Om))$
be solutions to the following initial boundary value problems: %
\begin{alignat}{10}  %
 u_t-\Delta u+cu=f, &{\ \ \ \ }& v_t-\Delta v+cv=g, &{\ }&\quad  {\mbox{in \ $\Om\times(0,T)$,}}\label{11.16-equation}\\%
 u(x,0)=u_0(x),\ \ \   &{\ \ \ \ }& v(x,0)=v_0(x),\ \ \  &{\ \ \ }& {\mbox{in \ $\Om$.}}\ \ \ \ \ \ \ \ \ \ \nonumber
\end{alignat} %
Then %
\be \label{11.17-equation} %
 u(\cdot ,t)=u^s (\cdot ,t),\quad
v(\cdot ,t)=v^s (\cdot ,t) \qquad {\mbox{for all $t\in (0,T)$}} %
\ee  %
 and  %
\be  \label{11.18-equation}   %
u(\cdot,t)\prec^s v(\cdot,t) \qquad {\mbox{for all $t\in (0,T)$.}}  %
\ee  %
\el  %

\medskip

To simplify some notations in our formulations and proofs, we
 will use the following definitions from \cite{BS}:  %

\bd  \label{11.01-definition} %
Let $\Omega \subset \mathbb{R}^n$ be a bounded open set, $c\ge 0$,
and $f\in L^2_+(\Omega)$. We say that $u$ solves the problem
$\mathbb{B}_1(\Omega,c,f)$ if $u$ is the solution to the following
boundary value problem:  %
\be  \label{11.20.1-equation}  %
u\in W^{1,2}_0(\Omega), \quad -\Delta u+cu=f \quad {\mbox{in
$\Omega$.}}  %
\ee  %
\ed  %

\bd  \label{11-1-definition}  %
Let $\Om\subset \mathbb{R}^n$ be a bounded open set, $c\ge 0$,
$f\in L^2_+(\Om)$ and $\gamma:\,\mathbb{R}^+_0\to \mathbb{R}^+_0$
be a continuous and nondecreasing function. We will say that
$\underline{u}$ {\it is a solution of problem
$\mathbb{B}_2(\Om,c,\gamma,f)$}, if $\underline{u}$ is the
nonnegative minimal solution of the following boundary value
problem:  %
\be \label{11.21-equation} %
u\in W^{1,2}_0(\Om), \quad u\ge 0, \quad -\Delta u+cu=\gamma(u)+f
\quad {\mbox{in \ $\Om$,}} %
\ee  %
that is,  %
\begin{enumerate}  %
\item[\bf{(i)}] %
$\underline{u}$ is a solution of the problem
(\ref{11.21-equation}), and %
\item[\bf{(ii)}] %
$0\le \underline{u}\le u$ for all other solutions $u$ of
(\ref{11.21-equation}). %
\end{enumerate}  %
\ed  %

For a brief discussion of important properties of the solutions to
the problems in Definitions~\ref{11.01-definition} and
\ref{11-1-definition}, we refer to \cite[Section~9]{BS}.

\bt  \label{11-2-Theorem}  %
 Let $(\cdot)^t$ denote a continuous $(k,n)$-Steiner
 symmetrization, $0\le t\le 1$. Let $\Omega\subset \mathbb{R}^n$ be a bounded open set,
 $c\ge 0$, $f\in L^2_+(\Omega)$, and let $\gamma$ be a Young function.
 For a fixed $t$, $0\le t\le 1$, let $g\in L^2_+(\Omega^t)$ be
 such that $g^t=g$ and $f\prec^t g$.  Let  $u$  and $v$ be the solutions to the problems
 $B_2(\Omega,c,\gamma,f)$ and $B_2 (\Omega^t,c,\gamma,g)$, respectively.
Then  %
\be \label{11.19-equation}  %
 v=v^t %
\ee %
 and %
\be \label{11.20-equation}  %
u\prec^t v.  %
\ee  %
\et  %

{\it Proof}. As well known, the assertion of this theorem holds
true for the $(k,n)$-Steiner symmetrization for any $k$, $1\le
k\le n$; see, for example,  Theorem~10.1 in \cite{BS}. Thus, in
view of the standard inductive argument we have to prove it for
the SC $1$-symmetrization only. Then, of course,
$-\infty<t<\infty$.

{\bf (1)}  First we assume that $\gamma \equiv 0$ and that $f$ is
a simple function with compact support in $\Omega$, see
Definition~\ref{Simple function}. For a fixed $t\in \mathbb{R}$,
let $\tilde v$ denote the solution to the problem
$\mathbb{B}_1(\Omega^t,c,f^t)$. The maximum principle tells us
that $v=v^t$ and $\tilde v={\tilde v}^t$. Furthermore, let $h$ be
an arbitrary function in $L^2_+(\Omega^t)$ satisfying $h=h^t$, and
let $w$ be the solution of the problem
$\mathbb{B}_1(\Omega^t,c,h)$. Since again $w=w^t$, $w\ge 0$, and
$f\prec^t g$, we
find after partial integration that %
$$  %
\int_{\Omega^t} \tilde v h^t\,dx = \int_{\Omega^t} wf^t\,dx\le
\int_{\Omega^t}wg\,dx=\int_{\Omega^t}vh^t\,dx,  %
$$  %
which means that %
\be  \label{11.30.1-equation} %
\tilde v \prec^t v. %
\ee  %

Next let $\Omega'$ be an open set such that $\supp(f) \subset
\Omega'\subset \overline{\Omega'}\subset \Omega$. By
Corollary~\ref{approximation of open sets by polarization} and
Lemma~\ref{approximation of a simple function by polarization}, we
can find a finite number of polarizations $P_i$ with polarizers
$H_i=\{y>y_i\}$,
$i=1,\ldots,N$, where  %
$$ %\be  \label{11.30.2-equation} %
y_1<y_2<\ldots<y_N\le t,  %
$$   %\ee  %
such that the closure of $\bigcirc_{i=1}^N P_i(\Omega')$ is in
$\Omega^t$ and  $\left(\bigcirc_{i=1}^N P_i f\right)^t=f^t$.

Let $\Omega_m=\bigcirc_{i=1}^m P_i \Omega'$ and
$f_m=\bigcirc_{i=1}^m P_i f$, $m=1,\ldots,N$. Let $u'$ and $u'_m$
be the solutions to the problems $\mathbb{B}_1(\Omega',c,f)$ and
$\mathbb{B}_1(\Omega'_m,c,f_m)$, respectively, $m=1,\ldots,N$.
Applying Theorem~9.1 \cite{BS} we conclude that %
\be  \label{11.30.3-equation} %
u' \prec_{H_1} u'_1\prec_{H_2} u'_2\prec_{H_3}\cdots \prec_{H_N}
u'_N,   %
\ee  %
where ``$\prec_H$'' denotes the partial order related to the
polarization with the polarizer $H$, see \cite[p. 1783]{BS}. Since
$\Omega'_N\subset \Omega^t$ and $f_N=f^t$, the solutions
$u'_N$ and $\tilde v$ satisfy the inequality %
$$ %
0\le u'_N\le \tilde v \quad {\mbox{a.e. in $\Omega_N$.}}  %
$$  %
Together with (\ref{11.30.1-equation}) and
(\ref{11.30.3-equation}) this implies that $u'\prec^t v$.

Now we choose open bounded sets $\Omega^k$ such that
$\overline{\Omega^k}\subset \Omega^{k+1}$, $k=1,2,\ldots$, and
$\cup_k \Omega^k=\Omega$. Let $u^t$ denote the solution of the
problem $\mathbb{B}_1(\Omega^k,c,f)$, $k=1,2,\ldots$. By the above
consideration we have $u^k\prec^t v$, $k=1,2,\ldots$. By the
convergence property of elliptic boundary value problems in
varying domains, see Lemma~A \cite{BS}, the sequence $u^k$
converges to $u$ in $L^2(\Omega)$. This proves the assertion in
the case under consideration.

{\bf (2)} Next we assume that $\gamma\equiv 0$ but $f$ is an
arbitrary function in $L^2_+(\Omega)$. Since simple functions are
dense in $L^2_+(\Omega)$, there is a sequence of simple functions
$f^k$, $k=1,2,\ldots$, with compact supports in $\Omega$, such
that $f^k\to f$ in $L^2(\Omega)$. Then we can find open sets
$\Omega^k$, $k=1,2,\ldots$, such that $\supp(f^k)\subset
\Omega^k\subset \overline{\Omega^k}\subset \Omega^{k+1}\subset
\Omega$ and $\cup_k \Omega^k=\Omega$.

Let $u^k$ and $v_k$ be the solutions to the problems
$\mathbb{B}_1(\Omega^k,c,f^k)$ and
$\mathbb{B}_1((\Omega^k)^t,c,(f^k)^t)$, respectively,
$k=1,2,\ldots$. By part {\bf (1)}, we have  %
\be \label{11.30.4-equation}  %
u_k\prec^t v_k, \quad k=1,2,\ldots  %
\ee %
By Lemma~\ref{Lemma - Approximating sequence in p-norm},
$(f^k)^t\to f^t$ in $L^2(\Omega^t)$ and also by
the monotonicity property of rearrangements,  %
$$  %
\closure((\Omega^k)^t)\subset (\Omega^{k+1})^t \subset \Omega^t
\quad {\mbox{and \ \ }} \cup_k (\Omega^k)^t=\Omega^t. %
$$  %
Therefore by the convergence property of elliptic boundary value
problems in varying domains, see Lemma~A \cite{BS}, we have  %
$$  %
u_k\to u \ \ \ {\mbox{ in $L^2(\Omega)$}} \quad {\mbox{and $v_k\to
v$ in $L^2(\Omega^t)$.}}  %
$$
This together with (\ref{11.30.4-equation}) proves the theorem in
the case $\gamma\equiv 0$.

\smallskip

{\bf (3)} Next let $\gamma\not\equiv 0$. According to Remark~9.2
in \cite{BS} we approximate $u$ and $v$ by solutions to the
problems $\mathbb{B}_1(\Omega,c,\gamma(u_{m-1})+f)$ and
$\mathbb{B}_1(\Omega^t,c,\gamma(v_{m-1})+g)$, respectively,
$m=1,2,\ldots$. Here $u_0\equiv v_0\equiv 0$.  Assume that we had
proved that $u_m\prec^tv_m$ for some $m$. Notice that for $m=0$
this is trivial. Then by equation (\ref{11.4-equation}) of
Theorem~\ref{11-1-Theorem}, we obtain that $\gamma(u_m)+f\prec^t
\gamma(v_m)+g$. By parts {\bf {(1)}} and {\bf{(2)}} of this proof
this means that also $u_{m+1}\prec^t v_{m+1}$, and we conclude by
induction.
  \hfill  $\Box$

 \medskip

Theorem~\ref{11-2-Theorem} and equivalence relations of
Theorem~\ref{11-1-Theorem} lead to the following corollary. %

\bc \label{11.2-Corollary}  %
 Let $u$ and $v$ be solutions defined in
 Theorem~\ref{11-2-Theorem}. Then for every Young function $j$, %
 \be \label{11.24-equation}  %
 \int_\Om j(u)\,dx \le \int_{\Om^t} j(v)\,dx,  %
 \ee  %
if the above integrals converge. In particular, %
 \be  \label{11.25-equation} %
 \|u\|_p\le \|v\|_p \quad \quad {\mbox{for all \ $1\le p\le
 \infty$.
 }}  %
 \ee  %
  \ec  %

\medskip

One might ask under which conditions the equality holds in
inequalities~(\ref{11.24-equation}) and (\ref{11.25-equation}),
and believe that equality is possible only --- roughly speaking
--- if $\Om$ possesses a partial symmetry. For Steiner
symmetrizations this result was proved in \cite{BS}. To prove this
uniqueness result  for the continuous $(k,n)$-Steiner
symmetrization, we restrict
ourselves to the case where $\Om$ is a domain. %

\bt \label{11-3-Theorem}  %
Let $(\cdot)^t$, $0\le t\le 1$, denote the continuous
$(k,n)$-Steiner symmetrization into some $(k,n)$-Steiner
symmetrization with respect to a plane $\Sigma$. Let $\Om$,
$\Omega^t$, $c$, $f$, $g$, $\gamma$, $u$, and $v$ be as in
Theorem~\ref{11-2-Theorem} and assume that $\Om$ is a bounded
domain and $f>0$ on $\Omega$. Assume that there is some Lipschitz
continuous Young function $j$
which for some $t$, $0\le t\le 1$, satisfies %
\be  \label{11.26-equation} %
\int_\Om j(u)\,dx =\int_{\Om^t} j(v)\,dx >0. %
\ee  %
Then  $\Om=\Om^t$ and $f=g$ modulo some translation in a
direction orthogonal to $\Sigma$.  %
\et  %

{\it Proof}. For  $(k,n)$-Steiner symmetrizations this result is
proved in Theorem~10.3 \cite{BS}. Thus, thanks to the standard
inductive argument, we have to prove the theorem for the SC
$1$-symmetrization only. Then $-\infty<t<\infty$ and we may assume
that $\Sigma=\{y=0\}$. Here $x=(x',y)\in \mathbb{R}^{n-1}\times
\mathbb{R}$.

Assume that for some fixed $t\in \mathbb{R}$, $\Omega^t$ is not a
translation of $\Omega$ in the direction of $y$-axis. It follows
from Corollary~\ref{Corollary 6.3} that we can find a polarizer
$H=\{y> t_0\}$ with $t_0<t$ such that either
$(\Omega_H)^t=\Omega^t$, $(f_H)^t=f^t$, $\Omega\not=\Omega_H$, and
$\sigma_H(\Omega)\not=\Omega_H$ or $(\Omega_H)^t=\Omega^t$,
$(f_H)^t=f^t$, $\Omega=\Omega_H$, $f\not= f_H$, and
$\sigma_H(f)\not=f_H$. Then, if $w$ is the solution to the problem
$\mathbb{B}_2(\Omega_H,c,\gamma,f_H)$, we conclude by Theorem~9.3
in \cite{BS} that  %
\be  \label{11.30.5-equation} %
\int_\Omega j(u)\,dx < \int_{\Omega_H} j(w)\,dx.  %
\ee  %
Further, since $(\Omega_H)^t=\Omega^t$ and $(f_H)^t=f^t$, we also
have $w\prec^t v$. By Corollary~\ref{11.2-Corollary} this means
that %
$$  %
\int_{\Omega_H} j(w)\,dx \le \int_{\Omega^t}j(v)\,dx, %
$$  %
which together with (\ref{11.30.5-equation}) contradicts
(\ref{11.26-equation}).

Now we assume that for some $t\in \mathbb{R}$,  $\Omega=\Omega^t$
modulo translation in the direction of the $y$-axis. Without loss
of generality we may assume that $\Omega=\Omega^t$. Assume in
addition that $f\not=f^t$. By Corollary~\ref{Corollary 6.3}, there
is a polarizer $H=\{y> t_1\}$ with $t_1<t$ such that
$\Omega_H=\Omega$ and $f\not=f_H$, and we can argue as before to
derive a contradiction to (\ref{11.26-equation}).

\smallskip

Thus, we have $f=f^t$ and it remains to show that $f^t=g$.

Assume that $f^t\not=g$. We set $\tilde f=\gamma(u)+f$ and $\tilde
g=\gamma(v)+g$. Since $u=u^t$, $v=v^t$ and $u\prec^t v$, it
follows from Theorem~\ref{11-1-Theorem} that %
$$  %
\tilde f \prec^t \tilde g, \quad \tilde f={\tilde f}^t
\not={\tilde g}^t=\tilde g,  %
$$  %
and also  %
\be  \label{11.30.6-equation} %
\int_{2t-y}^y \left( \tilde f(x',s)-\tilde g(x',s)\right)\,ds\le 0  %
\ee  %
for all $x'\in \mathbb{R}^{n-1}$ and every $y\ge t$, where the
inequality~(\ref{11.30.6-equation}) must be strict on a subset of
$\Omega$ of positive measure. Now let $h$ be an arbitrary function
in $L^2_+(\Omega)$ satisfying $h=h^t\not\equiv 0$. Then, if $w$ is
the solution to the problem $\mathbb{B}_1(\Omega,c,h)$, we
conclude that $w=w^t$. Moreover, the strong maximum principle
yields  %
\be  \label{11.30.7-equation} %
\left|\frac{\partial}{\partial y}w(x)\right|>0 \quad {\mbox{a.e.
in $\Omega$.}}  %
\ee  %

Now after partial integration and by using
(\ref{11.30.6-equation}) and (\ref{11.30.7-equation}) we obtain  %
\begin{alignat}{10}  \label{11.30.8-equation}  %
\int_\Omega (u-v)h\,dx &=&\int_\Omega w(\tilde f-\tilde g)\,dx \nonumber \ \ \ \ \ \ \ \ \ \ \ \ \ \ \ \ \ \ \ \ \ \ \ \ \ \ \ \ \ \ \ \ \ \ \ \   \ \ \ \ \ \ \ \ \ \ \ \ \ \ \ \ \ \ \ \ \ \ \   \\
&=& \frac12 \int_\Omega \left|\frac{\partial}{\partial y}
w(x',y)\right|\left(\int_{2t-y}^y(\tilde f(x',\tau)-\tilde
g(x',\tau))\,d\tau\right)\,dx'dy<0. \nonumber %
\end{alignat}  %
Since $j'$ is nondecreasing and $u=u^t$, we have $j'(u)=(j'(u))^t$
(see equation (3.6) in \cite{BS}), and in view of
(\ref{11.26-equation}) it follows that $j'(u)\not=0$. Therefore we
may take $h=j'(u)$ in the equation above. Because of the
convexity of $j$ we get then %
$$  %
\int_\Omega (j(u)-j(v))\,dx \le \int_\Omega j'(u)(u-v)\,dx <0,  %
$$  %
a contradiction.  The theorem is proved.  \hfill   $\Box$

\medskip

Next we prove that solutions of the above considered problems are
continuous from the left with respect to the parameter
 of symmetrization $t$.%

\bt \label{11-4-Theorem}  %
Let $(\cdot)^t$, $0\le t\le 1$, denote the continuous
$(k,n)$-Steiner symmetrization into some $(k,n)$-Steiner
symmetrization. Let $\Om,f,c,\gamma$ be defined as in
Theorem~\ref{11-2-Theorem} and let $t_m$, $m=1,2,\ldots $ be an
increasing sequence in the standard interval $I=[0,1]$
such that  $t_m \to t_0\in I$.  %
Further, let $v$ and $v_m$ be the positive minimal solutions to
the problems \  $\mathbb{B}_2(\Om^{t_0},c,\gamma,f^{t_0})$  and
\ $\mathbb{B}_2(\Om^{t_m},c,\gamma,f^{t_m})$, respectively. Then %
\be  \label{11.22-equation}  %
 v_m \longrightarrow v \qquad {\mbox{in  \ $W
 ^{1,2}({\mathbb{R}}^n)$.}} %
\ee  %
\et  %

{\it Proof}. {\bf (a)} First we prove the theorem for the SC
$1$-symmetrization. Then $-\infty<t<\infty$. Let $t_m$ be an
increasing sequence such that $t_m\to t_0\in \mathbb{R}$.  Since
$(f^{t_m})^{t_0}=f^{t_0}$ for $m=1,2,\ldots$, we may apply
Theorem~\ref{10-4-Theorem} to conclude that the functions $v_m$,
$m=1,2,\ldots$ are equibounded in $W^{1,2}({\mathbb{R}}^n)$.
Therefore we can find a subsequence $v_{m'} $ and a function $w\in
W^{1,2} ({\mathbb{R}}^n )$, such that  %
\be \label{11.22.1-equation} %
 v_{m'} \rightharpoonup  w \qquad
\mbox{weakly in } \ W^{1,2} ({\mathbb{R}}^n )  %
\ee  %
 and  %
\be \label{11.23-equation}  %
 v_{m'}\longrightarrow w \quad \mbox{
in } \ L^2 ({\mathbb{R}}^n ) \ \mbox{ and a.e. }   %
\ee   %
 By Lemma~\ref{Stability Lemma}, we have $\partial \Omega ^{t_0} \subset \partial
\Omega ^{t_m } +r_m \overline{B}_1 ,$ where $r_m=t_0-t_m$. Then
 $r_{m+1}\le r_m$ for $m=1,2,\ldots$ and $r_m\to 0$. Since $v_m\equiv 0$ in $\mathbb{R}^n\setminus \Omega^{t_m}$,
we conclude that $w\in W_0 ^{1,2}(\Omega ^{t_0})$. Thus we can
argue as in the proof of Theorem~\ref{11-2-Theorem} to derive that
$w=v$.

 Next since
$(f^{t_m} )^{t_{m+1}} =f^{t_{m+1}}$ we obtain from
Theorem~\ref{11-2-Theorem}  that   %
\be \label{11.000}   %
 v_1 \prec^{t_2 } v_2 \prec ^{t_3 } \ \ldots \ \prec ^{t_0} v.  %
\ee %  %
Together with (\ref{11.23-equation}) and the
equation~(\ref{11.25-equation}) of Corollary~\ref{11.2-Corollary}
this implies that the sequence of norms $\|v_m\|_2^2$ decreases
and %
$$   %
 \Vert v_m \Vert _2 ^2  \to \Vert v\Vert _2 ^2.  %
$$   %
By the uniform convexity of $L^2 ({\mathbb{R}}^n ) $ this means
that  %
\be  \label{11.30.9-equation}  %
 v_m  \longrightarrow v \qquad
\mbox{in } \ L^2 ({\mathbb{R}}^n ). %
\ee   %
 Now (\ref{11.30.9-equation}), combined with the fact that the functions
 $v_m $, $m=1,2,\ldots$ are equibounded in
$W^{1,2} ({\mathbb{R}}^n )$,  implies that  %
$$   %
 v_m \rightharpoonup v \qquad \mbox{weakly in } \
W^{1,2} ({\mathbb{R}}^n ).  %
$$   %
 Finally we have,
$$  %\begin{eqnarray*}
 \Vert \nabla v_m \Vert _2 ^2  +c\Vert v_m \Vert _2 ^2 =
\int\limits_{\Omega^{t_m } } \Big( \gamma (v_m )+f\Big) v_m \ dx
 \to \int\limits_{\Omega^{t_0}} \Big( \gamma
(v)+f\Big) v \ dx = \Vert \nabla v\Vert _2 ^2 +c\Vert v\Vert _2
^2.  %
$$   %\end{eqnarray*}
In view of the uniform convexity of $W^{1,2}
({\mathbb{R}}^n )$
this yields (\ref{11.22-equation}).

{\bf (b)}  Now we prove the theorem for the continuous
$(k,n)$-Steiner symmetrization and $t_0=1$. In this case,
$\Omega^{t_0}=\Omega^*$ and $f^{t_0}=f^*$; i.e. $\Omega^{t_0}$ and
$f^{t_0}$ are just the corresponding $(k,n)$-Steiner
symmetrizations of $\Omega$ and $f$, respectively. Similarly, we
have $(\Omega^{t_m})^*=\Omega^*$ and $(f^{t_m})^*=f^*$ for
$m=1,2,\ldots$ Applying Theorem~10.1 \cite{BS}, we conclude as
before that the functions $v_m$, $m=1,2,\ldots$ are equibounded in
$W^{1,2}(\mathbb{R}^n)$. Therefore we can find a subsequence
$v_{m'} $ and a function $w\in W^{1,2} ({\mathbb{R}}^n )$, such
that (\ref{11.22.1-equation}) and (\ref{11.23-equation}) remain
valid. By equation~(\ref{9.5-equation}) of
Theorem~\ref{9-2-Theorem}, we have $\partial \Omega^{t_0}\subset
\partial \Omega^{t_m}+\varepsilon_m \overline{B}_1$ with some
$\varepsilon_m>0$ such that $\varepsilon_m\to 0$ as $m\to \infty$.
Since $v_m\equiv 0$ in $\mathbb{R}^n\setminus \Omega^{t_m}$, we
conclude that $w\in W_0^{1,2}(\Omega^{t_0})$. As in part
\textbf{(a)}, this implies that $w\equiv v$.

According to Definition~\ref{continuous $(k,n)$-symmetrization},
$\Omega^{t_{m+1}}$ and $f^{t_{m+1}}$ are obtained from
$\Omega^{t_m}$ and $f^{t_m}$, respectively, after a finite number
of $(k-1,n)$-Steiner symmetrizations followed by an appropriate SC
$(k-1,n)$-Steiner symmetrization. Therefore, it follows from
Theorem~10.1 in \cite{BS} and Theorem~\ref{11-2-Theorem} of this
section that the relations (\ref{11.000}) remain valid in this
case as well. As in part \textbf{(a)} of this proof, the latter
yields~(\ref{11.22-equation}).

{\bf (c)} The theorem is proved for the SC one dimensional
symmetrization and for the case $t_0=1$, which corresponds to the
$(k,n)$-Steiner symmetrization. Now, the general case follows from
these two cases via the standard inductive argument. %
\hfill $\Box$ %

\medskip

Since the proofs of Theorems~\ref{11-2-Theorem} and
\ref{11-4-Theorem} depend only on the maximum principle, we can
derive similar results for parabolic problems. The proofs of
Theorems~\ref{11-5-Theorem} and \ref{11-6-Theorem} below are based
on the approximation scheme involving solutions of some elliptic
problems. This idea was used in \cite{ALT2} for the Schwarz
symmetrization and then in \cite{BS} for the $(k,n)$-Steiner
symmetrization. As we will see this method works also for
continuous symmetrizations. First, following \cite{BS}, we will define solutions of parabolic problems.%

\bd  \label{11-3-definition} %
 Let $\Om\subset \mathbb{R}^n$ be a bounded open set, $c\ge 0$, $T>0$,
 $f\in L_+^2(\Om\times (0,T))$, $\varphi\in L_+^2(\Om)$, and let
 $\gamma:\,\mathbb{R}_0^+\to \mathbb{R}_0^+$ be a globally
 Lipschitzian function. We say that $u$ {\it solves the problem
 $\mathbb{I}(\Om,T,c,\gamma,f,\varphi)$} if $u$ is a solution of
 the following initial boundary value problem:  %
\begin{eqnarray} %  %
& &u\in L^2(0,T;W_0^{1,2}(\Om))\cap C([0,T];L^2(\Om)), \quad \quad
\frac{\partial u}{\partial t}\in L_2([0,T];L^2(\Om)), \nonumber \\ %
& & u_t -\Delta u+cu=\gamma(u)+f
\quad {\mbox{in \ $\Om\times(0,T)$,}} \label{11.27-equation}\\
& & u(x,0)=\varphi(x) \quad {\mbox{in \ $\Om$.}} \nonumber %
\end{eqnarray}  %
\ed  %

\smallskip

Under the assumptions of Definition~\ref{11-3-definition} the
problem $\mathbb{I}(\Om,T,c,\gamma,f,\varphi)$ has a unique
nonnegative solution that can be approximated by the so-called
{\it method of discretization in time}, see \cite{Kac} or
\cite[Section 10]{BS}. To define this approximation, we choose
$N\in\mathbb{N}$. Then we divide the interval $(0,T)$ into $N$
subintervals $[t_{i-1},t_i]$,
where $t_i=iT/N$, and we set   %
\be  \label{11.30.10-equation}  %
f_i(x)=\frac{T}{N} \int_{t_{i-1}}^{t_i} f(x,s)\,ds, \quad
i=1,\ldots,N.  %
\ee  %

We put $u_0=\varphi$. Then let $u_i$ be the solution to the problem %
\be  \label{11.30.11-equation} %
\mathbb{B}_1(\Omega,c+(N/T),\gamma(u_{i-1})+f_i+(N/T)u_{i-1})  %
\ee  %
defined inductively for all $i=1,\ldots,N$.

Let $u^N(x,t)$ denote the function of $x\in \Omega$ and
$t\in[0,T]$ defined for $t_{i-1}\le t\le t_i$, $i=1,\ldots,N$ by  %
\be  \label{11.30.12-equation}  %
u^N(x,t)=u_{i-1}(x)+(t-t_{i-1})(N/T)(u_i(x)-u_{i-1}(x)).  %
\ee  %
The latter equation gives the desired approximation to the
solution of the problem~(\ref{11.27-equation}). Namely we have
(see \cite[Theorem 2.2.4,p.42 ff.]{Kac}):  %
$$  %
u^N(\cdot,t)\rightharpoonup u(\cdot,t) \quad {\mbox{weakly in
$W^{1,2}_0(\Omega)$ for all $t\in (0,T)$,}} %
$$  %
$$  %
\frac{\partial u^N}{\partial t} \rightharpoonup \frac{\partial
u}{\partial t} \quad {\mbox{weakly in $L^2([0,T];L^2(\Omega))$, }}  %
$$  %
and %
$$  %
u^N\to u \quad {\mbox{in $C(0,T;L^2(\Omega))$.}}  %
$$ %

\bt  \label{11-5-Theorem} %
Let $(\cdot)^s$, $0\le s\le 1$, denote a continuous
$(k,n)$-Steiner symmetrization. Let $\Om$, $c$, $T$, $f$,
$\gamma$, $\varphi$ and $u$ be as in
Definition~\ref{11-3-definition} and let $g\in
L_+^2(\Om^s\times(0,T))$, $\psi \in L_+^2(\Om^s)$ with
$f(\cdot,t)\prec^s g(\cdot,t)$ and $g(\cdot,t)=(g(\cdot,t))^s$ for
all $t\in(0,T)$, and $\varphi\prec^t \psi$, $\psi=\psi^t$. Let $v$
be the solution of the problem
$\mathbb{I}(\Om^t,T,c,\gamma,g,\psi)$.
Then  %
\be  \label{11.28-equation}  %
u(\cdot,t)\prec^s v(\cdot,t) \quad \quad {\mbox{for all $t\in(0,T)$}}  %
\ee  %
and %
\be \label{11.29-equation}  %
v(\cdot,t)=(v(\cdot,t))^s \quad \quad {\mbox{for all $t\in (0,T)$.}}  %
\ee  %
\et  %

{\it Proof}. Fix $0<s\le 1$. For any $N\in \mathbb{N}$ and $1\le
i\le N$, let $u_i$ be the solution to the
problem~(\ref{11.30.11-equation})  and let $v_i$ be the solution
to the problem~(\ref{11.30.11-equation}) with $\Omega$, $f$, and
$\varphi$ replaced by $\Omega^s$, $g$, and $\psi$, respectively.
By the assumptions we have $u_0\prec^s v_0$. If $N$ is large
enough then $\gamma(\tau)+N\tau/T$ is increasing and convex in
$\tau$. Since $f\prec^s g$, by (\ref{11.30.10-equation}) we also
have $f_i\prec^s g_i$ for $i=1,\ldots,N$.

Applying Theorem~\ref{11-2-Theorem}, we conclude that $u_i\prec^s
v_i$ for $i=1,\ldots,N$. By (\ref{11.30.12-equation}) it follows
that $u^N\prec^s v^N$. Then passing to the limit as $N\to \infty$
we obtain~(\ref{11.28-equation}).

Since $v_i(\cdot,t)=(v_i(\cdot,t))^s$ for all $t\in (0,T)$ and all
$i=1,\ldots,N$ we have $v^N(\cdot,t)=(v^N(\cdot,t))^s$. Passing to
the limit as $N\to \infty$, this gives~(\ref{11.29-equation}).
\hfill  $\Box$

\medskip

The solutions of the above parabolic problems  are continuous from
the left with respect to the parameter of symmetrization $s$.

\bt \label{11-6-Theorem}  %
Let $\Om$, $c$, $T$, $f$,$\gamma$, $\varphi$, and $u$ be as in
Theorem~\ref{11-5-Theorem}
and let $s_m$, $m=1,2,\ldots $ be an increasing sequence in $I=[0,1]$  such that  $s_m \to s_0\in I$.  %
Further, let $v$ and $v_m$ be the positive solutions of the
problems \  $\mathbb{I}(\Om^{s_0},T,c,\gamma,g^{s_0},\psi)$  and
\ $\mathbb{I}(\Om^{s_m},T,c,\gamma,g^{s_m},\psi))$, respectively. Then %
\be  \label{11.22.2-equation}  %
 v_m \longrightarrow v \qquad {\mbox{in  \ $W
 ^{1,2}({\mathbb{R}}^n)\times (0,T)$.}} %
\ee  %
\et  %

{\it Proof}. %
Using discretization in time as in the proof of
Theorem~\ref{11-5-Theorem}, we approximate $v$ and $v_m$,
$m=1,2,\ldots$ with solutions $v^{(i)}$ and $v_m^{(i)}$ to the
problem (\ref{11.30.11-equation}) for an appropriate initial data.
By Theorem~\ref{11-4-Theorem}, $v_m^{(i)}\to v^{(i)}$ in
$W^{1,2}(\mathbb{R}^n)$. Taking the limit as $N\to \infty$ and
using Theorem~\ref{11-5-Theorem}, we obtain
(\ref{11.22.2-equation}). %
\hfill $\Box$ %

\medskip

\noindent %
\textbf{Remark 11.1.}  %
 The major results in Sections~10 and 11 remain valid if we replace the operator $(-\Delta +c)\ $
 by any uniformly elliptic operator which is invariant under
 considered transformations. In the case of the continuous $k,n)$-Steiner symmetrization
 this is true for instance for operators of the type
$$
-  \sum\limits_{i=1} ^{n-k} \frac{ \partial }{\partial x_i } \Big(
\sum\limits_{j=1}^{n-k} a_{ij} (x')
 \frac{\partial }{\partial x_j } +b_i (x') \Big)
-\sum\limits_{i=n-k+1}^{n} \frac{\partial ^2 }{\partial x_i^2 }
+c(x') ,
$$
where the coefficients $a_{ij} ,b_i $ and $c$  are  bounded and
independent of $y=(x_{n-k+1 } ,\ldots ,x_n)$, $c$ is nonnegative
and
$$
\sum\limits_{i,j=1}^{n-k} a_{ij} (x') \xi _i \xi _j  \geq \lambda
\sum\limits_{i=1}^{n-k} \xi _i ^2 , \qquad \lambda >0.
$$

\section{Appendix}  \label{Appendix}  %
\setcounter{section}{12}  %
\setcounter{equation}{0}  %
In this appendix we prove Theorem~\ref{Sarvas approximation
lemma}. For compact sets $\Omega$ and even $j=2s$, equation
(\ref{3.15-equation}) is a part of Theorem~4.32 in \cite{Sa}. If
$j=2s-1$ is odd, we apply an even number of symmetrizations to the
set $\Omega_1$ and again obtain (\ref{3.15-equation}).

To prove (\ref{3.17-equation}) for $\Omega\in {\mathcal{F}}_n$, we
consider a non-empty slice $\Omega(x')$ that is compact in
$\mathbb{R}^k$. Applying (\ref{3.15-equation}) to symmetrizations
in slices we obtain %
\be  \label{12.1-equation}  %
\lim_{j\to \infty} d(\Omega_j(x'),\Omega^*(x'))=0.  %
\ee  %
Since $\Omega^*(x')$ is a $k$-dimensional ball and
${\mathcal{L}}^k(\Omega_j(x'))={\mathcal{L}}^k(\Omega^*(x'))$,
equation~(\ref{12.1-equation}) implies (\ref{3.17-equation}) in
the case of compact sets.

\smallskip

In the rest of this section, we work with open bounded sets
$\Omega$. Proving (\ref{3.16-equation}), we may assume without
loss of generality that all open sets under consideration belong
to the unit ball $B^{(n)}$.

First we prove three technical lemmas. In all these lemmas, by $S$
we denote the $(k,k)$-Steiner symmetrization in $\mathbb{R}^k$
with respect to the origin $x=(0,\ldots,0)$. Then let $S_1$ and
$S_2$ be $(k-1,k)$-Steiner symmetrizations, which approximate $S$
in the sense of Theorem~\ref{Sarvas approximation lemma}, and let
$\Sigma_1$ and $\Sigma_2$ be the symmetry planes (one dimensional)
of the symmetrizations $S_1$ and $S_2$, respectively. For
notational convenience we will assume  without loss of generality
that $\Sigma_1=\{x=(t\cos \gamma\pi,t\sin \gamma\pi,0,\ldots,0)\in
\mathbb{R}^k:\,-\infty<t<\infty\}$ for some irrational $\gamma\in
(0,1/2)$ and $\Sigma_2$ is the $x_1$-axis.

\smallskip

For given $R$, $x_1$, and $\rho$ such that $R>0$, $-R\le x_1\le
R$, and $0<\rho<R$,
let  %
$$  %
Z(x_1,R,\rho)=\overline{B^{(k)}_R} \cap \left(\cup_\zeta
B^{(k)}_\rho (\zeta)\right),  %
$$  %
where the union is taken over all balls $B^{(k-1)}_\rho (\zeta)$
centered at $\zeta=(\zeta_1,\ldots,\zeta_k)\in \mathbb{R}^k$ such
that $|\zeta|=R$ and $\zeta_1=x_1$. %Thus, $Z(x_1,R,\rho)$ is an intersection of a ball and a solid torus.

\bl  \label{Lemma-12.2}  %
Let $0<r<R_1<1$. %Under the assumptions of Lemma~\ref{Lemma-12.1}
Then there exists $\tau=\tau(r,R_1)$, $0<\tau<1$, such that for
every compact set $K$ satisfying the conditions $K\subset
\overline{B^{(k)}_R}$ for some $R\in [R_1,1]$ and
${\mathcal{L}}^k(K)\le {\mathcal{L}}^k(B^{(k)}_r)$ there is real
$x_1=x_1(K)$ such that
$-r\le x_1\le r$ and %
$$  %
\overline{B^{(k)}_R}\setminus K_2\supset Z(x_1,R,\tau R).  %
$$  %
\el  %

{\it Proof}. Since $K_2$ is Steiner symmetric with respect to
$\Sigma_2$ and since ${\mathcal{L}}^k(K_2)\le
{\mathcal{L}}^k(B^{(k)}_r)$ there is a point
$x_0=(x_1^0,x_2^0,0,\ldots,0)\in S^{(k)}_r\setminus K_2$ such that
$-r\le x_1^0\le r$ and $x_2^0=\sqrt{r^2-(x_1^0)^2}\ge 0$.

Since the slice $K_2(x_1^0)$ is a $(k-1)$-dimensional ball we have
$K_2(x_1^0)\subset B^{(k-1)}_{x_2^0}$. Since at the same time
$K_2(x_1^0)$ is the result of $S_2$ symmetrization of $K_1(x_1^0)$
we have  %
\begin{alignat}{10} \label{12.30-equation} %$$  %
{\mathcal{L}}^{k-1}(B^{(k)}_R(x_1^0)\setminus
K_1)&=&{\mathcal{L}}^{k-1}(B^{(k)}_R(x_1^0)\setminus K_2)\ \ \ \ \ \ \ \ \ \ \ \ \ \ \ \ \ \ \  \\ %
&\ge& {\mathcal{L}}^{k-1}(B^{(k)}_{R_1}(x_1^0)\setminus
B^{(k-1)}_{x_2^0})\ge c(r,R_1), \nonumber %
\end{alignat} %\ee  %$$  %
where %
$$  %
c(r,R_1)=\min_{|x_1^0|\le r}
{\mathcal{L}}^{k-1}(B^{(k)}_{R_1}(x_1^0)\setminus
B^{(k-1)}_{x_2^0})>0.  %
$$  %

For $x\in B^{(k)}_R\setminus \Sigma_1$, let $\Gamma_1(x)$ denote
the ray through $x$ that is orthogonal to $\Sigma_1$ and has its
origin at some point $\tilde x\in \Sigma_1\cap B^{(k)}_R$. Let
$I(x)$ be a closed segment of $\Gamma_1(x)$ joining $x$ and
$\partial B^{(k)}_R$. We note that $I(x)\cap K_1=\emptyset$ if
$x\not\in K_1$. Let %
$$  %
J(x_1^0)=\cup_x I(x),  %
$$  %
where the union is taken over all $x\in B^{(k)}_R(x_1^0)\setminus
K_1(x_1^0)$. Let $A=A(K,x_1^0)$ be the maximal segment of the
$x_1$-axis
such that $x_1^0\in A$ and if $x_1\in A$ then  %
\be  \label{12.2-equation}  %
{\mathcal{L}}^{k-1}(\mathbb{R}^k(x_1)\cap J(x_1^0))\ge
(1/2)c(r,R_1).  %
\ee  %
Using (\ref{12.30-equation}), one can show that there is a constant $t=t(r,R_1)>0$ such that %
\be  \label{12.3-equation} %
{\mathcal{L}}^1(A(K,x_1^0))\ge t(r,R_1)  %
\ee  %
for  every compact set $K$ satisfying the assumptions of the
lemma.

Finally, it is not difficult to see that (\ref{12.2-equation}) and
(\ref{12.3-equation}) imply the lemma.  \hfill $\Box$

\smallskip

The following lemma allows us to control how fast the
approximation process is for compact sets.

\bl  \label{Lemma-12.1}  %
 Let $0<r<R_1<1$. There exist a positive integer $N=N(r,R_1)$ and real
$\beta=\beta(r,R_1)$, $0<\beta<1$,  such that  for all $j\ge N$ and all $R$ such that  $R_1\le R\le 1$, we have %
$$  %
K_j\subset \overline{B^{(k)}_{\beta R}}  %
$$  %
for every compact set $K\subset \overline{B^{(k)}_R}$ such that
${\mathcal{L}}^k(K)\le {\mathcal{L}}^k(B^{(k)}_r)$.  %
\el  %

{\it Proof}. %
We may assume without loss of generality that
\be \label{12.31-equation}  %
K\subset \overline{B^{(k)}_R}\setminus Z(x_1,R,\tau R)  %W(x_1):=
\ee   %
with $\tau=\tau(r,R_1)>0$ defined in Lemma~\ref{Lemma-12.2} and
some $x_1=x_1(K)$ such that $-r\le x_1\le r$. Indeed, if $K$ does
not satisfy (\ref{12.31-equation}), then we replace $K$ with its
second symmetrization $K_2$, which by Lemma~\ref{Lemma-12.2}
satisfies the required condition.

Let $C_R$ be the circle in the %two dimensional
plane
${\mathcal{E}}=\{x\in \mathbb{R}^k:\, x=(x_1,x_2,0,\ldots,0)\}$.
Using complex notation, we will write $Re^{i\theta}$ for
$x=(R\cos\theta,R\sin\theta,0,\ldots,0)\in C_R$.

It follows from (\ref{12.31-equation}) that there is a
sufficiently small constant $a_0=a_0(r,R_1)>0$ such that for every
compact set $K$ satisfying the conditions of the lemma there is a
point $Re^{i\theta_0}$ with $\theta_0=\theta_0(K)$ such that the
intersection $C_R\cap Z(x_1,R,\tau R)$ contains an arc $\alpha$
centered at $Re^{i\theta_0}$ with the angular measure $\ge a_0$.

In addition, it is not difficult to see that there exists a
constant $\nu=\nu(r,R_1)$, $0<\nu<1$, such that the set $C_R\cap
Z(x_1,R,\nu\tau R)$ contains a subarc $\alpha'\subset \alpha$
again centered at $Re^{i\theta_0}$, which has the angular measure
$\ge (2/3)a_0$.  %

Let ${\rm{Ref}}_1$ and ${\rm{Ref}}_2$ denote the reflections in
the plane ${\mathcal{E}}$ with respect to the lines
$l_1=\{z=te^{i\gamma\pi}:\,-\infty<t<\infty\}$ and
$l_2=\{z=t:\,-\infty<t<\infty\}$, respectively. Let
${\mathcal{R}}^s=({\rm{Ref}}_1\circ {\rm{Ref}}_2)^s$,
$s=2,3,\ldots$. Then
${\mathcal{R}}^s(Re^{i\theta_0})=Re^{i(2s\gamma\pi +\theta_0)}$,
$s=2,3,\ldots$.

Since $\gamma$ is irrational, the set
$\{Re^{2is\gamma\pi}\}_{s=2}^\infty$ is dense in $C_R$; see, for
example, \cite[Lemma 3.25]{Sa}. Therefore there are indices
$s_1,\ldots,s_{N_1}$ such that for any $\theta_0$, the points
$Re^{i(2s_j\gamma\pi+\theta_0)}$, $j=1,\ldots,N_1$, divide $C_R$
into $N_1$ arcs, each of which has the angular measure $\le
a_0/4$. Let $N=\max\{s_1,\ldots,s_{N_1}\}$. Let %
$$  %
\alpha_s={\mathcal{R}}^s(\alpha), \quad
\alpha'_s={\mathcal{R}}^s(\alpha'), \quad s=2,3,\ldots,N,  %
$$  %
where $\alpha$ and $\alpha'$ are the arcs defined above. By our
choice of the points $Re^{i(2s_j\gamma\pi+\theta_0)}$, we
have  %
$$  %
C_R=\cup_{s=1}^N \alpha'_s=\cup_{s=1}^N \alpha_s.  %
$$  %

Let %
$$  %
W(x_1)=\overline{B^{(k)}_R}\setminus Z(x_1,R,\tau R), \quad
W'(x_1)=\overline{B^{(k)}_R}\setminus Z(x_1,R,\nu\tau R),  %
$$  %
where $x_1=x_1(K)$ and  $\tau$ and $\nu$ are positive constants
defined in the
beginning of this proof. For $s=1,2,\ldots$ and $x_1=x_1(K)$, let  %
$$  %
W_s(x_1)=(S_2\circ S_1)^s(W(x_1)), \quad W'_s(x_1)=(S_2\circ
S_1)^s(W'(x_1)).  %
$$  %

It follows from (\ref{12.31-equation}) and our construction that %
$$  %
K_{2N}\subset W_N(x_1)\subset W'_N(x_1)  %
$$  %
and %
\be  \label{12.6-equation}  %
W'_N(x_1)\cap S^{(k)}_R=\emptyset.  %
\ee  %

Let %
$$  %
\beta(x_1)=\inf \{\beta>0:\, W_N(x_1)\subset
\overline{B^{(k)}}_{\beta R}\}.  %
$$  %
By (\ref{12.6-equation}), $0<\beta(x_1)<1$ for all $x_1$ such that
$-r\le x_1\le r$. Let $\tilde \beta=\sup_{|x_1|\le r} \beta(x_1)$.
Our construction above works for every compact set $K$ satisfying
the assumptions of the lemma. Choosing $K=\overline{B_r^{(k)}}$,
we obtain that  $\tilde \beta\ge r/R\ge r>0$.

Thus, to complete the proof, we have to show that $\tilde\beta<1$.
If not, we can find a sequence $x^{(s)}_1\in [-r,r]$,
$s=1,2,\ldots$, such that $x^{(s)}_1\to x^{(0)}_1$ and
$\beta(x^{(s)}_1)\to 1$ as $s\to \infty$. For all sufficiently
large $s$ we have %
$$  %
Z(x^{(0)}_1,R,\nu\tau R)\subset Z(x^{(s)}_1,R,\tau R).  %
$$  %
Therefore for all such $s$ we have  %
$$  %
W'_N(x^{(0)}_1)\supset W_N(x^{(s)}_1).  %
$$  %
By (\ref{12.6-equation}), $W'_N(x^{(0)}_1)\subset
\overline{B^{(k)}_{\beta'R}}$ for some $0<\beta'<1$. Hence for all
sufficiently large $s$  we have  %
$$  %
W_N(x^{(s)}_1)\subset \overline{B^{(k)}_{\beta'R}} %
$$  %
contradicting the assumption that $\beta(x^{(s)})\to 1$ as $s\to
\infty$. The proof of the lemma is complete.  \hfill  $\Box$

%\medskip

\bl  \label{Lemma-12.3}  %
Let $S$, $S_1$, and $S_2$ be the symmetrizations as  in
Lemmas~\ref{Lemma-12.2} and \ref{Lemma-12.1}. Let $0<r<R$ and
$0<\rho <R$. Then there is
a constant $c_1=c_1(r,\rho,R)>0$ such that %
\be   \label{12.3.1-equation}  %
{\mathcal{L}}^k(\Omega_1\setminus B^{(k)}_r)\ge c_1  %
\ee  %
for every open set $\Omega\subset \mathbb{R}^k$ such that
$\Omega_2$ contains some point $x_0=(x_0',y_0)$ with $x'_0\in
\mathbb{R}$, $y_0\in \mathbb{R}^{k-1}$ such that $|x_0|\ge R$ and
$|y_0|\ge \rho$.  %
\el  %

{\it Proof}.  Let $\Omega$ satisfies the assumptions of the lemma
and let %
$$  %
c_2=c_2(r,\rho,R)=\min_{x'_0} {\mathcal{L}}^{k-1}((\{x'_0\}\times
B^{(k-1)}_\rho)\setminus B^{(k)}_r),  %
$$  %
where the minimum is taken over all $x'_0$ such that
$(x'_0)^2+\rho^2\ge R^2$. It is clear that $c_2>0$. Since
$\Omega_2(x'_0)$ is a $(k-1)$-dimensional ball with the radius
$\ge \rho$, we have ${\mathcal{L}}^{k-1}(\Omega_2(x'_0)\setminus
B^{(k)}_r)\ge c_2$.

For $x\not\in \Sigma_1$, let $\hat I(x)$ denote the segment
orthogonal to $\Sigma_1$ that joins $x$ and $\Sigma_1$.  For $x\in
\mathbb{R}^k\setminus (B^{(k)}_r\cup \Sigma_1)$, we set $I(x)=\hat
I(x)\setminus B^{(k)}_r$. We note that $I(x)\subset \Omega_1$,
whenever $x\in \Omega_1$. For $x'_0\in \mathbb{R}$ defined above, let  %
$$  %
J_1(x'_0)=\cup_x I(x),  %
$$  %
where the union is taken over all $x=(x'_0,y)\in \mathbb{R}^k
\setminus B^{(k)}_r$ such that $y\in \Omega_1(x'_0)$. Let
$A_1(\Omega,x'_0)$ be the maximal closed segment of the $x_1$-axis
such
that $x'_0\in A_1(\Omega,x'_0)$ and %
\be  \label{12.4-equation}  %
{\mathcal{L}}^{k-1}((J_1(x'_0))(x_1))\ge c_2/2,  %
\ee   %
where $(J_1(x'_0))(x_1)$ denotes the $(k-1)$-slice of $J_1(x'_0)$
at $x_1$.  It is not difficult to see that there is a constant
$t_1=t_1(r,\rho,R)>0$ such that  %
\be  \label{12.5-equation}  %
{\mathcal{L}}^1(A_1(\Omega,x'_0))\ge t_1  %
\ee  %
for every domain $\Omega$ satisfying the assumptions of the lemma.

Since $J_1(x'_0)\subset \Omega_1$, equations (\ref{12.4-equation})
and (\ref{12.5-equation}) imply~(\ref{12.3.1-equation}).  \hfill
$\Box$

\medskip

{\it Proof of Theorem~\ref{Sarvas approximation lemma} for open
sets}. %
The proof is by contradiction. Suppose that
equation~(\ref{3.16-equation}) does not hold for some open set
$\Omega\subset B^{(k)}$.  Then we have to consider the following
two cases.

{\bf (i)} There are real number $\varepsilon_0>0$ and a sequence
$x_s\to x_0$ such that $x_s\in \partial \Omega_{m_s}$ for some
subsequence of indices $m_s$ and %
\be  \label{12.7-equation}  %
d(x_s,\partial \Omega^*)\ge \varepsilon_0 \quad {\mbox{for
$s=0,1,\ldots$}}  %
\ee  %

{\bf (ii)}  There are real number $\varepsilon_0>0$ and a sequence
$x_s\to x_0$ such that $x_s\in \partial \Omega^*$ for all $s$ and
there is a subsequence of indices $m_s$ such that %
\be  \label{12.8-equation}  %
d(x_s,\partial \Omega_{m_s})\ge \e_0 \quad {\mbox{for
$s=1,2,\ldots$}}  %
\ee  %

In each of these cases we may assume without loss of generality
that all indices $m_s$ are odd and $m_s\ge 3$. Then
$\Omega_{m_s}=S_1(\Omega_{m_s-1})$.

To prove {\bf (i)}, we consider two subcases.  %

{\bf (1)}  We first assume that $B^{(n)}_{\e_0}(x_0)\subset
\Omega^*$. Let $x_s=(x'_s,y_s)$, $s=0,1,\ldots$ Then there are
constants $\delta_0>0$ and $R>|y_0|$ such that %
\be  \label{12.9-equation}  %
\Omega^*(x') \supset B^{(k)}_R \quad {\mbox{for all $x'\in
B^{(n-k)}_{\delta_0}(x'_0)$.}}  %
\ee  %

Next we  estimate how much the slice $\Omega_{m_s-1}(x'_s)$
differs from the ball $B^{(k)}_R$. Let $L_1(x',y)$ and $L_2(x',y)$
be the symmetrizing planes through $(x',y)\in
\mathbb{R}^{n-k}\times \mathbb{R}^k$ of $S_1$ and $S_2$,
respectively. Then $L_j(x',y)$ can be represented in the form
$L_j(x',y)=\{(x',t):\, t\in L_j(y)\}$, where $L_j(y)$ is a
$(k-1)$-dimensional plane in $\mathbb{R}^{k}$,  which does not
depend on $x'$.

Since $y_s\to y_0$ and $R>|y_0|$, we have that
$|y_s|<R_1=(R+|y_0|)/2$ for all sufficiently large $s$. Since
$\Omega_{m_s}$ is the $(k-1)$-dimensional Steiner symmetrization
of $\Omega_{m_s-1}$ with respect to $S_1$, the set
$\Omega_{m_s}(x'_s)\cap L_1(y_s)$ is a $(k-1)$-dimensional ball in
the corresponding $(k-1)$-dimensional plane. Since $(x'_s,y_s)\in
\partial
\Omega_{m_s}$, the latter implies that %
$$  %
\Omega_{m_s}(x'_s)\cap L_1(y_s)\subset B_{R_1}^{(k)}\cap
L_1(y_s)\subset B_R^{(k)}\cap L_1(y_s). %
$$ %
This implies that there is $\delta_1=\delta_1(R,|y_0|)>0$ such
that %
$$  %
{\mathcal{L}}^{k-1}((B_R^{(k)}\cap L_1(y_s))\setminus
(\Omega_{m_s}(x'_s)\cap L_1(y_s)))\ge \delta_1  %
$$ %
for all sufficiently large $s$. Since
$\Omega_{m_s}=S_1(\Omega_{m_s-1})$, the latter inequality implies that %
\be  \label{12.11-equation} %
{\mathcal{L}}^{k-1}((B_R^{(k)}\cap L_1(y_s))\setminus
(\Omega_{m_s-1}(x'_s)\cap L_1(y_s)))\ge \delta_1  %
\ee  %
for all sufficiently large $s$.

Now since $\Omega_{m_s-1}(x'_s)$ omits a set of positive
$(k-1)$-dimensional measure in $B_R^{(k)}\cap L_1(y_s)$ and at the
same time $\Omega_{m_s-1}(x'_s)\cap L_2(y_s)$ is a
$(k-1)$-dimensional ball in $L_2(y_s)$, we can argue as in the
proof of Lemma~\ref{Lemma-12.3} (cf. how (\ref{12.3.1-equation})
follows from (\ref{12.4-equation}) and (\ref{12.5-equation})) to
deduce that there is $\delta_2=\delta_2(R,|y_0|)>0$ such that %
\be  \label{12.12-equation}  %
{\mathcal{L}}^k(B_R^{(k)} \setminus \Omega_{m_s-1}(x'_s))\ge
\delta_2  %
\ee  %
for all sufficiently large $s$.

Next we show that (\ref{12.12-equation}) leads to a contradiction.
The slice $\Omega(x'_0)$ is an open set in $B^{(k)}$. Therefore
for arbitrary small $\e>0$ there is a compact set $F\subset
\Omega(x'_0)$ such that %
\be \label{12.13-equation}  %
{\mathcal{L}}^k(\Omega(x'_0)\setminus F)<\e. %
\ee  %
Then for sufficiently small $\e_1>0$, we have
$B^{(n-k)}_{\e_1}(x'_0)\times F\subset \Omega$.
By the monotonicity property of symmetrizations, we have %
$$  %
((\{x'_0\}\times F)_{m_s-1})(x'_0)\subset \Omega_{m_s-1}(x'_s),  %
$$  %
where $((\{x'_0\}\times F)_{m_s-1})(x'_0)$ denotes the slice of
$(\{x'_0\}\times F)_{m_s-1}$ at $x'_0$. Applying
equation~(\ref{3.17-equation}) of Theorem~\ref{Sarvas
approximation lemma} to the compact sets $\{x'\}\times F$ and
$\{x'\}\times S(F)$, we obtain the following limit relation for measure in slices    %
\be  \label{12.14-equation}  %
{\mathcal{L}}^k((((\{x'\}\times F)_j)(x')) \bigtriangleup
(((\{x'\}\times
S(F)))(x')))\to 0 \quad {\mbox{as $j\to \infty$.}}  %
\ee  %

Now, (\ref{12.9-equation}), (\ref{12.13-equation}), and
(\ref{12.14-equation}) imply that %
$$  %
{\mathcal{L}}^k(B^{(k)}_R \setminus ((\{x'_s\}\times
F)_{m_s-1})(x'_s))\le 2\e  %
$$  %
for all sufficiently large $s$. Since $((\{x'_s\}\times
F)_{m_s-1})(x'_s)\subset \Omega_{m_s-1}(x'_s)$, the latter
inequality contradicts (\ref{12.12-equation}). This completes the
proof of the theorem in the case under consideration.  %

\smallskip

{\bf (2)}  In the second case, we assume that
$B^{(n)}_{\e_0}(x_0)\cap \Omega^*=\emptyset$. We recall that
$x_s=(x'_s,y_s)\to x_0=(x'_0,y_0)$ as $s\to \infty$.

If $y_0=(0,\ldots,0)$, then $\Omega^*(x')=\emptyset$ for all
$x'\in \mathbb{R}^{n-k}$ sufficiently close to $x'_0$. Hence,
$\Omega(x')=\emptyset$ and therefore $\Omega_j(x')=\emptyset$ for
all such $x'$ and all $j=1,2,\ldots$ It is easily seen that the
latter contradicts our assumptions that $x_s\in\partial
\Omega_{m_s}$ and $x_s\to x_0$ as $s\to \infty$.

\smallskip

Assume now that $|y_0|=R_0>0$. Since $B_{\e_0}^{(n)}(x_0)\cap
\Omega^*=\emptyset$ there is $\delta_0>0$ such that for every
$x'\in\overline{B_{\delta_0}^{(n-k)}}(x'_0)$ the slice
$\Omega^*(x')$ is either an empty set or an open $k$-dimensional
ball $B_{R(x')}^{(k)}$ with the radius $R(x')$ such that %
\be \label{12.32-equation}  %
0<R(x')\le \rho<R_0  %
\ee  %
with some $\rho$ independent of $x'\in
\overline{B_{\delta_0}^{(n-k)}}(x'_0)$.
In particular, (\ref{12.32-equation}) shows that %
\be  \label{12.15-equation}  %
0<{\mathcal{L}}^k(\Omega(x'))={\mathcal{L}}^k(\Omega^*(x'))\le
{\mathcal{L}}^k(B^{(k)}_\rho)  %
\ee  %
for all $x'\in \overline{B_{\delta_0}^{(n-k)}}(x'_0)$.

Let $\e_0>0$ be fixed and sufficiently small. For every $x'\in
\overline{B_{\delta_0}^{(n-k)}}(x'_0)$, we choose   a
$k$-dimensional compact set $K(x')$ such that $K(x')\subset
\Omega(x')$ and
\be  \label{12.17-equation}  %
{\mathcal{L}}^k(\Omega(x')\setminus K(x'))\le \e_0.  %
\ee  %
By (\ref{12.15-equation}),  we have %
\be  \label{12.18-equation}  %
{\mathcal{L}}^k(K(x'))\le {\mathcal{L}}^k(B^{(k)}_\rho)  %
\ee  %
for all $x'\in \overline{B_{\delta_0}^{(n-k)}}(x'_0)$.

 Let $\e_2>0$ be sufficiently small and let $\e_1>0$ be such that
$\rho+\varepsilon_1<R_0$ and %
\be \label{12.33-equation}  %
{\mathcal{L}}^k(B^{(k)}_{\rho+\e_1}\setminus B^{(k)}_\rho)<\e_2  %
\ee %
for all $\rho \le R_0$.

For $j=1,2,\ldots$ and $x'\in
\overline{B_{\delta_0}^{(n-k)}}(x'_0)$, let $K^j(x')$ denote the
slice at $x'$ of the $j$-th successive symmetrization of the set
$\{x'\}\times K(x')$ defined by formulas (\ref{3.13-equation}) and
(\ref{3.14-equation}). Alternatively, $K^j(x')$ can be obtained by
applying appropriate $(k-1)$-dimensional symmetrizations to the
set $K(x')$ in $\mathbb{R}^k$.

Since $K(x')\subset B^{(k)}$ for all $x'$ and since $K(x')$
satisfies (\ref{12.18-equation}) for all $x'\in
\overline{B_{\delta_0}^{(n-k)}}(x'_0)$, we can apply
Lemma~\ref{Lemma-12.1} with $r=\rho$ and $R_1=\rho+\varepsilon_1$
to the compact sets $K(x')$. This implies that there exists a
positive integer $N= N(\rho,\varepsilon_1)$ such that %
\be \label{12.19-equation} %
K^j(x')\subset B^{(k)}_{\rho+\e_1}  %
\ee  %
for all $x'\in \overline{B_{\delta_0}^{(n-k)}}(x'_0)$ and all
$j\ge N$.

Now, combining (\ref{12.17-equation}), (\ref{12.33-equation}), and
(\ref{12.19-equation}), we obtain that  %
\be  \label{12.16-equation}  %
{\mathcal{L}}^k(\Omega_j(x')\setminus B^{(k)}_\rho)\le
\e_0+\e_2  %
\ee  %
for all $x'\in \overline{B_{\delta_0}^{(n-k)}}(x'_0)$ and all
$j\ge N$.

 \smallskip

 Now we return to the sequence $x_s=(x'_s,y_s)\to x_0=(x'_0,y_0)$.
 Since $x_s\in \partial \Omega_{m_s}$, for every $s$
we can find a point $\hat x_s=(\hat x'_s,\hat y_s)\in
\Omega_{m_s}$ such that $\hat x_s\to x_0$ as $s\to \infty$. Then,
of course, $\Omega({\hat x}'_s)\not=\emptyset$ and therefore
$\Omega^*({\hat x}'_s)\not=\emptyset$.

Suppose first that $x_0\not\in \Sigma_2$. Then
$d=d(x_0,\Sigma_2)>0$. Now we can apply Lemma~\ref{Lemma-12.3}
with $r=\rho$, $R=(R_0+\rho)/2$, and $\rho=d$ and with the domains
$\Omega_1$ and $\Omega_2$ in that lemma replaced by the domains
$\Omega_{m_s-1}({\hat x}'_s)$ and $\Omega_{m_s}({\hat x}'_s)$,
respectively. By Lemma~\ref{Lemma-12.3} there exists a constant
$c_1=c_1(\rho,d,R_0)>0$ such that  %
$$ %
{\mathcal{L}}^k(\Omega_{m_s-1}(x'_s)\setminus B^{(k)}_\rho)\ge c_1  %
$$  %
for all sufficiently large $s$, which obviously contradicts
(\ref{12.16-equation}) if $\e_0$ and $\e_2$ are chosen
sufficiently small. This proves the theorem in the case under
consideration if $x_0\not \in \Sigma_2$.

\smallskip

Suppose now that $x_0\in \Sigma_2$. For every $s=1,2,\ldots$, the
symmetrizing plane $L_2(\hat y_s)$ contains some point $\tilde
x_s=(\hat x'_s,\tilde y_s)\in \partial \Omega_{m_s-1}$. Selecting
a subsequence if necessary, we may assume that $\tilde x_s \to
{\tilde x}^{(1)}=(x'_0,{\tilde y}^{(1)})$.

Since $B^{(n)}_{\e_0}(x_0)\cap \Omega^*=\emptyset$ and $x_0\in
\Sigma_2$, it follows
from the definition of $(k,n)$-Steiner symmetrization that  %
$$  %
B^{(n)}_{\e_0}({\tilde x}^{(1)})\cap \Omega^*=\emptyset.  %
$$  %

Now, if ${\tilde x}^{(1)}\not \in \Sigma_1$, to complete the
proof, we can apply our argument above replacing the plane
$\Sigma_2$ and the sequence of points $x_s\in \partial
\Omega_{m_s}$ with the plane $\Sigma_1$ and the sequence ${\tilde
x}_s\in \partial \Omega_{m_s-1}$, respectively.

In the case ${\tilde x}^{(1)}\in \Sigma_1$, we continue our
construction to find points ${\tilde x}^{(2)}\in \Sigma_2$,
${\tilde x}^{(3)}\in \Sigma_1$, $\ldots$ The sequence of points
${\tilde x}^{(j)}$ will be finite if ${\tilde x}^{(2m)}\not \in
\Sigma_2$ or ${\tilde x}^{(2m-1)}\not\in \Sigma_1$ for some $m\ge
1$. Otherwise, the sequence ${\tilde x}^{(j)}$ will contain
infinite number of terms. If it is finite, say $j=1,\ldots,N$,
then we apply our previous argument to the point ${\tilde
x}^{(N)}$ and to the plane $\Sigma_i$, where $i=1$ if $N$ is odd
and $i=2$ if $N$ is even.

Assume now that the constructed sequence of points ${\tilde
x}^{(j)}=(x'_0,{\tilde y}^{(j)})$, $j=1,2,\ldots$, is infinite. By
our construction we have $ \quad |{\tilde y}^{(j+1)}|\ge |{\tilde
y}^{(j)}|\,\sec(\gamma\pi)$ for all $j=1,2,\ldots$.  Therefore, %
$$  %
|{\tilde x}^{(j)}|\to \infty \quad \quad {\mbox{as $j\to
\infty$.}} %
$$ %
Since for every $j$ there is an index $m(j)$ and a point
$z^{(j)}\in \Omega_{m(j)}$ such that $|z^{(j)}-{\tilde
x}^{(j)}|\le 1$ the latter limit relation contradicts our
assumption that $\Omega\subset B^{(n)}$.

This completes the proof of the theorem in the case {\bf (i)}.

\medskip

In the case {\bf (ii)} the proof is simpler. As in the case {\bf
(i)}, we consider two subcases.  %

{\bf (1)} First, we suppose that $B^{(n)}_{\e_0}(x_0)\cap
\Omega_{m_s}=\emptyset$ for some sufficiently small $\e_0>0$ and
some infinite subsequence of indices $m_s$, $s=1,2,\ldots$

Fix $\delta_0>0$ sufficiently small. Then let $K$ be a compact subset of $\Omega$ such that %
\be  \label{12.20-equation}  %
{\mathcal{L}}^n(\Omega\setminus K)<\delta_0.  %
\ee  %
Let $K^*=S(K)$ and let $K_j$, $j=1,2,\ldots$, be successive
symmetrizations of $K$ defined by (\ref{3.13-equation})
and (\ref{3.14-equation}). Since $K_{m_s}\subset \Omega_{m_s}$, we have %
\be  \label{12.21-equation}  %
B^{(n)}_{\e_0}(x_0)\cap K_{m_s}=\emptyset.
\ee  %
Since $x_0\in \partial \Omega^*$ there is $\delta_1>0$ such that %
$$  %
{\mathcal{L}}^n(B^{(n)}_{\e_0}(x_0)\cap \Omega^*)\ge \delta_1.  %
$$  %
This together with (\ref{12.20-equation}) implies that  %
\be  \label{12.22-equation} %
{\mathcal{L}}^n(B^{(n)}_{\e_0}(x_0)\cap K^*)\ge \delta_2  %
\ee  %
for some $\delta_2>0$ if $\delta_0>0$ in (\ref{12.20-equation}) is
sufficiently small.

Applying equation~(\ref{3.17-equation}) to the compact set $K$, we
obtain  %
\be  \label{12.23-equation} %
{\mathcal{L}}^n(K_{m_s} \bigtriangleup K^*)\to 0 \quad {\mbox{as
$s\to \infty$.}}  %
\ee  %

One can easily see that equations (\ref{12.22-equation}) and
(\ref{12.23-equation}) contradict (\ref{12.21-equation}).  This
proves the theorem in the case under consideration.

\smallskip

{\bf (2)}  Suppose now that $\varepsilon_0>0$ and $x_0=(x'_0,y_0)$ are such that %
\be \label{12.24-equation} %
B^{(n)}_{\e_0}(x_0)\subset \Omega_{m_s} %
\ee  %
for some infinite subsequence $m_s$, $s=1,2,\ldots$ This implies
that, $\Omega_{m_s}(x'_0)\not=\emptyset$ and therefore
$\Omega(x'_0)\not=\emptyset$ and $\Omega^*(x'_0)\not=\emptyset$.
Since $\Omega^*(x'_0)$ is a $k$-dimensional ball there is  $r>0$
such that $\Omega^*(x'_0)=B^{(k)}_r$. Since $x_0\not\in \Omega^*$,
we have
$r\le |y_0|$. This implies that %
\be  \label{12.25-equation} %
{\mathcal{L}}^k((B^{(n)}_{\e_0}(x_0)\setminus \Omega^*)(x'_0))\ge
\delta_1  %
\ee  %
for some $\delta_1>0$. Here $(B^{(n)}_{\e_0}(x_0)\setminus
\Omega^*)(x'_0)$ denotes the slice of
$B^{(n)}_{\e_0}(x_0)\setminus \Omega^*$ at $x'_0$.   For every
arbitrary small $\delta_2>0$ there exists a
$k$-dimensional compact set $K\subset \Omega(x'_0)$ such that %
\be  \label{12.26-equation}  %
{\mathcal{L}}^k(\Omega(x'_0)\setminus
K)={\mathcal{L}}^k(\Omega^*(x'_0)\setminus K^*)<\delta_2,  %
\ee  %
where $K^*$ denotes the $(k,k)$-Steiner symmetrization of $K$.

Applying equation~(\ref{3.17-equation}) of Theorem~\ref{Sarvas
approximation lemma} to $K$, we obtain  %
\be  \label{12.27-equation}  %
{\mathcal{L}}^k(K_s \bigtriangleup K^*)\to 0 \quad {\mbox{as $s\to
\infty$,}}  %
\ee  %
where $K_s=((\{x'_0\}\times K)_s)(x'_0)$ denotes the slice at
$x'_0$ of the symmetrized set $(\{x'_0\}\times K)_s$ defined by
formulas (\ref{3.13-equation}) and (\ref{3.14-equation}).

By the monotonicity property of symmetrization, we have
$K_{m_s}\subset \Omega_{m_s}(x'_0)$. This and (\ref{12.26-equation}) imply that  %
\be  \label{12.28-equation}  %
{\mathcal{L}}^k(\Omega_{m_s}(x'_0)\setminus K_{m_s})<\delta_2.  %
\ee  %

Finally, (\ref{12.24-equation}), (\ref{12.25-equation}), and
(\ref{12.28-equation}) imply that there is a constant $\delta_3>0$
such that %
\be  \label{12.29-equation}  %
{\mathcal{L}}^k(K_{m_s}\setminus \Omega^*(x'_0))\ge \delta_3 %
\ee  %
for all sufficiently large $s$. Now one can easily see that
(\ref{12.29-equation}) contradicts (\ref{12.26-equation}) and
(\ref{12.27-equation}) if $\delta_2>0$ in (\ref{12.26-equation})
is sufficiently small.

The proof of equation (\ref{3.16-equation}) of Theorem~\ref{Sarvas
approximation lemma} is finished.

\medskip

To prove (\ref{3.17-equation}) for $\Omega\in
{\mathcal{G}}_{n,b}$, we fix $x'\in \mathbb{R}^{n-k}$ such that
the slice $\Omega(x')$ is not empty. Considering the restrictions
of the symmetrizations $S$, $S_1$, and $S_2$ to the slice
$\mathbb{R}^{n}(x')$, we obtain %
\be  \label{12.34-equation}  %
\lim_{j\to \infty} d(\partial
\Omega_j(x'),\partial\Omega^*(x'))=0. %
\ee  %

Since $\Omega^*(x')$ is an open $k$ dimensional ball and
${\mathcal{L}}^k(\Omega^*(x'))={\mathcal{L}}^k(\Omega_j(x'))$,
equation (\ref{12.34-equation}) implies (\ref{3.17-equation}) in
the case of bounded open sets.

The proof of Theorem~\ref{Sarvas approximation lemma} is now
complete.   \hfill $\Box$

\medskip

\noindent%
\textbf{Remark 12.1.} %
One can easily show that (\ref{3.17-equation}) remains valid even
for unbounded sets $\Omega$ if the measure of the corresponding
slice is finite, i.e. if ${\mathcal{L}}^k(\Omega(x'))<\infty$.

In contrast, simple examples of unbounded domains $\Omega$ with a
finite measure, ${\mathcal{L}}^n(\Omega)<\infty$, show that
(\ref{3.16-equation}) is not true for unbounded open sets in
general.

%%%%%%%%%%%%%%%%%%%%%%%%%%%%%%%%%%%%%%%%%%%%%%%%%%%%%%%%%%%%%%%%%%%%%%%
%%%%%%%%%%%%%%%%%%%%%%%%%%%%%%%%%%%%%%%%%%%%%%%%%%%%%%%%%%%%%%%%%%%%%%%%
\end{document}